\pdfoutput=1
\documentclass{amsart}

\usepackage[T1]{fontenc}
\DeclareUnicodeCharacter{2032}{'}
\usepackage{dsfont} 
\usepackage{amsmath, amsthm, amsxtra, amsfonts, amssymb}
\usepackage{latexsym, ifthen, multicol, hyperref, mathdots, array}
\usepackage[dvipsnames]{xcolor}
\usepackage{bm} 
\usepackage{verbatim, comment}
\usepackage[tracking,final]{microtype} 
\usepackage[font={small, normalfont}, labelfont={}, figurename=Figure]{caption}
\usepackage{subcaption} 

\usepackage{mathtools, etoolbox}
\usepackage{stmaryrd} 
\usepackage{tikz-cd}
\usepackage{float} 
\usepackage{enumitem} 
	\setlist[enumerate,1]{label={\upshape(\arabic*)}}
\usepackage[export]{adjustbox} 
\usepackage[notcite, notref, final]{showkeys}

\usepackage[textsize = tiny, backgroundcolor = white, textcolor = orange, bordercolor = orange]{todonotes}
\usepackage{soul} 

\usepackage{import}
\usepackage{xifthen}
\usepackage{pdfpages}
\usepackage{transparent}

\theoremstyle{plain}
\numberwithin{equation}{section} 
\newtheorem{theorem}[equation]{Theorem}
\newtheorem{lemma}[equation]{Lemma}

\newtheorem{proposition}[equation]{Proposition}

\newtheorem{corollary}[equation]{Corollary}

\theoremstyle{definition}
\newtheorem{definition}[equation]{Definition}
\newtheorem{example}[equation]{Example}

\newtheorem{remark}[equation]{Remark}

\newcommand{\Cat}{\mathsf}
\newcommand{\var}[1]{{\operatorname{#1}}}

\makeatletter
\providecommand*{\twoheadrightarrowfill@}{%
	\arrowfill@\relbar\relbar\twoheadrightarrow
}
\providecommand*{\twoheadleftarrowfill@}{%
	\arrowfill@\twoheadleftarrow\relbar\relbar
}
\providecommand*{\xtwoheadrightarrow}[2][]{%
	\ext@arrow 0579\twoheadrightarrowfill@{#1}{#2}%
}
\providecommand*{\xtwoheadleftarrow}[2][]{%
	\ext@arrow 5097\twoheadleftarrowfill@{#1}{#2}%
}
\makeatother

\makeatletter
\newcommand*{\relrelbarsep}{.386ex}
\newcommand*{\relrelbar}{%
	\mathrel{%
		\mathpalette\@relrelbar\relrelbarsep
	}%
}
\newcommand*{\@relrelbar}[2]{%
	\raise#2\hbox to 0pt{$\m@th#1\relbar$\hss}%
	\lower#2\hbox{$\m@th#1\relbar$}%
}
\providecommand*{\rightrightarrowsfill@}{%
	\arrowfill@\relrelbar\relrelbar\rightrightarrows
}
\providecommand*{\leftleftarrowsfill@}{%
	\arrowfill@\leftleftarrows\relrelbar\relrelbar
}
\providecommand*{\xrightrightarrows}[2][]{%
	\ext@arrow 0359\rightrightarrowsfill@{#1}{#2}%
}
\providecommand*{\xleftleftarrows}[2][]{%
	\ext@arrow 3095\leftleftarrowsfill@{#1}{#2}%
}
\makeatother

\vbadness=\maxdimen

\def\ddefloop#1{\ifx\ddefloop#1\else\ddef{#1}\expandafter\ddefloop\fi}
\def\ddef#1{\expandafter\def\csname cat#1\endcsname{\ensuremath{\mathcal{#1}}}}
\ddefloop ABCDEFGHIJKLMNOPQRSTUVWXYZ{YD} \ddefloop
\def\ddef#1{\expandafter\def\csname ds#1\endcsname{\ensuremath{\mathds{#1}}}}
\ddefloop ABCDEFGHIJKLMNOPQRSTUVWXYZabcdefghijklmnopqrstuvwxyz \ddefloop
\def\ds1{\mathds{1}}
\newcommand{\mathsc}[1]{{\normalfont\textsc{#1}}}

\newcommand{\dual}{{}^\vee \mspace{-2mu}}
\newcommand{\odual}{\vee}

\newcommand{\End}{\mathsc{End}}
\newcommand{\Nat}{\mathsc{Nat}}
\newcommand{\Dinat}{\mathsc{Dinat}}

\author{Anh Tuong Nguyen}
\address{Department of Mathematics, University of Illinois Urbana-Champaign, Urbana, Illinois 61801, USA}
\email{anguye76@illinois.edu}
\title{Coend elements of a braided Hopf algebra}
\date{}

\begin{document}
\begin{abstract}
Let $H$ be a Hopf algebra in a braided rigid monoidal category $\mathcal{V}$ admitting a coend $C$. We define a ``coend element'' of $H$ to be a morphism from $C$ to $H$. We then study certain coend elements of $H$, which generalize important elements (e.g., pivotal and ribbon elements) of a finite dimensional Hopf algebra over a field. This builds on prior work of Bruguières and Virelizier \cite{BV-12:double-Hopf-mnd} on $R$-matrices of braided Hopf algebras. As an application, we provide another description for pivotal and ribbon structures on the category $\mathcal{V}_H$ of $H$-modules.
\end{abstract} 


\maketitle
\setcounter{tocdepth}{2}

\bibliographystyle{alpha}

\section{Introduction}

The purpose of this work is to generalize results for certain special (pivotal, ribbon, etc.) elements of a finite dimensional Hopf algebra over a field \( \dsk \) by replacing the category  \( \Cat{Vec}^{\var{fd}} \) of finite dimensional \( \dsk \)-vector spaces with an arbitrary braided rigid monoidal category \( \catV \) admitting a coend. 

To explain the motivation of this work, recall that a \emph{modular tensor category} is a nondegenerate braided finite tensor category equipped with a ribbon structure. We do not assume that such a category is semisimple, instead preferring to refer to semisimple modular tensor categories as \emph{modular fusion categories}. The importance of modular fusion categories arises from its connection to various fields including (extended) 3-dimensional topological quantum field theories and quantum invariants of 3-manifolds \cite{BDSPV-15:mtc-extended-tft, Turaev-16:invariants-knots-3-mflds}, conformal field theories and vertex operator algebras \cite{MS-89:CFT, Huang-05:VOA}, and topological quantum computing \cite{Wang-10:top-quantum-computation}. Recently, there have been increased interests in the non-semisimple setting, see for example \cite{KL-01:non-ss-TQFT, CGPM-14:invariants-3-mflds-non-ss-cat, DRGG-22:3d-tqft-non-ss-mtc}.

Given their importance, it is therefore desirable to produce constructions of modular fusion categories. There are many general constructions, some of which are modularization \cite{Bruguieres-00:modularization, Muger-00:galois-theory-BTC}, local modules \cite{Pareigis-95:braiding-dyslexia, KO-02:McKay-corresp}, the Drinfeld center of a spherical fusion category \cite{Muger-03:quantum-double-tensor-cat}; see also non-semisimple generalizations in \cite{Shimizu-19:non-deg-BFTC, LW-22:non-ss-MTC-local-mod, LW-22:non-ss-MTC-rel-center}. Classically, concrete examples of modular fusion categories arise from quantum groups at roots of unity \cite{BK-00:tensor-cat-modular-functors}; these quantum groups are examples of modular Hopf algebras as defined in \cite[Chapter IX.5]{Turaev-16:invariants-knots-3-mflds}.

As a general principle, Hopf algebras are useful in constructing modular categories because there is a one-to-one correspondence between additional structures on the category \( \catM_H \) of finite dimensional representations of \( H \), and special elements of \( H \). From this observation, \( R \)-matrices, ribbon and spherical elements are defined in \cite{Drinfeld-87:quantum-grps, RT-90:rbn-graphs, BW-99:spherical-cat} and correspond to braidings, ribbon structures, and spherical structures, respectively. In particular, \( \catM_H \) is modular if \( H \) is equipped with an \( R \)-matrix \( R \) and a ribbon element \( t \) such that \( (H, R) \) is factorizable, see e.g.\ \cite{Takeuchi-01:MTC-HA}. For instance, the Drinfeld double \( D(\dsk G) \) of the group algebra of a finite group \( G \) has a canonical \( R \)-matrix and a ribbon element which makes \( \catM_{D(\dsk G)} \) a modular tensor category \cite[Section~3.2]{BK-00:tensor-cat-modular-functors}.

\begin{table}[H]
{\small
\begin{center}
\begin{tabular}{c|c}
Structures or properties of \( \catM_H \) & Special elements or properties of \( H \) \\ 
\hline 
Braiding \( \sigma \) for \( \catM_H \) & \( R \)-matrix \( R \in H \otimes H \) \\
Pivotal structure \( \phi \) on \( \catM_H \) & Pivotal element \( p \in H \) \\
Ribbon structure (twist) \( \theta \) on \( (\catM_H, \sigma) \) & Ribbon element (twist) \( t \) for \( (H, R) \) \\
Nondegeneracy of \( (\catM_H, \sigma) \) & \( (H, R) \) is factorizable

\end{tabular} 
\caption{Classical bijective correspondence between structures or properties of \( \catM_H \) and special elements or properties of \( H \).}
\label{Tbl:TK-duality-Vec}
\end{center}}
\end{table}

In the effort to find more examples of modular tensor categories arising from Hopf algebras, one could try to replace \( \catM_H \) with \( \catV_H \), the category of modules over a Hopf algebra in a braided finite tensor category \( \catV \), and generalize the right hand side of Table~\ref{Tbl:TK-duality-Vec}.

In \cite{BV-12:double-Hopf-mnd}, Bruguières and Virelizier gave a definition of an \( R \)-matrix for a Hopf algebra \( H \) in a braided rigid monoidal category \( \catV \), provided that the category \( \catV \) admits a coend \( C \). Here, \( C \) is a certain object in \( \catV \) that is equipped with a universal dinatural transformation \( \rho_X: \dual X \otimes X \to C \), where \( \dual X \) is the left dual of \( X \), see Section~\ref{Sec:coend} for more details. An \( R \)-matrix for \( H \) is then defined as a morphism \( \mathfrak{R}: C \otimes C \to H \otimes H \) satisfying certain axioms (Definition~\ref{Def:C-qst-HA}). This generalizes the classical notion of \( R \)-matrices in \( \Cat{Vec}^\var{fd} \) as certain elements of \( H \otimes H \), since \( C = \dsk \) when \( \catV =  \Cat{Vec}^\var{fd} \). Under this definition, they recover the 1-1 correspondence between \( R \)-matrices of \( H \) and braidings on \( \catV_H \) (Theorem~\ref{Res:C-R-mat-1-1}).

By utilizing similar techniques as in \cite{BV-12:double-Hopf-mnd}, we can make sense of the notion of an ``element'' of a Hopf algebra \( H \) in \( \catV \) as follows.

\begin{definition}
	Let \( H \) be a Hopf algebra in a braided rigid monoidal category \( \catV \) admitting a coend \( C \). A \emph{coend element} (or \( C \)-\emph{element}) of \( H \) is a morphism \( h: C \to H \). 
\end{definition}

For example, when \( \catV = \Cat{Vec}^\var{fd} \), a \( \dsk \)-element of \( H \) is an element of \( H \) in the usual sense. We can now state the main results of this paper. Note that \( \catV_H \) denotes the category of right \( H \)-modules in \( \catV \). 

\begin{theorem}[Theorems~\ref{Res:C-piv-1-1}, \ref{Res:C-rbn-1-1}]
Let \( \catV \) be a braided rigid monoidal category admitting a coend \( C \), and let \( H \) be a Hopf algebra in \( \catV \). 
\begin{enumerate}[label=\upshape(\alph*)]
\item Consider \( C \)-\emph{pivotal elements} of \( H \) defined in \textup{Definition~\ref{Def:C-piv}}. There is a canonical bijection between the set \( \mathsf{CPiv}(H) \) of \( C \)-pivotal elements of \( H \) and the set \( \mathsf{Piv}(\catV_H) \) of pivotal structures of \( \catV_H \).
\item  Assume that \( H \) is quasitriangular with an \( R \)-matrix \( \mathfrak{R} \) as in \cite{BV-12:double-Hopf-mnd}, and consider \( C \)-\emph{balanced (resp. ribbon) elements} of \( H \) defined in \textup{ Definition~\ref{Def:C-rbn}}.
\begin{enumerate}[label=\upshape(\roman*)]
\item  There is a canonical bijection between the set \( \mathsf{CBal}(H) \) of \( C \)-balanced elements of \( H \) and the set \( \mathsf{Bal}(\catV_H) \) of balanced structures of \( \catV_H \).
\item The bijection in (i) further restricts to a bijection between the subset \( \mathsf{CRbn}(H) \) of \( C \)-ribbon elements of \( H \) and the subset \( \mathsf{Rbn}(\catV_H) \) of ribbon structures of \( \catV_H \).
\end{enumerate}
\end{enumerate}
\end{theorem}

Thus, we obtain a generalization to Table~\ref{Tbl:TK-duality-Vec} as shown in Table~\ref{Tbl:TK-duality-V}. In particular, we also describe a characterization of the nondegeneracy of \( \catV_H \), which is explained in Section~\ref{Sec:factorizable-HA}. 

\begin{table}[ht]
{\small
\begin{center}
\begin{tabular}{c|c}
Structures or properties of \( \catV_H \) & Special elements or properties of \( H \)  \\ 
\hline
Braiding \( \sigma_{M, N}: M \otimes N \to N \otimes M \) & \( R \)-matrix \( \mathfrak{R}: C \otimes C \to H \otimes H \) \cite{BV-12:double-Hopf-mnd} \\
Pivotal structure \( \phi_M: M \to M \dual \dual \) & Pivotal element \( \mathfrak{p}: C \to H \) \\
Ribbon structure \( \theta_M: M \to M \) & Ribbon element \( \mathfrak{t}: C \to H \)  \\
Nondegeneracy & \( (H, \mathfrak{R}) \) is factorizable
\end{tabular} 
\caption{Generalized bijective correspondence between structures or properties of \( \catV_H \) and special \( C \)-elements or properties of \( H \).}
\label{Tbl:TK-duality-V}
\end{center}}
\end{table}

Addtionally, there is a relation between these generalized pivotal and ribbon elements via the Drinfeld element \( \mathfrak{u} \) and two other elements \( \mathfrak{q}_\mu \) and \( \mathfrak{c}_\mu \) as seen in the following theorem, which generalizes well-known results over vector spaces, see e.g., \cite{Drabant-02:balanced-cat-HA} or Section~\ref{Sec:elts-Vec}.

\begin{theorem}[Theorems~\ref{Res:C-Drinfeld-bal-piv}, \ref{Res:C-Drinfeld-rbn-piv}]
Let \( \catV \) be a braided rigid monoidal category admitting a coend \( C \), and let \( (H, \mathfrak{R}) \) be a quasitriangular Hopf algebra in \( \catV \) as in \cite{BV-12:double-Hopf-mnd}.
\begin{enumerate}[label=\upshape(\alph*)]
\item Consider the \( C \)-Drinfeld element \( \mathfrak{u} \) as defined in \textup{Definition~\ref{Def:C-Drinfeld}}. This element induces a bijection \( \mathsf{CBal}(H) \cong \mathsf{CPiv}(H) \) via \( \mathfrak{t} \mapsto \mathfrak{t} \mathfrak{u} \). 
\item Under the above correspondence, a \( C \)-balanced element \( \mathfrak{t} \) is ribbon if and only if one of the following two equivalent conditions is satisfied:
\begin{enumerate}[label=\upshape(\roman*)]
\item  \( \mathfrak{t}^{-2} = \mathfrak{c}_\mu \),
\item the corresponding pivotal element \( \mathfrak{p} \) satisfies \( \mathfrak{p}^2 = \mathfrak{q}_\mu \),
\end{enumerate}
where \( \mathfrak{q}_\mu \) and \( \mathfrak{c}_\mu \) are certain \( C \)-elements defined in \textup{Definition~\ref{Def:C-q-c}}. The product on \( \var{Hom}(C, H) \) used in this result is defined in \textup{Figure~\ref{Fig:C-conv-prod}(b)}.
\end{enumerate}
\end{theorem}

To obtain the results above, there is an important intermediate gadget that provides the bridge between \( C \)-elements of \( H \) and structures on \( \catV_H \). These are elements of \( \Nat(1_\catV, 1_\catV \otimes H) \), i.e., natural transformations \( \{ \alpha_X: X \to X \otimes H \}_{X \in \catV} \). It is well-known that the endofunctor \( T_H = 1_\catV \otimes H \) is a monad on \( \catV \), such that the category of modules of \( T_H \) coincides with the category \( \catV_H \) of modules over \( H \). Let \( U_H: \catV_H \to \catV \) denote the forgetful functor from \( \catV_H \) to \( \catV \). The key observation in \cite{BV-07:Hopf-mnd} is that there is a bijection \( \End(U_H) \cong \Nat(1_\catV, T_H) \), which generalizes the bijection \( \End(U_H) \cong H, \alpha \mapsto \alpha_H(1_H) \) when \( H \) is a finite dimensional Hopf algebra over \( \dsk \).

In fact, the theory of ``monadic elements'' \( \alpha \in \Nat(1_\catV, T) \) makes sense for arbitrary Hopf monads \( T \) and is the subject of the paper \cite{BV-07:Hopf-mnd}, even for categories that do not necessarily admit a coend. However, for Hopf monads of the type \( T_H \) arising from a Hopf algebra \( H \) in \( \catV \), provided that the category \( \catV \) admits a coend \( C \), one can use graphical calculus to obtain a more direct generalization of elements of \( H \) as morphisms \( C \to H \). 

The paper is organized as follows. In Section 2, we introduce the general background material needed for this paper, including monoidal categories, braided Hopf algebras, and coends. In Section 3, we summarize the main results of \cite{BV-07:Hopf-mnd} concerning elements of a Hopf monad \( T \) in the special context that \(  T = 1_\catV \otimes H \) is the Hopf monad induced by a Hopf algebra in a braided rigid monoidal category \( \catV \). At the end of Section 3, we modify the theory slightly to relax the pivotal assumption on \( \catV \), which is necessary for arbitrary Hopf monads but not for Hopf monads induced by Hopf algebras. This leads us to Section 4, where we further assume that \( \catV \) admits a coend \( C \). Using graphical calculus and universal properties of the coend \( C \), we display the \( C \)-element versions of special elements such as pivotal and ribbon elements of \( H \) in this setting, and prove their generalized properties and relations. Finally, we give some concluding remarks including a few applications of our work in Section 5. 

\section{Preliminaries}

The goal of this section is to introduce the general background material, including basic definitions and results, as well as notations and conventions that will be used throughout this paper. In Section~\ref{Sec:monoidal-cat}, we define the setting where all of our objects live, which are monoidal categories endowed with various additional structures. In Section~\ref{Sec:braided-HA-HM}, we introduce the first main object, that of a braided Hopf algebra. We also discuss its generalization in the form of a Hopf monad on a rigid monoidal category. In Section~\ref{Sec:coend}, we introduce the second main object: a certain universal object called a coend \( C \) of a rigid monoidal category. Next, we briefly discuss all classical results about special elements of a Hopf algebra over a field \( \dsk \) in Section~\ref{Sec:elts-Vec}, which form the basis of our work. Finally, in Section~\ref{Sec:C-qst-HA}, we introduce the first result in this direction that forms the inspiration of this work: a generalization of the \( R \)-matrix obtained by Bruguières and Virelizier in \cite{BV-12:double-Hopf-mnd}.

\subsection{Monoidal categories} \label{Sec:monoidal-cat}

We review some general facts about monoidal categories, which will be used extensively throughout. 

\subsubsection{Conventions for categories}

We assume that the reader is familiar with standard category concepts as presented in \cite{MacLane-98:cat, Riehl-16:cat-theory}. Unless otherwise specified, all categories are small.

Let \( \catC \) and \( \catD \) be categories. The category of functors and natural transformations from \( \catC \) to \( \catD \) is denoted by \( [\catC, \catD] \). Given two functors \( F, G: \catC \to \catD \), the set of natural transformations from \( F \) to \( G \) is denoted by \( \Nat(F, G) \) and the monoid of natural endomorphisms of \( F \) is denoted by \( \End(F) \). Given a category \( \catC \), the opposite category of \( \catC \) with all morphisms reversed is denoted by \( \catC^\var{op} \). 

\subsubsection{Monoidal categories}
For a thorough introduction to monoidal categories and graphical calculus, the reader can refer to \cite{TV-17:monoidal-cat-TFT}, \cite{EGNO-16:tensor-cat}, and \cite{Selinger-11:survey-graphical}. On account of MacLane's coherence theorem, we assume that all monoidal categories are strict. Given a monoidal category \( (\catC, \otimes, \ds1) \), the category \( \catC^{\var{op}} =  (\catC^\var{op}, \otimes, \ds1) \) is also monoidal. We also have a monoidal category  \( \catC^{\otimes \var{op}} = (\catC, \otimes^\var{op}, \ds1) \) where \( \catC \) is equipped with the reversed tensor product \( X \otimes^\var{op} Y = Y \otimes X \) for all \( X, Y \in \catC \). 

A \emph{monoidal functor} between monoidal categories \( (\catC, \otimes, \ds1)  \) and \( (\catD, \underline{\otimes}, \underline{\ds1}) \) is a triple \( (F, F^{(2)}, F^{(0)}) \), where \( F: \catC \to \catD \) is a functor, \( F^{(2)}: F(-) \underline{\otimes} F(-) \to F(- \otimes -) \) is a natural transformation, and \( F^{(0)}: \underline{\ds1} \to F(\ds1)  \) is a distinguished morphism, satisfying coherence axioms. A monoidal functor is \emph{strong} if \( F^{(2)} \) and \( F^{(0)} \) are isomorphisms. A \emph{comonoidal functor} is a monoidal functor \( C^{\var{op}} \to D^{\var{op}} \). 

A \emph{monoidal natural transformation} between two monoidal functors \( F, G: \catC \to \catD \) is a natural transformation \( \eta: F \to G \) such that \( \eta_{X \otimes Y} F^{(2)}_{X, Y} = G^{(2)}_{X, Y} (\eta_X \otimes \eta_Y) \) for all \(  X, Y \in \catC \), and \( \eta_\ds1 F^{(0)} = G^{(0)} \). The set of monoidal natural transformations between two monoidal functors \( F \) and \( G \) is denoted by \( \Nat_{\otimes}(F, G) \), and the set of monoidal natural endomorphisms of a monoidal functor \( F \) is denoted by \( \End_{\otimes}(F) \). 

\subsubsection{Rigidity}

Let \( \catC \) be a monoidal category. A \emph{duality pairing} in \( \catC \) is a quadruple \( (X, Y, e, c) \) consisting of objects \( X, Y \) of \( \catC \) and morphisms \( e: X \otimes Y \to \ds1 \) and \( c: \ds1 \to Y \otimes X \), such that the two compositions
\begin{equation*}
	X \xrightarrow{X \otimes c} X \otimes Y \otimes X \xrightarrow{e \otimes X} X, \quad
	Y \xrightarrow{c \otimes Y} Y \otimes X \otimes Y \xrightarrow{Y \otimes e} Y
\end{equation*}
are both identity morphisms. In this case, \( X \) is called a \emph{left dual} of \( Y \), \( Y \) a \emph{right dual} of \( X \), and \( e \) and \( c \) are the \emph{evaluation} and \emph{coevaluation} maps respectively. A monoidal category \( \catC \) is \emph{rigid} if every object \( X \in \catC \) is part of distinguished duality pairings \( (\dual X, X, \var{ev}^l_X, \var{coev}^l_X) \) and \( (X, X \dual, \var{ev}^r_X, \var{coev}^r_X) \). The objects  \( \dual X \) and \( X \dual \) are called the \emph{left dual} and the \emph{right dual} of \( X \), respectively, and \( \var{ev}^l_X \), \( \var{coev}^l_X \) (resp. \( \var{ev}^r_X \), \( \var{coev}^r_X \)) are called the left (resp. right) evaluation and coevaluation maps. 

It is well-known that there are unique isomorphisms between any two left (or right) duals of an object \( X \) respecting the evaluation and coevaluation maps. In particular, we will abstain from writing the following canonical isomorphisms
\begin{equation}
\label{Eqn:monoidal-identification}
\begin{split}
	(\dual X) \dual \cong X \cong \dual (X \dual) &, \quad \dual \ds1 \cong \ds 1 \cong \ds1 \dual, \\
	\dual (X \otimes Y) \cong \dual Y \otimes \dual X &, \quad (X \otimes Y) \dual \cong Y \dual \otimes X \dual.
\end{split}
\end{equation}

\begin{figure}[ht]
    \centering
    \scalebox{0.8}{%
    \def\svgwidth{\columnwidth}
    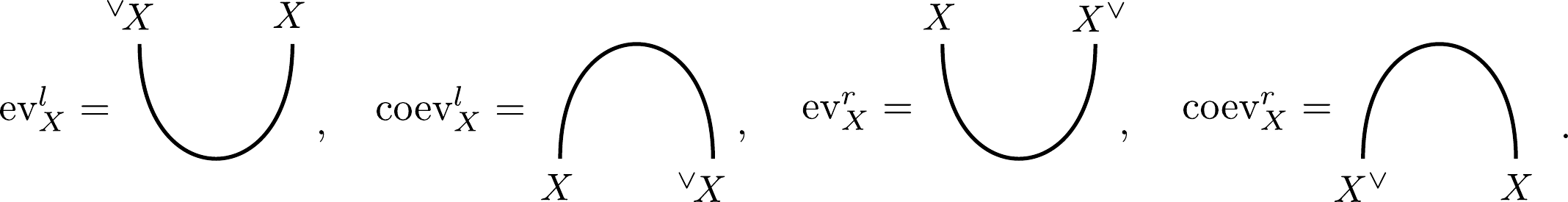
}
	\caption{Left and right evaluation and coevaluation maps of \( X \in \catC \).}
	\label{Fig:(co)ev}
\end{figure}

We follow the top-to-bottom convention for the graphical calculus of monoidal categories. For any \( X \in \catC \), the left and right evaluation and coevaluation maps are depicted as in Figure~\ref{Fig:(co)ev}.

\begin{figure}[ht]
    \centering
    \scalebox{.5}{%
    \def\svgwidth{\columnwidth}
    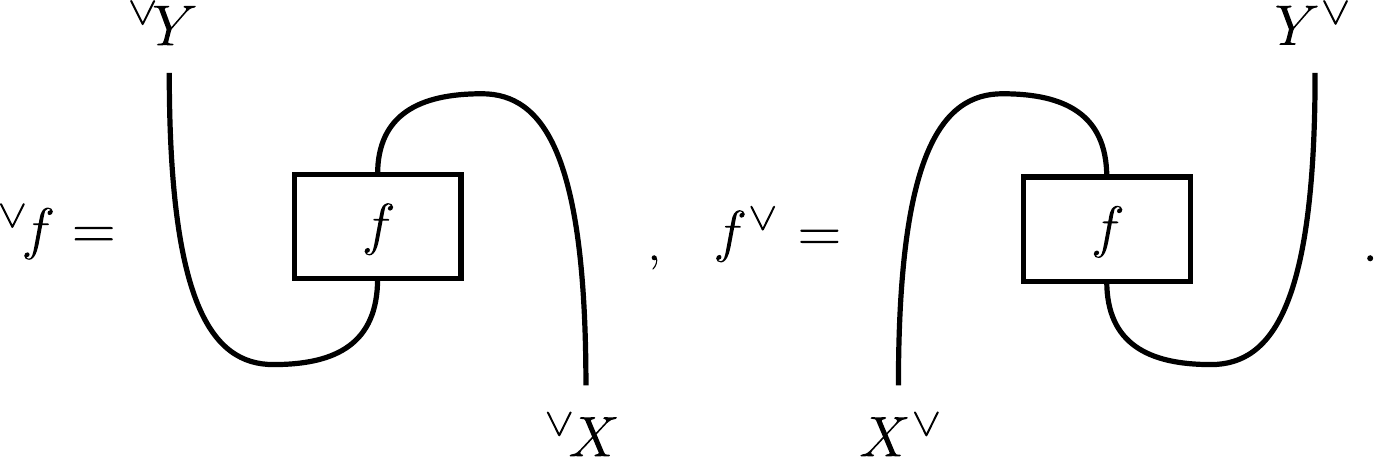
}
	\caption{Left and right dual of a morphism \( f: X \to Y \).}
	\label{Fig:dual-morphism}
\end{figure}

Furthermore, for any morphism \( f: X \to Y \) we define \( \dual f \) and \( f \dual \) as in Figure~\ref{Fig:dual-morphism}, such that \( (-) \dual \) and \( \dual (-) \) are strong monoidal functors \( \catC^{\var{op}} \to \catC^{\otimes \var{op}} \). 

\subsubsection{Conjugation by duality functors}
\label{Sec:conj-duality}

Let \( \catC, \catD \) be rigid monoidal categories. For any functor \( F: \catC \to \catD \), define \( F^!: \catC \to \catD \) by \( F^!(X) = \dual F( X \dual)  \)  and \( F^!(f) = \dual F(f \dual)  \). For any natural transformation \( \alpha: F \to G \), where \( F, G \) are functors \( \catC \to \catD \), define \( \alpha^!_X = {}^{\odual} \! \alpha_{X \dual}: G^!(X) \to F^!(X) \). This defines a functor 
\begin{equation*}
	(-)^!: [\catC, \catD] \to [\catC, \catD]^{\var{op}}, \qquad F \mapsto F^!, \quad \alpha \mapsto \alpha^!,
\end{equation*}
which we call the \emph{conjugation by duality functor}. Similarly, one can define a functor \( {}^! (-): [\catC, \catD]^{\var{op}} \to [\catC, \catD] \) which is quasi-inverse to \( (-)^! \). Furthermore, for composable functors \( F \) and \( G \) between rigid monoidal categories, we have \( (G \circ F)^! = G^! \circ F^! \), and for  composable natural transformations \( \alpha \) and \( \beta \) between functors between rigid monoidal categories, we have \( (\beta \circ \alpha)^! = \beta^! \circ \alpha^! \). 

\begin{remark}
\label{Rem:antipode-End(F)}
It follows that for any endofunctor \( F \) on a rigid monoidal category, the conjugation by duality functor \( (-)^! \) is an antiautomorphism of \( \End(F) \) with inverse \( {}^! (-) \).
\end{remark}

\begin{lemma}[{\cite[Lemma 3.4]{BV-07:Hopf-mnd}}]
\label{Res:mon-trans-rigid}
Let \( F, G: \catC \to \catD \) be strong monoidal functors and let \( \alpha: F \to G \) be a monoidal natural transformation. If \( \catC \) is rigid, then \( \alpha \) is an isomorphism with \( \alpha^{-1} = \alpha^! = {}^! \alpha \). \qed
\end{lemma}

\subsubsection{Pivotal structures}

A rigid monoidal category \( \catC \) is \emph{pivotal} if it is equipped with a \emph{pivotal structure},  i.e., a monoidal natural transformation \( \phi_X: X \to X \dual \dual \) such that \( \phi_{X \otimes Y} = \phi_X \otimes \phi_Y \) and \( \phi_\ds1 = \var{id}_\ds1 \), up to the identifications~(\ref{Eqn:monoidal-identification}). The set of pivotal structures of \( \catC \) is denoted by \( \mathsf{Piv}(\catC) \). 

\begin{remark}
\label{Rem:piv-invertible}
By Lemma~\ref{Res:mon-trans-rigid}, any pivotal structure \( \phi \) on a rigid monoidal category \( \catC \) is invertible with \( \phi^{-1}_X = {}^{\odual} \! \phi_{X \mspace{-1mu} \dual} = \phi_{\dual \mspace{-1mu} X}^{\odual} \) for all \( X \in \catC \). In particular, \( \phi^{!!} = \phi \). The invertibility of \( \phi \) also implies that \( \phi_\ds1 = \var{id}_\ds1 \) if we identify \( \ds1 = \ds1 \dual \dual \), since \( \phi_\ds1^{\otimes 2} = \phi_\ds1 \). 
\end{remark}

\subsubsection{Braided categories}

Let \( \catC \) be a monoidal category. A \emph{lax braiding} on \( \catC \) is a natural transformation \( \sigma_{X, Y}: X \otimes Y \to Y \otimes X \) satisfying
\begin{equation*}
	\sigma_{X, Y \otimes Z} = (\var{id}_Y \otimes \sigma_{X, Z})(\sigma_{X, Y} \otimes \var{id}_Z), \quad 
	 \sigma_{X \otimes Y, Z} = (\sigma_{X, Z} \otimes \var{id}_Y)(\var{id}_X \otimes \sigma_{Y, Z}),
\end{equation*}
and \( \sigma_{X, \ds1} = \var{id}_X = \sigma_{\ds1, X} \), for all \( X, Y, Z \in \catC \). A \emph{braiding} is a lax braiding that is also invertible. A monoidal category \( \catC \) is \emph{braided} if it is equipped with a braiding. 

\begin{figure}[ht]
    \centering
    \begin{subfigure}{0.5\textwidth}
        \centering
        \scalebox{0.8}{%
    \def\svgwidth{\columnwidth}
    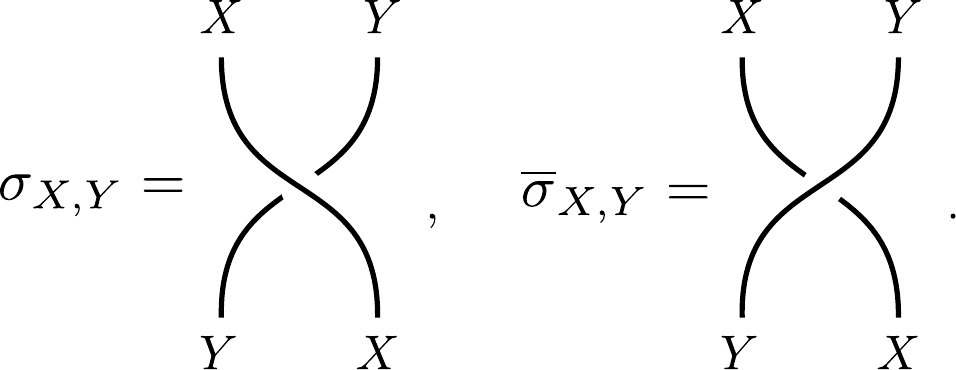
}
        \caption{}
    \end{subfigure}\hfill
    \begin{subfigure}{0.5\textwidth}
        \centering
        \scalebox{0.9}{%
    \def\svgwidth{\columnwidth}
    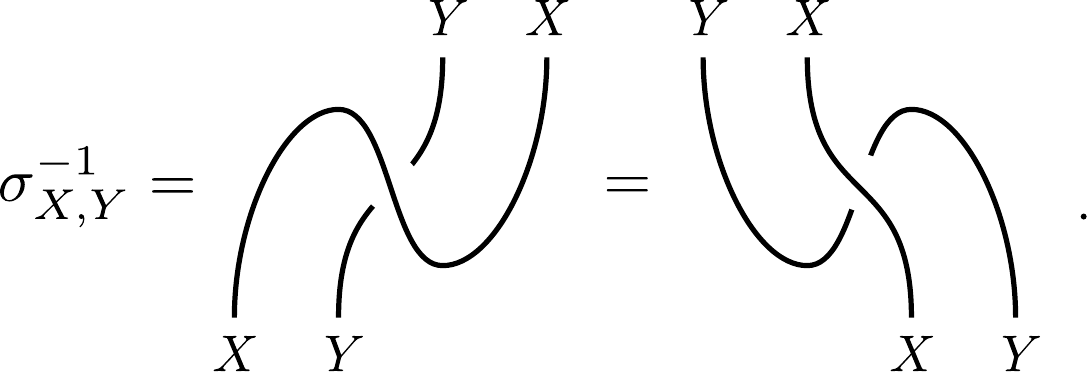
}
        \caption{}
    \end{subfigure}
    \caption{Graphical calculus of a braiding on a (rigid) monoidal category.}
    \label{Fig:braiding}
\end{figure}

If \( \catC = (\catC, \sigma) \) is a braided monoidal category, then the reverse braiding \( \bar{\sigma}_{X, Y} = \sigma^{-1}_{Y, X} \) is also a braiding on \( \catC \). Graphically, the braiding \( \sigma \) and its reverse \( \bar{\sigma} \) are represented in Figure~\ref{Fig:braiding}(a).

\begin{remark}
\label{Rem:braiding-rev}
When \( \catC \) is rigid, every lax braiding \( \sigma \) on \( \catC \) is invertible, and therefore also a braiding, with \( \sigma^{-1}_{X, Y} \) given as in Figure~\ref{Fig:braiding}(b). Furthermore, the invertibility of \( \sigma \) readily implies that \( \sigma_{X, \ds1} = \var{id}_X = \sigma_{\ds1, X} \) for all \( X \in \catC \).
\end{remark}

\subsubsection{Ribbon categories}

A \emph{balanced category} is a braided monoidal category \( (\catC, \sigma) \) equipped with a \emph{balanced structure}, which is defined as a natural endomorphism \( \theta \in \End(1_\catC) \) satisfying \( \theta_{X \otimes Y} =  (\theta_X \otimes \theta_Y) \, \sigma_{Y, X} \, \sigma_{X, Y} \) for all \( X, Y \in \catC \), and \( \theta_\ds1 = \var{id}_\ds1 \). A balanced structure on \( (\catC, \sigma) \) is sometimes also called a \emph{twist}. If \( \catC \) is furthermore rigid, then the twist \( \theta \) is said to be \emph{self-dual} if it satisfies \( \theta = \theta^! \). In this case, we say that \( \theta \) is a \emph{ribbon structure} and \( \catC \) is a \emph{ribbon category}. The set of balanced (resp. ribbon) structures on \( \catC \) will be denoted by \( \mathsf{Bal}(\catC) \) (resp. \( \mathsf{Rbn}(\catC) \)).

\subsubsection{Drinfeld morphisms}

\label{Sec:Drinfeld-maps}

In this section, we fix a braided rigid monoidal category \( (\catC, \sigma) \). Using graphical calculus, the results in this section are straightforward to verify; see also \cite[Appendix A.2]{HPT-16:cat-trace} for similar results on balanced and pivotal structures of \( \catC \).

\begin{figure}[ht]
    \centering
     \scalebox{.6}{%
    \def\svgwidth{\columnwidth}
    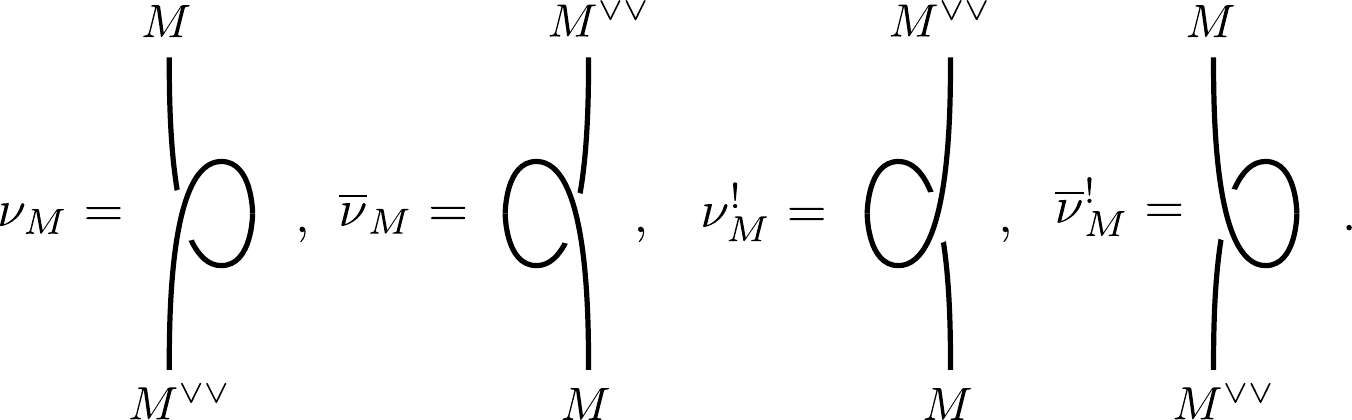
}
	\caption{Drinfeld morphisms in a braided rigid monoidal category.}
	\label{Fig:Drinfeld-maps}
\end{figure}

\begin{definition}
Define natural transformations \( \nu, \bar{\nu}, \nu^!, \bar{\nu}^! \) as in Figure~\ref{Fig:Drinfeld-maps}. We refer to these morphisms collectively as \emph{Drinfeld morphisms}, with \( \nu \) called the (\emph{right}) \emph{Drinfeld morphism}.
\end{definition}

\begin{lemma}
\label{Res:Drinfeld-map-properties}
The morphism \( \nu \) is invertible with inverse \( \nu^{-1} = \bar{\nu} \). Moreover, \( \nu \) satisfies the relation
\begin{equation}
\label{Eqn:balancing-relation}
	\nu_{M \otimes N} = (\nu_M \otimes \nu_N) \, \bar{\sigma}_{N,M} \, \bar{\sigma}_{M, N}
\end{equation}
for all \( M, N \in \catC \). The morphism \( \nu^! \) is the conjugation of \( \nu \) by duality functors, with inverse \( \bar{\nu}^! \), such that \( \nu^{!!}  = \nu \). \qed
\end{lemma}

\begin{lemma}
\label{Res:Drinfeld-bij}
The Drinfeld morphism \( \nu \) induces 
a bijection
\begin{equation}
\label{Eqn:Drinfeld-bij}
	 \mathsf{Bal}(\catC) \cong \mathsf{Piv}(\catC), \quad
	 \theta \mapsto \nu \circ \theta
\end{equation}
between balanced and pivotal structures of \( \catC \). \qed 
\end{lemma}

\begin{remark}
\label{Rem:twist-invertible}
Since any pivotal structure \( \phi \) is invertible with \( \phi^{-1} = \phi^! \), and the Drinfeld isomorphism \( \nu \) is invertible, it follows that any balanced/ribbon structure \( \theta \) on a rigid monoidal category is invertible as well: if \( \theta = \bar{\nu} \phi \) then \( \theta^{-1} = \phi^{-1} \nu \). Furthermore, in this case the condition \( \theta_\ds1 = \var{id}_\ds1 \) is automatic, since \( \theta_\ds1^{\otimes 2} = \theta_\ds1 \). 
\end{remark}

\begin{figure}[ht]
    \centering
    \scalebox{.55}{%
    \def\svgwidth{\columnwidth}
    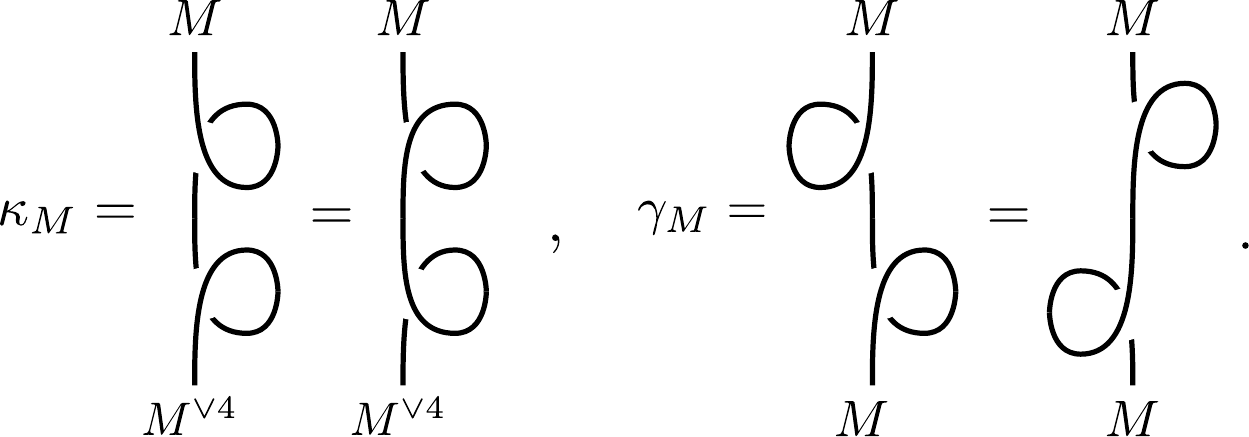
}
	\caption{The morphisms \( \kappa \) and \( \gamma \).}
	\label{Fig:Drinfeld-maps-2}
\end{figure}

\begin{definition}
\label{Def:Drinfeld-maps-2}
Define two natural transformations \( \kappa \) and \( \gamma \) by
\begin{equation*}
	\kappa_M =  \nu_{M \dual \dual} \bar{\nu}^!_M = \bar{\nu}^!_{M \dual \dual} \nu_M, \qquad \gamma_M = (\nu \nu^!)_{\dual \dual M}= (\nu^! \nu)_M ,
\end{equation*}
graphically displayed in Figure~\ref{Fig:Drinfeld-maps-2}. The equalities follow from the fact that the Drinfeld morphisms satisfy \( (-)^{!!} = \var{id} \), see Lemma~\ref{Res:Drinfeld-map-properties}.

\end{definition}

\begin{lemma}
\label{Res:Drinfeld-rbn-piv}
Let \( \catC \) be a braided rigid monoidal category, and let \( \kappa \) and \( \gamma \) be defined as in \textup{Figure~\ref{Fig:Drinfeld-maps-2}}. Under the correspondence between pivotal and balanced structures \textup{(\ref{Eqn:Drinfeld-bij})}, a pivotal structure \( \phi \) corresponds to a ribbon structure \( \theta \) if and only if one of the following equivalent conditions hold:
\begin{enumerate}[label=\upshape(\roman*)]
\item \( \phi^2 = \kappa \).
\item \( \theta^{-2} = \gamma \). \qed  
\end{enumerate}
\end{lemma}

\subsection{Braided Hopf algebras and Hopf monads}

 \label{Sec:braided-HA-HM}

\subsubsection{Braided Hopf algebras}

\begin{figure}[ht]
    \centering
    \begin{subfigure}[b]{0.75\textwidth}
        \centering
        \scalebox{1}{%
    \def\svgwidth{\columnwidth}
    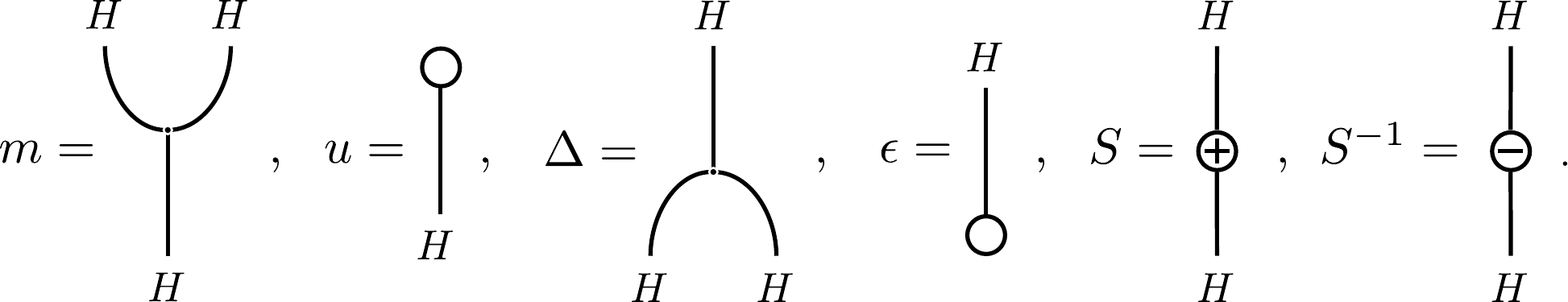
}
        \caption{}
    \end{subfigure}\hfill
    \begin{subfigure}[b]{0.2\textwidth}
        \centering
        \scalebox{0.7}{%
    \def\svgwidth{\columnwidth}
    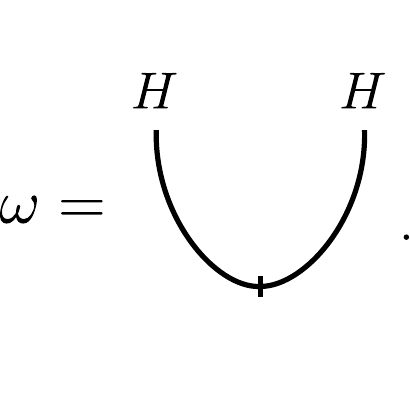
}
        \caption{}
    \end{subfigure}
    \caption{Graphical calculus for structure morphisms of \( H \).}
    \label{Fig:H-structure-maps}
\end{figure}

For the definitions of algebras, coalgebras, bialgebras and Hopf algebras in a braided monoidal category \( \catV \), see \cite[Chapter 6]{TV-17:monoidal-cat-TFT}. For a Hopf algebra \( H \), we represent its product \( m \), unit \( u \), coproduct \( \Delta \), counit \( \epsilon \), antipode \( S \) and its inverse \( S^{-1} \) graphically as in Figure~\ref{Fig:H-structure-maps}(a).

\begin{figure}[ht]
    \centering
    \scalebox{1}{%
    \def\svgwidth{\columnwidth}
    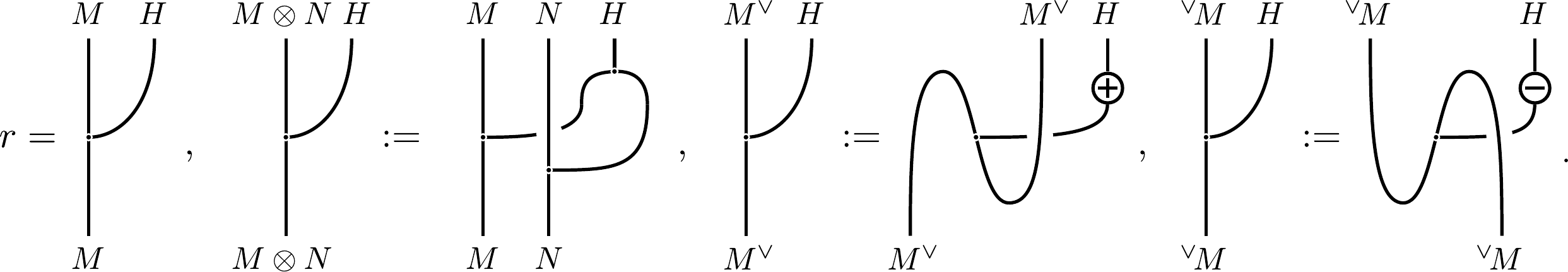
}
	\caption{The monoidal structure of \( \catV_H \).}
	\label{Fig:H-action}
\end{figure}

A right \( H \)-module is a pair \( (M, r) \) where \( M \in \catV \) and \( r: M \otimes H \to M \) is a morphism in \( \catV \) satisfying the usual associativity and unital axioms. It is well-known that the category \( \catV_H \) of right \( H \)-modules is a rigid monoidal category; the \( H \)-action \( r \), the monoidal structure and the dual actions are displayed in Figure~\ref{Fig:H-action}.

\begin{definition}[{\cite[Section 6.2.3]{TV-17:monoidal-cat-TFT}}]
	Let \( H = (H, m, u, \Delta, \epsilon) \) be a bialgebra in a braided monoidal category \( \catV \). A \emph{bialgebra pairing} for \( H \) is a morphism \( \omega: H \otimes H \to \ds1 \) in \( \catV \) such that 

\begin{equation*}
\begin{alignedat}{2}
	\omega (m \otimes \var{id}_H)  &= \omega(\var{id}_H \otimes \omega \otimes \var{id}_H)(\var{id}_{H \otimes H} \otimes \Delta), &\qquad \omega(u \otimes \var{id}_H) &= \epsilon, \\
	\omega (\var{id}_H \otimes m)  &= \omega(\var{id}_H \otimes \omega \otimes \var{id}_H)(\Delta \otimes \var{id}_{H \otimes H}), &\qquad \omega(\var{id}_H \otimes u) &= \epsilon.
\end{alignedat}
\end{equation*}
We depict a bialgebra pairing \( \omega \) graphically as in Figure~\ref{Fig:H-structure-maps}(b). Bialgebra pairings for Hopf algebras are called \emph{Hopf pairings}.

A bialgebra pairing \( \omega \) is \emph{non-degenerate} if there exists a \emph{copairing}, i.e., a morphism \( \Omega: \ds1 \to H \otimes H \) satisfying the snake equations
\begin{equation*}
	(\var{id}_H \otimes \omega)(\Omega \otimes \var{id}_H) = \var{id}_H = (\var{id}_H \otimes \Omega)(\omega \otimes \var{id}_H).
\end{equation*}
Equivalently, a bialgebra pairing \( \omega \) is a morphism whose left and right adjoint maps \( H \to H \dual \) and \( H \to \dual H \) are algebra homomorphisms, and it is nondegenerate if these maps are also isomorphisms.
\end{definition}

\subsubsection{Hopf monads}
Let \( \catC \) be a category. A \emph{monad} on \( \catC \) is an algebra \( (T, \mu, \eta) \) in the strict monoidal category \( (\var{End}(\catC), \circ, 1_\catC) \) of endofunctors of \( \catC \). Given a monad \( (T, \mu, \eta) \) on \( \catC \), a \( T \)-\emph{module} is a pair \( (M, r) \) consisting of an object \( M \in \catC \) and a morphism \( r: T(M) \to M \) in \( \catC \) such that \( 	r \, T(r) = r \, \mu_M \)  and \( r \, \eta_M = \var{id}_M \). Given two \( T \)-modules \( (M, r) \) and \( (N, s) \), a \( T \)\emph{-module morphism} \( f: (M,r) \to (N, s) \) is a morphism  \( f: M \to N \) of underlying objects in \( \catC \) such that \( fr = s \, T(f) \). The collection of \( T \)-modules and \( T \)-module morphisms forms a category, denoted by \( \catC_T \).

A monad \( (T, \mu, \eta) \) on a monoidal category \( \catC \) is a \emph{bimonad} if \( T = (T, T^{(2)}, T^{(0)}) \) is a comonoidal functor, and \( \mu \) and \( \eta \) are comonoidal natural transformations. If \( \catC \) is a rigid monoidal category and \( T \) is a bimonad, we say that \( T \) is a \emph{Hopf monad} if it is equipped with two natural transformations \( s^l_X: T({}^\vee T(X)) \to {}^\vee X \) and \( s^r_X: T(T(X)^\vee) \to X^\vee \) for \( X \in \catC \),  satisfying \cite[Equations (20)-(23)]{BV-07:Hopf-mnd}. They are called the \emph{left antipode} and the \emph{right antipode} for \( T \), respectively.

\begin{figure}[ht]
    \centering
    \scalebox{.6}{%
    \def\svgwidth{\columnwidth}
    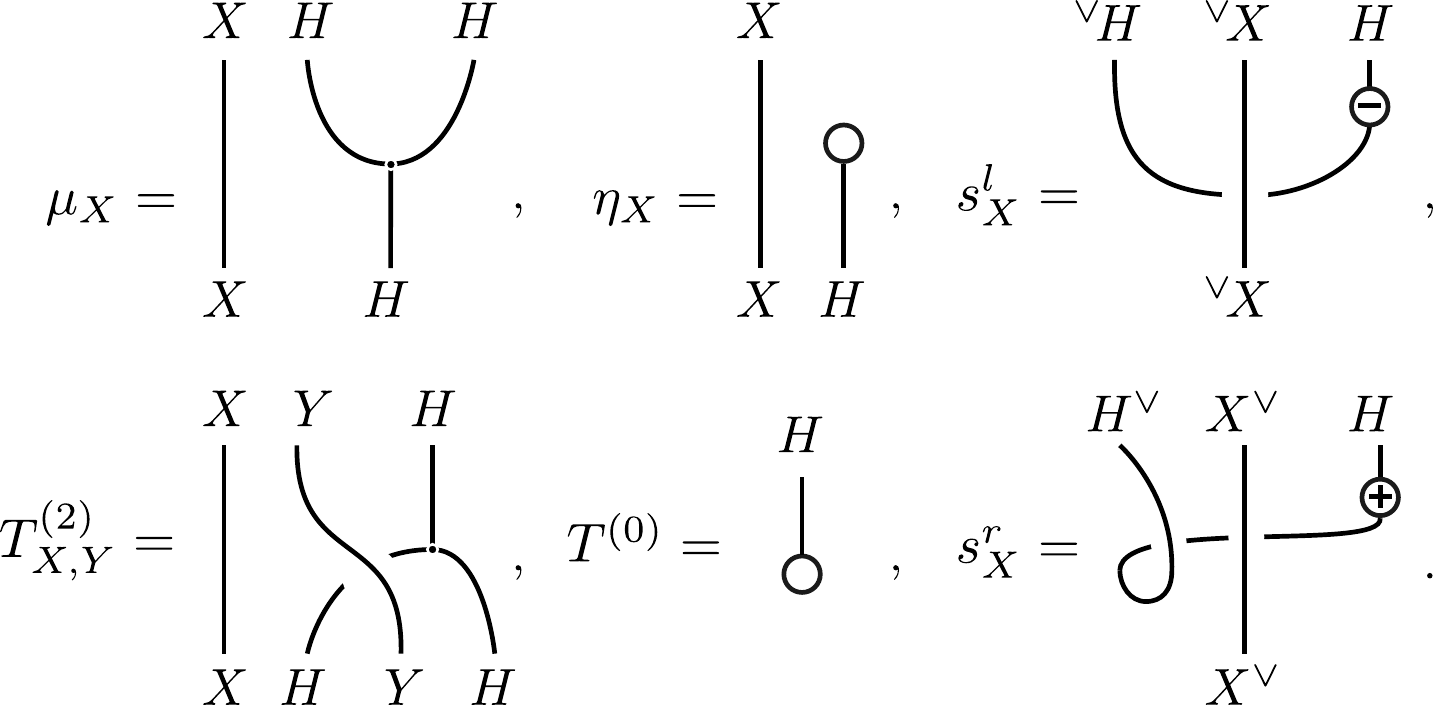
}
	\caption{Structure morphisms of the Hopf monad \( T = 1_\catV \otimes H \).}
	\label{Fig:ind-Hopf-mnd}
\end{figure}

\begin{example}
\label{Eg:ind-Hopf-mnd}
Let \( \catV \) be a braided rigid monoidal category and let \( H \) be a Hopf algebra in \( \catV \). The endofunctor \( T = 1_\catV \otimes H \) is a Hopf monad with structure morphisms presented in Figure~\ref{Fig:ind-Hopf-mnd} \cite[Example 2.4]{BV-12:double-Hopf-mnd}. Moreover, the category of \( T \)-modules coincides with the category of right \( H \)-modules in \( \catV \). 
\end{example}

\subsection{Coends}
\label{Sec:coend}

The references for this section are \cite[Sections IX.5 - IX.8]{MacLane-98:cat} and \cite[Sections 6.4 - 6.6]{TV-17:monoidal-cat-TFT}.

Let \( \catC, \catD \) be categories and let \( F: \catC^{\var{op}} \times \catC \to D \) be a functor. For an object \( D \) in \( \catD \), a \emph{dinatural transformation} from \( F \) to \( D \) is a family \( d = \{ d_X: F(X, X) \to D \}_{X \in \catC} \) of morphisms in \( \catD \) such that for every morphism \( f: X \to Y \) in \( \catC \), the diagram
\begin{equation*}
\begin{tikzcd}
	F(Y, X) & F(X, X) \\
	F(Y, Y) & D
	\arrow[from=1-1, to=1-2, "{F(f, \var{id}_X)}"]
	\arrow[from=1-1, to=2-1, "{F(\var{id}_Y, f)}"']
	\arrow[from=1-2, to=2-2, "d_X"]
	\arrow[from=2-1, to=2-2, "d_Y"']
\end{tikzcd}
\end{equation*}
commutes. Let \( \Dinat(F, D) \) denote the set of dinatural transformations from \( F \) to \( D \). 

A \emph{coend} of a functor \( F: \catC^{\var{op}} \times \catC \to \catD \) is a pair \( (C, \iota) \) where \( C \) is an object of \( \catD \) and \( \iota \in \Dinat(F, C) \) such that every dinatural transformation \( d \) from \( F \) to an object \( D \) of \( \catD \) factors as \( d = f \circ \iota \)  for a unique morphism \( f: C \to D \). In other words, we have a bijection
\begin{equation}
\label{Eqn:coend-dinat}
	\var{Hom}_\catD(C, D) \to \Dinat(F, D), \qquad f \mapsto f \circ \iota.
\end{equation}

\subsubsection{The coend of a rigid monoidal category} 

Let \( \catV \) be a rigid monoidal category, and let \( F: \catV^{\var{op}} \times \catV \to \catV \) denote the functor defined on objects by \( (X, Y) \mapsto \dual X \otimes Y \). We say that \( \catV \) \emph{admits a coend} if a coend of the functor \( F \) exists. We often denote this coend by \( (C, \iota) \) or simply \( C \), with \( C \) an object of \( \catV \), and \( \iota \) the dinatural transformation with \( X \)-component \( \iota_X: \dual X \otimes X \to C \) for \( X \in \catV \). 

Using the snake identity, the proof of the following lemma is straightforward.

\begin{lemma}
Let \( \catV \) be a rigid monoidal category, and let \( F \) be the functor defined as above. For any object \( D \in \catV \), there is a bijection
\begin{equation}
\label{Eqn:coend-bij}
	\Dinat(F, D) \to \Nat(1_\catV, 1_\catV \otimes D), \quad
	\phi \mapsto \tilde{\phi}_X := (\var{id}_X \otimes \phi_X)(\var{coev}^l_X \otimes \var{id}_X) 
\end{equation}
with inverse \( \psi \mapsto \bar{\psi}_X = (\var{ev}^l_X \otimes \var{id}_X)(\var{id}_X \otimes \psi_X) \). \qed 
\end{lemma}

By using the above lemma, in particular by composing the two bijections (\ref{Eqn:coend-dinat}) and (\ref{Eqn:coend-bij}), we see that equivalently, if \( (C, \iota) \) is a coend of \( \catV \) and \( \delta = \tilde{\iota} \), then for any \( D \in \catV \) we have a bijection
\begin{equation}
	\label{Eqn:coend-fp}
	\Sigma^{(1)}_D: \var{Hom}(C, D) \to \Nat(1_\catV, 1_\catV \otimes D), \qquad f \mapsto (\var{id}_X \otimes f) \, \delta_X.
\end{equation}
The natural transformation \( \delta_X: X \to X \otimes C \) is called the \emph{universal coaction} associated to the coend \( C \). Even though \( \iota \) and \( \delta \) carry the same data, working with \( \delta \) is often slightly more convenient. In terms of graphical calculus, we color \( C \)-strands gray, with the universal coaction \( \delta_X: X \to X \otimes C \) graphically depicted as
\begin{equation*}
	\scalebox{0.15}{%
    \def\svgwidth{\columnwidth}
\begingroup%
  \makeatletter%
  \providecommand\color[2][]{%
    \errmessage{(Inkscape) Color is used for the text in Inkscape, but the package 'color.sty' is not loaded}%
    \renewcommand\color[2][]{}%
  }%
  \providecommand\transparent[1]{%
    \errmessage{(Inkscape) Transparency is used (non-zero) for the text in Inkscape, but the package 'transparent.sty' is not loaded}%
    \renewcommand\transparent[1]{}%
  }%
  \providecommand\rotatebox[2]{#2}%
  \newcommand*\fsize{\dimexpr\f@size pt\relax}%
  \newcommand*\lineheight[1]{\fontsize{\fsize}{#1\fsize}\selectfont}%
  \ifx\svgwidth\undefined%
    \setlength{\unitlength}{98.92440874bp}%
    \ifx\svgscale\undefined%
      \relax%
    \else%
      \setlength{\unitlength}{\unitlength * \real{\svgscale}}%
    \fi%
  \else%
    \setlength{\unitlength}{\svgwidth}%
  \fi%
  \global\let\svgwidth\undefined%
  \global\let\svgscale\undefined%
  \makeatother%
  \begin{picture}(1,1.2433736)%
    \lineheight{1}%
    \setlength\tabcolsep{0pt}%
    \put(0,0){\includegraphics[width=\unitlength,page=1]{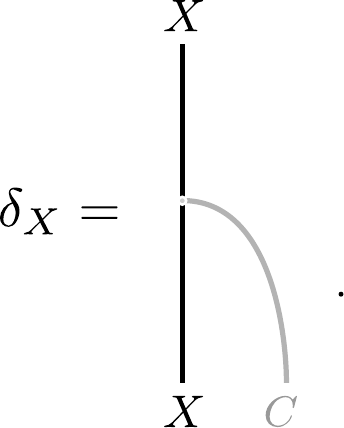}}%
  \end{picture}%
\endgroup%

}
\end{equation*}

\subsubsection{Extended factorization property of the coend of a braided rigid monoidal category}

In this section, we further assume that the rigid monoidal category \( \catV \) has a braiding \( \sigma \).

\begin{lemma}
\label{Res:coend-EFP}
	Let \( (\catV, \sigma) \) be a braided rigid monoidal category that admits a coend \( (C, \delta) \). For all \( n \geq 1 \), if there is an object \( D \) with a natural transformation 
\begin{equation*}
	\alpha = \{ \alpha_{X_1, \dots, X_n}: X_1 \otimes \cdots \otimes X_n \to X_1 \otimes \cdots \otimes X_n \otimes D \}_{(X_1, \dots, X_n) \in \catV^n}
\end{equation*}
then there exists a unique map \( f: C^{\otimes n} \to D \) such that the following diagram 
\begin{equation*}
\begin{tikzcd}
	X_1 \otimes \cdots \otimes X_n \ar[r, "\alpha_{X_1, \dots, X_n}"] \ar[d, "\delta_{X_1} \otimes \cdots \otimes \, \delta_{X_n}"'] & X_1 \otimes \cdots \otimes X_n \otimes D \\
	X_1 \otimes C \otimes \cdots \otimes X_n \otimes C \ar[r, "\beta"] & X_1 \otimes \cdots \otimes X_n \otimes C^{\otimes n} \ar[u, "\var{id} \otimes \, f"']
\end{tikzcd}
\end{equation*}
commutes. Here \( \beta \) is obtained by repeatedly applying the braiding \( \sigma \) of \( \catV \). In other words, for all \( n \geq 1 \) and all \( D \in \catV \), there is a bijection
\begin{align}
\begin{split}
	\Sigma^{(n)}_D: \var{Hom}(C^{\otimes n}, D) &\to  \Nat(\otimes^{(n-1)}, \otimes^{(n-1)} \otimes D) \\
	f &\mapsto (\var{id} \otimes f)\, \beta \, ( \otimes^{(n-1)} \delta^{\times n}),
\end{split}
\end{align}
where \( \otimes^{(n-1)} \) is the functor \( (X_1, \dots, X_n) \mapsto X_1 \otimes \cdots \otimes X_n \). 
\end{lemma}

\begin{proof}
	The lemma follows from induction on \( n \), the Fubini Theorem for coends \cite[Section IX.8]{MacLane-98:cat}, and the naturality of the braiding. See \cite[Lemma~5.4]{BV-12:double-Hopf-mnd} for the proof of a more general result.
\end{proof}

\begin{figure}[ht]
    \centering
     \scalebox{.85}{%
    \def\svgwidth{\columnwidth}
    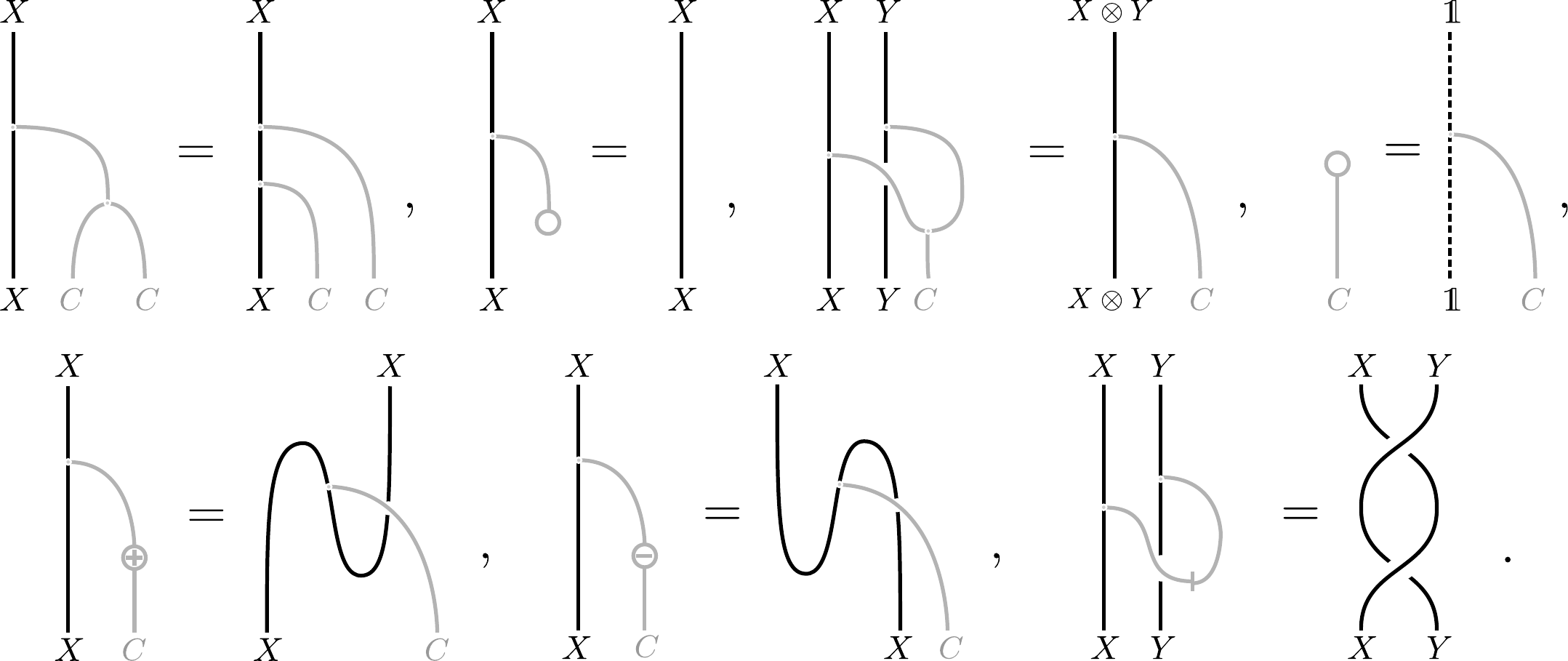
}
	\caption{Defining properties of the canonical Hopf algebra structure on the coend \( C \).}
	\label{Fig:coend-properties}
\end{figure}

As a consequence of the (extended) factorization property of coends, we obtain:

\begin{corollary}[{\cite[Sections 6.4, 6.5]{TV-17:monoidal-cat-TFT}}]
\label{Res:coend-structure-maps}
Let \( (\catV, \sigma) \) be a rigid braided category that admits a coend \( C \). Then \( C \) has a canonical Hopf algebra structure \( (C, m, u, \Delta, \epsilon, S) \) and a canonical Hopf pairing \( \omega: C \otimes C \to \ds1 \) compatible with the universal coaction \( \delta \) and the braiding \( \sigma \) as illustrated in  \textup{Figure~\ref{Fig:coend-properties}}. \qed
\end{corollary}

\begin{remark}[{\cite[Remark 8.2]{BV-12:double-Hopf-mnd}}]
\label{Rem:coend-self-coaction}
The universal coaction of \( C \) on itself can be expressed in terms of its Hopf algebra structure by \( \delta_C = (\var{id} \otimes m)(\sigma \otimes \var{id})(S \otimes \Delta)\Delta \). 
\end{remark}

We also define two variants \( \overline{\omega} \) and \( \underline{\omega} \) of \( \omega \), depicted graphically by
\begin{equation*}
	\scalebox{0.3}{%
    \def\svgwidth{\columnwidth}
\begingroup%
  \makeatletter%
  \providecommand\color[2][]{%
    \errmessage{(Inkscape) Color is used for the text in Inkscape, but the package 'color.sty' is not loaded}%
    \renewcommand\color[2][]{}%
  }%
  \providecommand\transparent[1]{%
    \errmessage{(Inkscape) Transparency is used (non-zero) for the text in Inkscape, but the package 'transparent.sty' is not loaded}%
    \renewcommand\transparent[1]{}%
  }%
  \providecommand\rotatebox[2]{#2}%
  \newcommand*\fsize{\dimexpr\f@size pt\relax}%
  \newcommand*\lineheight[1]{\fontsize{\fsize}{#1\fsize}\selectfont}%
  \ifx\svgwidth\undefined%
    \setlength{\unitlength}{200.38381629bp}%
    \ifx\svgscale\undefined%
      \relax%
    \else%
      \setlength{\unitlength}{\unitlength * \real{\svgscale}}%
    \fi%
  \else%
    \setlength{\unitlength}{\svgwidth}%
  \fi%
  \global\let\svgwidth\undefined%
  \global\let\svgscale\undefined%
  \makeatother%
  \begin{picture}(1,0.26875576)%
    \lineheight{1}%
    \setlength\tabcolsep{0pt}%
    \put(0,0){\includegraphics[width=\unitlength,page=1]{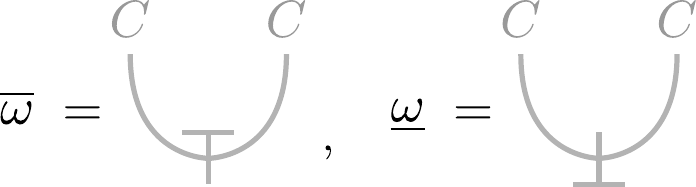}}%
  \end{picture}%
\endgroup%

}
\end{equation*}
which will allow us to freely switch between the braiding \( \sigma \) and its reverse braiding \( \overline{\sigma} \), as follows. The pairing \( \underline{\omega}  \) is the convolution inverse of \( \omega \) in the sense that
\begin{equation*}
	(\omega \otimes \underline{\omega} ) (\var{id} \otimes \sigma \otimes \var{id})(\Delta \otimes \Delta) = \epsilon \otimes \epsilon = (\underline{\omega} \otimes \omega) (\var{id} \otimes \sigma \otimes \var{id})(\Delta \otimes \Delta).
\end{equation*}
The pairing \( \overline{\omega} \) is the canonical Hopf pairing for \( C \) when \( \catV \) is braided with the reverse braiding \( \overline{\sigma} \). 

\begin{figure}[ht]
    \centering
    \scalebox{0.74}{%
    \def\svgwidth{\columnwidth}
    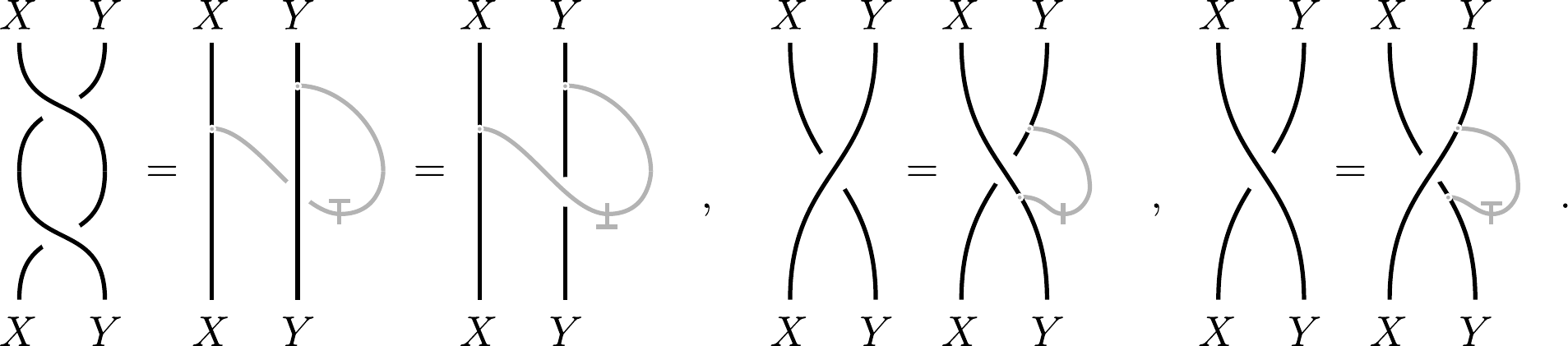
}
    \caption{The variants \( \overline{\omega} \) and \( \underline{\omega} \) of \( \omega \).}
    \label{Fig:braid-switch}
\end{figure}

\begin{lemma}
\label{Res:braid-switch}
Let \( \omega \) be the canonical Hopf pairing of the coend \( C \).
\begin{enumerate}[label=\upshape(\alph*)]
\item The diagrams in \textup{Figure \ref{Fig:braid-switch}} hold.
\item \( \overline{\omega} = \omega(S \otimes \var{id}) = \omega(\var{id} \otimes S^{-1})\overline{\sigma} \), \, \( \underline{\omega} = \omega(S^{-1} \otimes \var{id}) = \omega(S \otimes \var{id})\sigma \).
\end{enumerate}
\end{lemma}

\begin{figure}[ht]
    \centering
    \scalebox{.75}{%
    \def\svgwidth{\columnwidth}
    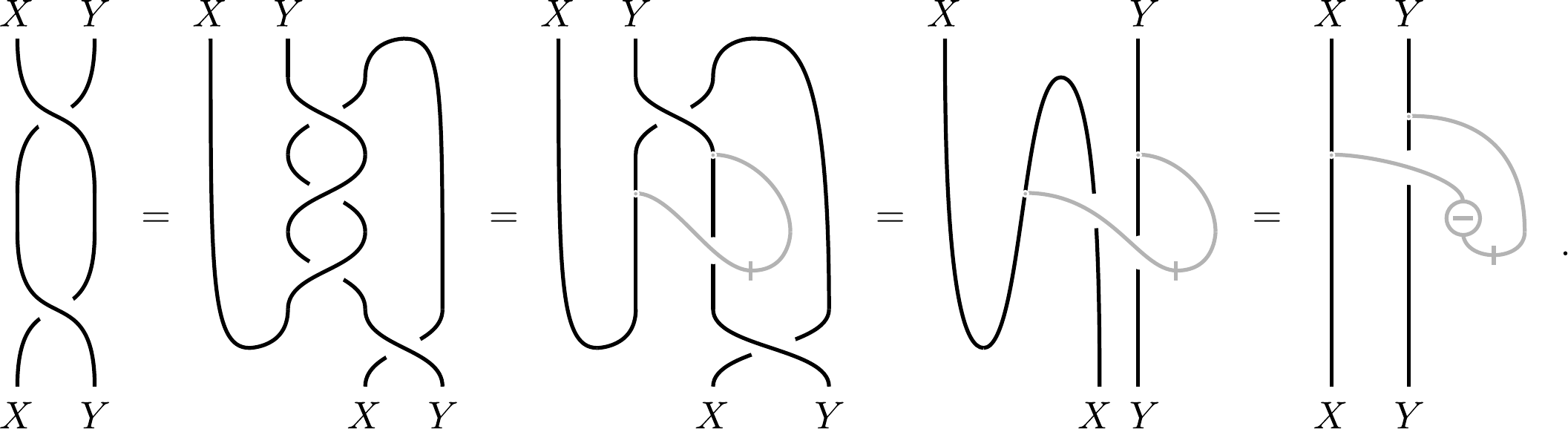
}
	\caption{Proof that \( \underline{\omega} = \omega(S^{-1} \otimes \var{id}) \).}
	\label{Fig:H-pairing-pf}
\end{figure}

\begin{proof}
Part (a) is immediate from the definitions of  \( \overline{\omega} \) and \( \underline{\omega} \). For part (b), the key is that the braiding \( \sigma \) can be expressed using \( \sigma^{-1} \) as shown in Figure~\ref{Fig:braiding}(b). The proof for one of the formulas of \( \underline{\omega} \) is shown in Figure~\ref{Fig:H-pairing-pf}, and the proofs of the other formulas are completely analogous.
\end{proof}

\begin{remark}
For some graphical calculus diagrams, we rotate certain morphisms such as the Hopf pairing \( \omega \) and the mulitplication map \( m_C \) following a \( C \)-coaction by 90 degrees counterclockwise, see for instance the first diagram in Figure~\ref{Fig:coend-H-mod} below. 
\end{remark}

\subsubsection{Examples of coends}

\begin{example}
Let \( \catC \) be a finite tensor \( \dsk \)-category, i.e., \( \catC \) is equivalent to the category \( \catM_A \) of finite dimensional modules over a finite dimensional \( \dsk \)-algebra \( A \), such that the bifunctor \( \otimes \) is bilinear on morphisms, and \( \var{End}(\ds1) \cong \dsk \). In this case, it is known that the coend of \( \catC \) exists, see \cite[Section~5.1.3]{KL-01:non-ss-TQFT} or \cite[Theorem~3.34]{Shimizu-17:unimod-FTC}. If \( \catC \) is addtionally semisimple so that \( \catC \) is a fusion category, we can exhibit the coend as
\begin{equation*}
	C = \bigoplus_{i \in \var{Irr}(\catC)} \dual i \otimes i,
\end{equation*}
where \( i \) runs over representatives of simple objects of \( \catC \), see \cite[Section~6.4.4]{TV-17:monoidal-cat-TFT}. When \( \catC = \Cat{Vec}^\var{fd} \), the coend is \( \dsk \) with the universal dinatural transformation given by the evaluation maps. 
\end{example}

\begin{figure}[ht]
	\centering
	\scalebox{.48}{%
    \def\svgwidth{\columnwidth}
    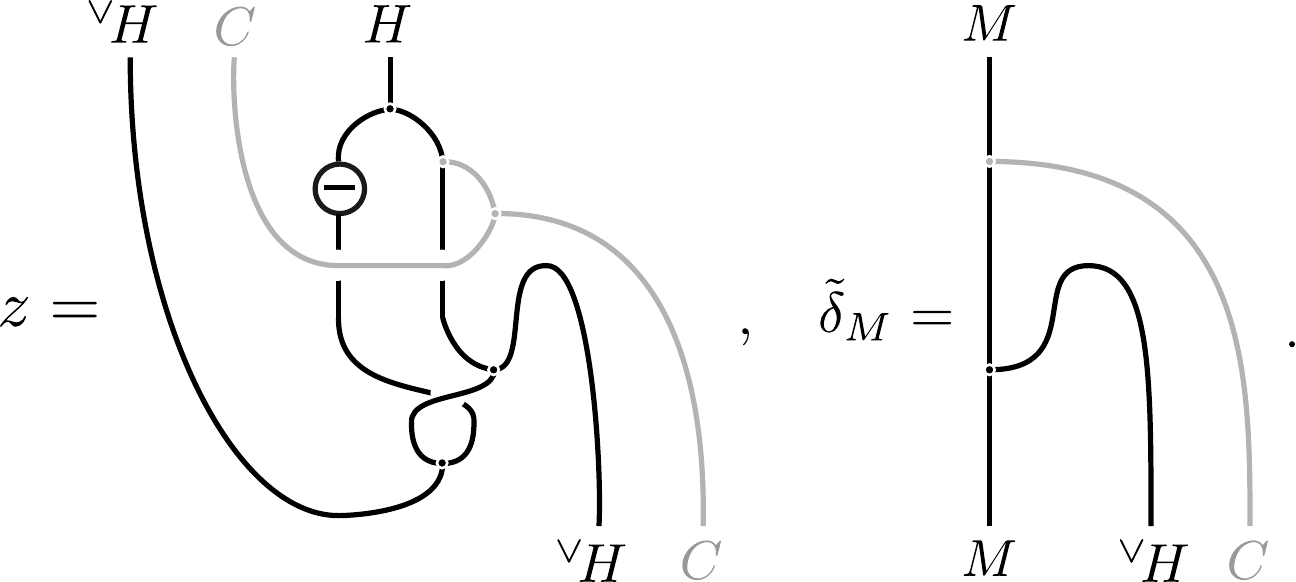
}
	\caption{The \( H \)-action and the universal coaction for the coend of \( \catV_H \).}
	\label{Fig:coend-H-mod}
\end{figure}

\begin{example}[{\cite[Section~6.3]{BV-12:double-Hopf-mnd}}]
\label{Eg:coend-H-mod}
Let \( \catV \) be a braided rigid monoidal category admitting a coend \( C \) with universal coaction \( \delta \), and let \( H \) be a Hopf algebra in \( \catV \). Then the coend of \( \catB = \catV_H \) exists and is given by \( Z(H) = \dual H \otimes C \). This is a right \( H \)-module with \( H \)-action \( z: Z(H) \otimes H \to Z(H) \) and universal coaction \( \tilde{\delta} \) given by Figure~\ref{Fig:coend-H-mod}.
\end{example}

\subsection{\texorpdfstring{Special elements of Hopf algebras over \( \dsk \)}{Special elements of Hopf algebras over k}}

\label{Sec:elts-Vec}

We summarize certain results on special elements of bialgebras and Hopf algebras over a field \( \dsk \) that we generalize in this paper, see e.g.\ \cite{Drabant-02:balanced-cat-HA, Radford-12:HA}. In this section, \( \otimes = \otimes_\dsk \), and we denote the symmetric braiding in \( \Cat{Vec}^\var{fd} \) by \( \tau \). We emphasize that since we work with right modules rather than left modules, certain definitions may appear different from the literature. 

\begin{definition}
Let \( H \) be a finite dimensional bialgebra over \( \dsk \).
\begin{enumerate}[label=(\roman*)]
\item An \( R \)-\emph{matrix} for \( H \) is an element \( R = R_i \otimes R^i \in H \otimes H \) satisfying 
	\begin{enumerate}[leftmargin=2.25\parindent]
		\item[(R1)] \( R \Delta(h) = \tau_{H,H}\Delta(h) R \) for all \( h \in H \).
		\item[(R2)] \( (\var{id} \otimes \Delta)(R) = R_i R_j \otimes R^j \otimes R^i \).
		\item[(R3)] \( (\Delta \otimes \var{id})(R) = R_i \otimes R_j \otimes R^iR^j \).
		\item[(R4)] \( (\var{id} \otimes \epsilon)(R) = (\epsilon \otimes \var{id})(R) = 1 \).
	\end{enumerate}
\item A \emph{grouplike element} of \( H \) is an element \( g \in H \) such that
\begin{enumerate}[leftmargin=2.25\parindent]
	\item[(G1)] \( \Delta(g) = g \otimes g \).
	\item[(G2)] \( \epsilon(g) = 1 \).
\end{enumerate}
\item If \( H \) is furthermore a Hopf algebra, we define a \emph{pivotal element} of \( H \) to be an element \( p \in H \) such that
\begin{enumerate}[leftmargin=2.25\parindent]
	\item[(P1)] \( p \) is a grouplike element, i.e., \( p \) satisfies (G1) and (G2).
	\item[(P2)] \( S^2(h) = p^{-1}hp \) for all \( h \in H \).
\end{enumerate}
\item Let \( R \) be an \( R \)-matrix for \( H \). A \emph{balanced element} (or \emph{twist}) of \( H \) with respect to \( R \) is an element \( t \in H \) such that
	\begin{enumerate}[leftmargin=2.25\parindent]
		\item[(T1)] \( t \) is central: \( th = ht \) for all \( h \in H \).
		\item[(T2)] \( \epsilon(t) = 1 \).
		\item[(T3)] \( \Delta(t) = (t \otimes t)(R_{21}R) \), where \( R_{21} = \tau_{H, H} R \).
	\end{enumerate}
	A balanced element is called a \emph{ribbon element} if it further satisfies
	\begin{enumerate}[leftmargin=2.25\parindent]
			\item[(T4)] \( S(t) = t \). 
	\end{enumerate}
\end{enumerate}
\end{definition}

\begin{proposition}
\label{Res:classic-Tannaka-duality}
Let \( H \) be a finite dimensional bialgebra over \( \dsk \).
\begin{enumerate}[label=\upshape(\roman*)]
\item There is a bijective correspondence between \( R \)-matrices for \( H \) and lax braidings on \( \catM_H \), given by \( R \mapsto (-) \cdot R_i \otimes (-) \cdot R^i \), and \( \sigma \mapsto \sigma_{H,H}(1_H \otimes 1_H) \).
\item If \( H \) is furthermore a Hopf algebra, then there is a bijective correspondence between pivotal elements of \( H \) and pivotal structures on \( \catM_H \), given by \( p \mapsto (-) \cdot p \), and \( \phi \mapsto \phi_H(1_H) \), where we have identified \( M \dual \dual \cong M \) as vector spaces for any \( H \)-module \( M \).
\item Suppose \( R \) is an \( R \)-matrix for \( H \). Then there is a bijective correspondence between balanced elements of \( H \) and balanced structures on \( \catM_H \), given by \( t \mapsto (-) \cdot t \), and \( \theta \mapsto \theta_H(1_H) \). Furthermore, a balanced element \( t \) is a ribbon element if and only if the corresponding balanced structure \( \theta \) is ribbon. \qed
\end{enumerate}
\end{proposition}

\begin{remark}
Let \( H \) be a finite dimensional Hopf algebra over \( \dsk \).
\begin{enumerate}[label=\upshape(\alph*)]
\item By Remark~\ref{Rem:braiding-rev} and Proposition~\ref{Res:classic-Tannaka-duality}(a), any \( R \)-matrix \( R \) for \( H \) is invertible, with \( R^{-1} = (\var{id} \otimes S^{-1})(R) \). As a result, the axiom \textup{(R4)} is automatically satisfied. Furthermore, if \( \sigma \) is the braid on \( \catM_H \) corresponding to \( R \), then the reverse braid \( \overline{\sigma} \) is given by \( \overline{R} = \tau_{H,H} R^{-1} \).
\item Similarly, any grouplike element \( g \in H \) is invertible with \( g^{-1} = S(g) \). As a result, the axiom \textup{(G2)} is automatically satisfied. 
\item If \( R \) is an \( R \)-matrix for \( H \), then any twist for \( H \) with respect to \( R \) is invertible by Remark~\ref{Rem:twist-invertible} and Proposition~\ref{Res:classic-Tannaka-duality}(c). As a result, the axiom \textup{(T2)} is automatically satisfied. 
\end{enumerate}
\end{remark}

\begin{definition}
Let \( H \) be a finite dimensional Hopf algebra over \( \dsk \) with an \( R \)-matrix \( R = R_i \otimes R^i \). We define \( u = R_i S(R^i) \), \( q = uS(u)^{-1} \), \( c = uS(u) \), and call \( u \) the (\emph{right}) \emph{Drinfeld element}.
\end{definition}

\begin{remark}
\label{Rem:Drinfelt-map-elt}
The morphisms \( \nu, \kappa \) and \( \gamma \) on \( \catM_H \) (see Section~\ref{Sec:Drinfeld-maps}) are given by the right actions of \( u, q \) and \( c \), respectively, if we identify \( M \cong M \dual \dual \cong M \dual \dual \dual \dual \) as vector spaces. As a result, the elements \( u, q, c \) are invertible, with \( u^{-1} \) corresponding to \( \bar{\nu} \), and therefore given by \( u^{-1} = S^2(R_i)R^i \). 
\end{remark}

The following corollary is a straightforward consequence of Lemma~\ref{Res:Drinfeld-bij}, Lemma~\ref{Res:Drinfeld-rbn-piv}, and Remark~\ref{Rem:Drinfelt-map-elt}. 

\begin{corollary}
\label{Res:Drinfeld-elt-bij}
Let \( (H, R) \) be a finite dimensional quasitriangular Hopf algebra over \( \dsk \). Then there is a one-to-one correspondence between the pivotal elements of \( H \) and the balanced elements of \( H \), given by \( p \mapsto u^{-1}p \) and \( t \mapsto ut \). Further, under the above correspondence, a pivotal element \( p \) maps to a balanced element \( t \) that is ribbon if and only if one of the following two equivalent conditions are satisfied:
\begin{enumerate}[label=\upshape(\roman*)]
\item \( p^2 = q \).
\item \( t^{-2} = c \). \qed
\end{enumerate}
\end{corollary}

\subsection{Quasitriangular braided Hopf algebras}
\label{Sec:C-qst-HA}

Let \( \catV \) be a braided rigid monoidal category admitting a coend \( C \), and let \( H \) be a bialgebra in \( \catV \). 

\begin{figure}[ht]
    \centering
    \scalebox{1}{%
    \def\svgwidth{\columnwidth}
    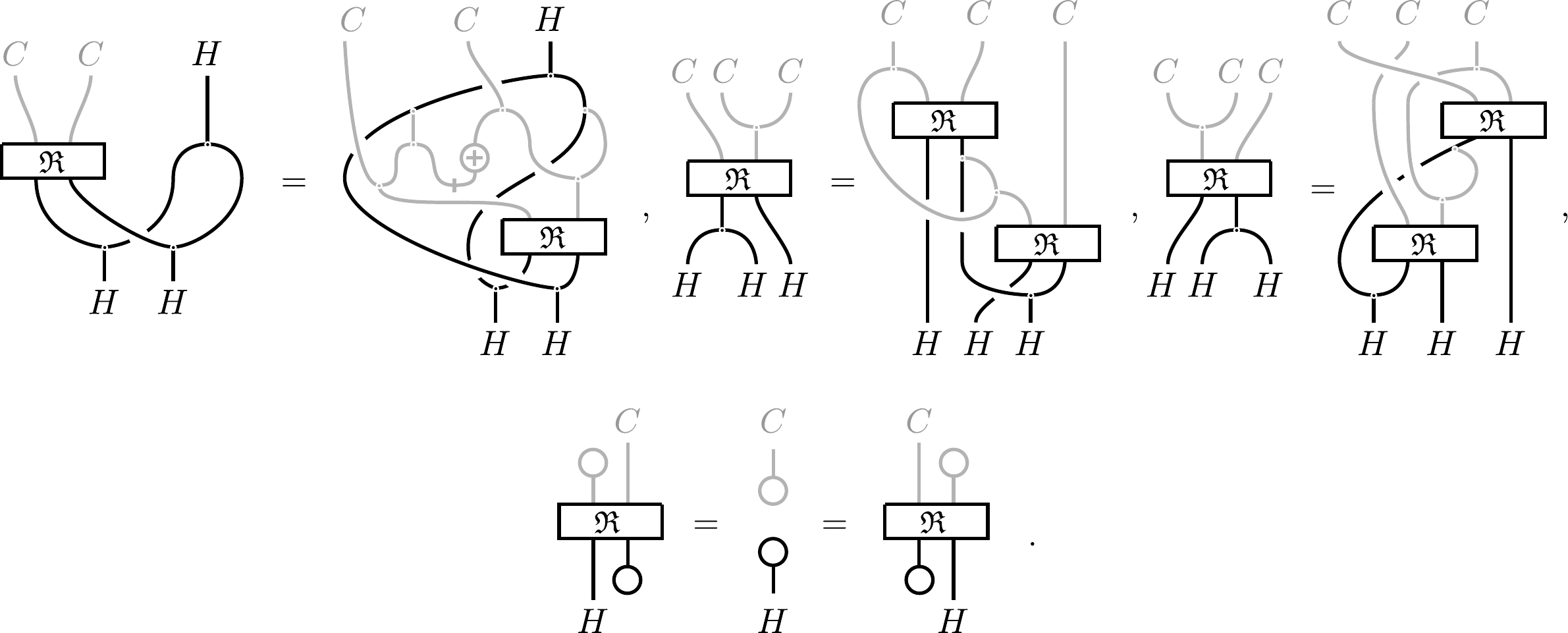
}
	\caption{Axioms of an \( R \)-matrix for a braided bialgebra or Hopf algebra.}
	\label{Fig:C-R-mat}
\end{figure}

\begin{definition}[{\cite[Section 8.6]{BV-12:double-Hopf-mnd}}]
\label{Def:C-qst-HA}
	A morphism \( \mathfrak{R}: C \otimes C \to H \otimes H \) satisfying the 4 diagrams of Figure~\ref{Fig:C-R-mat}  is called an \( R \)-\emph{matrix} for \( H \). A \emph{quasitriangular bialgebra} or \emph{quasitriangular Hopf algebra} is a bialgebra or Hopf algebra equipped with an \( R \)-matrix. 
\end{definition}

\begin{figure}[ht]
    \centering
    \scalebox{.27}{%
    \def\svgwidth{\columnwidth}
    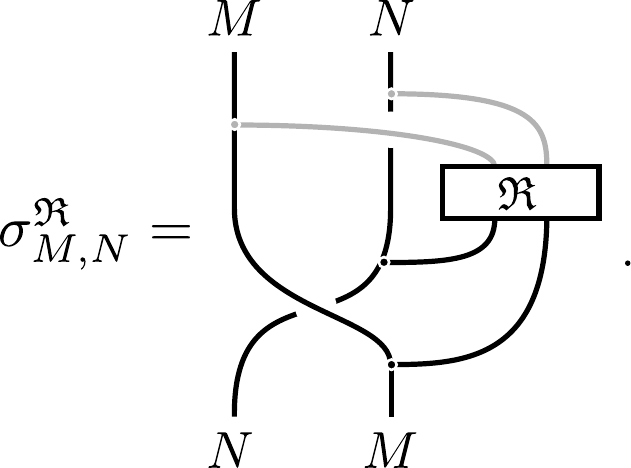
}
	\caption{The lax braiding on \( \catV_H \) induced by an \( R \)-matrix \( \mathfrak{R} \).}
	\label{Fig:C-R-mat-braiding}
\end{figure}

\begin{theorem}[{\cite[Section 8.6]{BV-12:double-Hopf-mnd}}]
\label{Res:C-R-mat-1-1}
Given an \( R \)-matrix \( \mathfrak{R}: C \otimes C \to H \otimes H \), let \( \sigma = \sigma^\mathfrak{R} \) be defined as in \textup{Figure~\ref{Fig:C-R-mat-braiding}} for any two right \( H \)-modules \( M \) and \( N \). Then \( \sigma \) is a lax braiding on \( \catV_H \) and the assignment \( \mathfrak{R} \mapsto \sigma^\mathfrak{R} \) is a bijection between \( R \)-matrices for \( H \) and lax braidings on \( \catV_H \). \qed 
\end{theorem}

We note that when \( H \) is a Hopf algebra, the lax braiding induced by an \( R \)-matrix for \( \mathfrak{R} \) is invertible, see Remark~\ref{Rem:C-R-invertible}.

\section{Elements of a Hopf monad induced by a Hopf algebra}
\label{Sec:M-elts}

This section is a summary of the relevant results of \cite{BV-07:Hopf-mnd}, with the exception of new materials in Section~\ref{Sec:remove-pivotal}, in the special context of a Hopf monad \( T_H \) induced by a Hopf algebra \( H \) in a braided rigid monoidal category \( \catV \). See Example~\ref{Eg:ind-Hopf-mnd} for the structure morphisms of this Hopf monad. We note that while the monadic perspective is more general, all the results in this section can be obtained by direct computation using graphical calculus. 

\subsection{\texorpdfstring{The monoid \( \Nat(1_\catV, 1_\catV \otimes H) \) of monadic elements}{The monoid Hom(1, 1 x H) of monadic elements}}

Let \( \catV \) be a monoidal category and let \( H = (H, m, u) \) be an algebra in \( \catV \). The category of modules over the monad \( T_H = 1_\catV \otimes H \) coincides with the category of right modules over \( H \) and hence is denoted by \( \catV_H \). Let \( U_H: \catV_H \to \catV \) and \( F_H: \catV \to \catV_H \) denote the forgetful functor and the free module functor, respectively. Recall that \( (F_H, U_H) \) is an adjunction with unit and counit given by
\begin{equation*}
\begin{alignedat}{2}
	 &\eta: 1_\catV \to U_H F_H   , &\qquad  \eta_X &= \var{id}_X \otimes u , \\
	 &\epsilon: F_H U_H \to 1_{\catV_H}, &\qquad \epsilon_{(M, r)} &= r,
\end{alignedat}
\end{equation*}
respectively, such that \( T_H = U_H F_H \). As a result, we can use graphical calculus to obtain the following lemma:

\begin{lemma}[{\cite[Lemma 1.3]{BV-07:Hopf-mnd}}]
\label{Res:mnd-adjunction}
Let \( \catD \) be any category and let \( F, G: \catV \to \catD \) be functors. We have mutually inverse bijections
\begin{align*}
\pushQED{\qed} 
(-)^{\sharp (F, G)}: \Nat(F, GT_H) &\leftrightarrows \Nat(FU_H, GU_H): (-)^{\flat (F, G)} \\
g &\mapsto g^{\sharp (F, G)}_{(M, r)} = G(r)g_M \\
 f_{XH } F(\eta_X) = f^{\flat (F, G)}_X &\mapsfrom f. \qedhere
 \popQED
\end{align*} 
\end{lemma}

\subsubsection{Convolution product}

\begin{figure}[ht]
    \centering
    \scalebox{.4}{%
    \def\svgwidth{\columnwidth}
    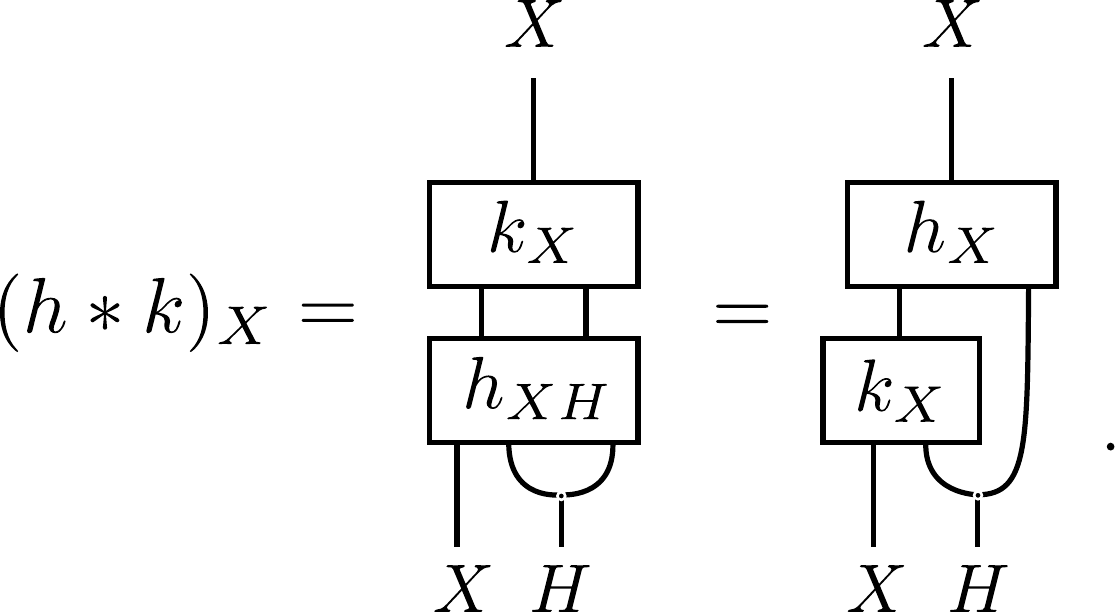
}
	\caption{Convolution product of \( \Nat(1_\catV, 1_\catV \otimes H) \) }
	\label{Fig:mnd-conv-prod}
\end{figure}

As a corollary of Lemma~\ref{Res:mnd-adjunction} when \( \catD = \catV \) and \( F = G = 1_\catV \), we obtain a bijection \( \Nat(1_\catV, 1_\catV \otimes H) \cong \End(U_H) \). In this case, \( \End(U_H) \) is also a monoid with function composition as multiplication. Transporting this monoid structure on \( \End(U_H) \) to \( \Nat(1_\catV, 1_\catV \otimes H) \), we obtain:

\begin{definition}
For \( h, k \in \Nat(1_\catV, 1_\catV \otimes H) \), define the \emph{convolution product} \( h * k \) as displayed in Figure~\ref{Fig:mnd-conv-prod}. Note the second equality follows from the naturality of \( h \).
\end{definition}

\begin{lemma}
\label{Res:mnd-bij-End(U)}
Let \( \catV \) be a monoidal category, and let \( H \) be an algebra in \( \catV \). The convolution product \( * \) gives \( \Nat(1_\catV, 1_\catV \otimes H) \) a monoid structure with unit \( 1_\catV \otimes u \), such that the mutually inverse bijections 
\begin{align}
\begin{split}
\label{Eqn:mnd-bijection}
(-)^\sharp: \Nat(1_\catV, 1_\catV \otimes H) &\leftrightarrows  \End(U_H): (-)^\flat \\
h &\mapsto h^\sharp_{(M, r)} = r h_M \\
 f_{XH} \eta_X = f	^\flat_X &\mapsfrom f
\end{split}
\end{align}
are isomorphisms of monoids. 
\qed
\end{lemma}

\begin{remark}
When \( \catV = \mathsf{Vec}^{\var{fd}} \), we have the more familiar mutually inverse bijections
\begin{align}
\label{Eqn:recon-bij}
\begin{split}
 H &\leftrightarrows  \End(U_H) \\
h  &\mapsto (-) \cdot h \\
 f_{(H, m)}(1_H) &\mapsfrom f
\end{split}
\end{align}
which is part of the theory of reconstructing \( H \) from its category of representations using the forgetful functor. In an abstract monoidal category we cannot speak of elements of \( H \), therefore the left hand side of (\ref{Eqn:recon-bij}) does not make sense. However, we always have the bijection (\ref{Eqn:mnd-bijection}). Thus, the elements of \( \Nat(1_\catV, 1_\catV \otimes H) \) can be seen as a substitute for the elements of \( H \) in this more general setting. For this reason, we will call these elements of \( \Nat(1_\catV, 1_\catV \otimes H) \) the \emph{monadic elements} of \( H \). 
\end{remark}

The investigation of certain special (pivotal, ribbon, etc.) monadic elements of \( H \) is the main focus of this section. The maps \( (-)^\sharp \) and \( (-)^\flat \), which replace the classical bijections in~(\ref{Eqn:recon-bij}), will play a key role. 

\subsubsection{Central elements}

\begin{figure}[ht]
    \centering
    \scalebox{.5}{%
    \def\svgwidth{\columnwidth}
    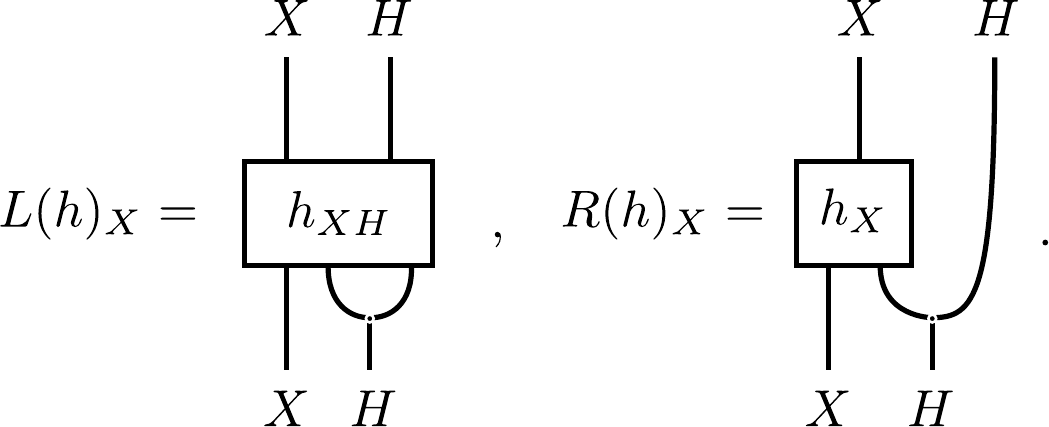
}
	\caption{Left and right multiplication in \( \Nat(1_\catV, 1_\catV \otimes H) \).}
	\label{Fig:mnd-L-R-mult}
\end{figure}

\begin{definition}
\label{Def:mnd-central-elt}
For a monadic element \( h \in \Nat(1_\catV, 1_\catV \otimes H) \), define the maps \( L_h, R_h \in \End(1_\catV \otimes H) \) as in Figure~\ref{Fig:mnd-L-R-mult}. If \( L_h = R_h \) then we say that \( h \) is \emph{central}. 
\end{definition}

For all \( k \in \Nat(1_\catV, 1_\catV \otimes H) \), we have \( L_h \circ k = h * k = R_k \circ h \), and hence \( L_h \) and \( R_h \) stand for \emph{left} and \emph{right multiplication} by \( h \), respectively. In particular, if \( h \) is central then \( h \) commutes with all elements in the monoid \( \Nat(1_\catV, 1_\catV \otimes H) \).

When \( \catV = \mathsf{Vec} \), the central elements of \( H \) are the ones whose multiplication maps are \( H \)-linear on all \( H \)-modules. We have an analogous result for the monadic central elements:

\begin{lemma}[{\cite[Lemma 1.5]{BV-07:Hopf-mnd}}]
For \( h \in \Nat(1_\catV, 1_\catV \otimes H) \), the following are equivalent:
\begin{enumerate}[label=\upshape(\roman*)]
\item \( h \) is central, i.e., \( L_h = R_h \),
\item \( h^\sharp_{(M,r)} \) is \( H \)-linear for all \( (M, r) \in \catV_H \), or equivalently, \( h^\sharp \in \End(U_H) \) lifts to an element in \( \End(1_{\catV_H}) \). \qed
\end{enumerate}
\end{lemma}

\begin{remark}
The proof of the preceding lemma follows from Lemma~\ref{Res:mnd-adjunction} with \( F = T_H \) and \( G = 1_\catV \). 
\end{remark}

\subsubsection{Monadic grouplike elements}

\begin{figure}[ht]
    \centering
    \scalebox{.6}{%
    \def\svgwidth{\columnwidth}
    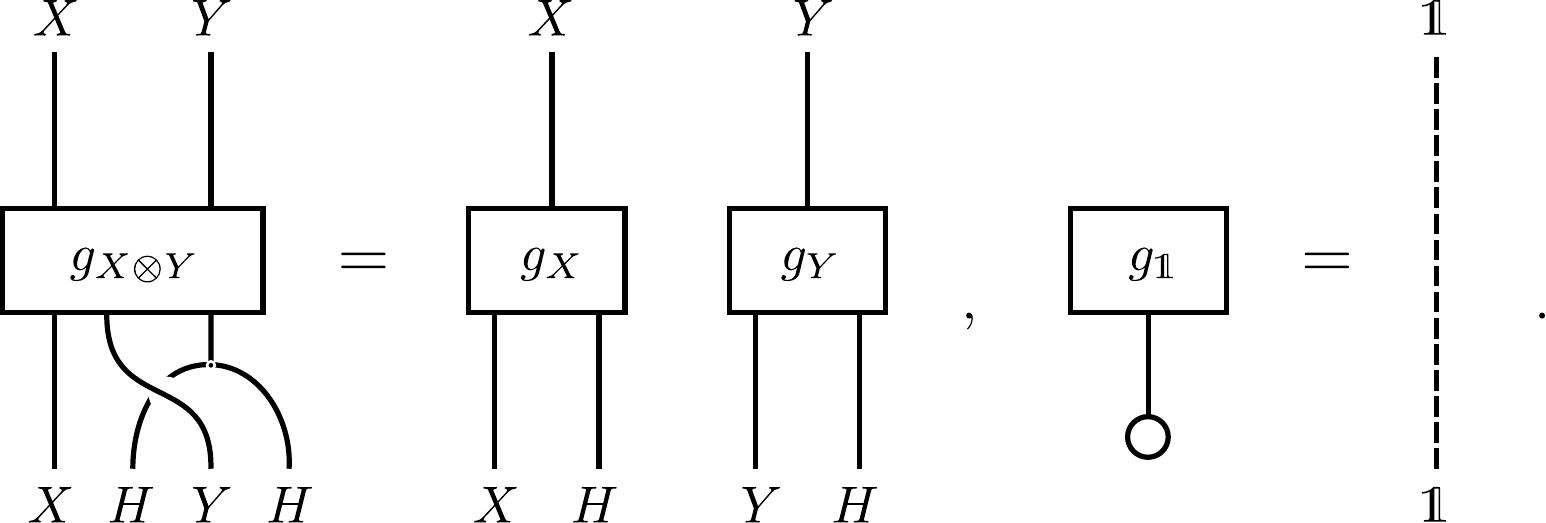
}
	\caption{Axioms for a monadic grouplike element \( g \).}
	\label{Fig:mnd-grplike}
\end{figure}

\begin{definition}
\label{Def:mnd-grplike}
Let \( H \) be a bialgebra in a braided monoidal category \( \catV \). A monadic element \( g \in \Nat(1_\catV, 1_\catV \otimes H) \) is \emph{grouplike} if \( g \) satisfies the two diagrams in Figure~\ref{Fig:mnd-grplike}. The set of monadic grouplike elements of \( H \) is denoted by \( \mathsf{MGrp}(H) \). 
\end{definition}

The following result is again a generalization of a well-known property of grouplike elements of a \( \dsk \)-bialgebra \( H \). 

\begin{lemma}[{\cite[Lemma 3.20]{BV-07:Hopf-mnd}}]
\label{Res:mnd-grplike}
The isomorphism~\textup{(\ref{Eqn:mnd-bijection})} restricts to an isomorphism \( \mathsf{MGrp}(H) \cong \End_\otimes(U_H) \), i.e., a monadic element \( g \) is grouplike if and only if \( g^\sharp \) is a monoidal natural endomorphism. \qed
\end{lemma}

\begin{remark}
\label{Rem:mnd-grplike-invertible}
By Lemmas~\ref{Res:mon-trans-rigid} and~\ref{Res:mnd-grplike}, when \( \catV \) is rigid and \( H \) is a Hopf algebra, any monadic grouplike element \( g \) is convolution invertible with \( \bar{g} = S(g) = S^{-1}(g) \). In particular, the second axiom \( \epsilon_H  g_\ds1 = \var{id}_\ds1 \) in Figure~\ref{Fig:mnd-grplike} is automatically satisfied.
\end{remark}

\subsubsection{Antipodes}

Let \( \catV \) be a braided rigid monoidal category. We have seen that there is an isomorphism of monoids \( \Nat(1_\catV, 1_\catV \otimes H) \cong \End(U_H) \). When \( H \) is a Hopf algebra, \( U_H \) is an endofunctor on the rigid monoidal category \( \catV_H \), so there is an antiautomorphism \( (-)^! \) on \( \End(U_H) \) as in Remark~\ref{Rem:antipode-End(F)}. We can further transport \( (-)^! \) to \( \Nat(1_\catV, 1_\catV \otimes H) \) to obtain the following definition.

\begin{definition}
\label{Def:mnd-antipode} 
Define maps \( S, S^{-1}: \Nat(1_\catV, 1_\catV \otimes H) \to \Nat(1_\catV, 1_\catV \otimes H) \) by \( S(f) =({}^!(f^\sharp))^\flat \) and \( S^{-1}(f) = ((f^\sharp)^!)^\flat \).
\end{definition}

\begin{figure}[ht]
    \centering
    \scalebox{.65}{%
    \def\svgwidth{\columnwidth}
    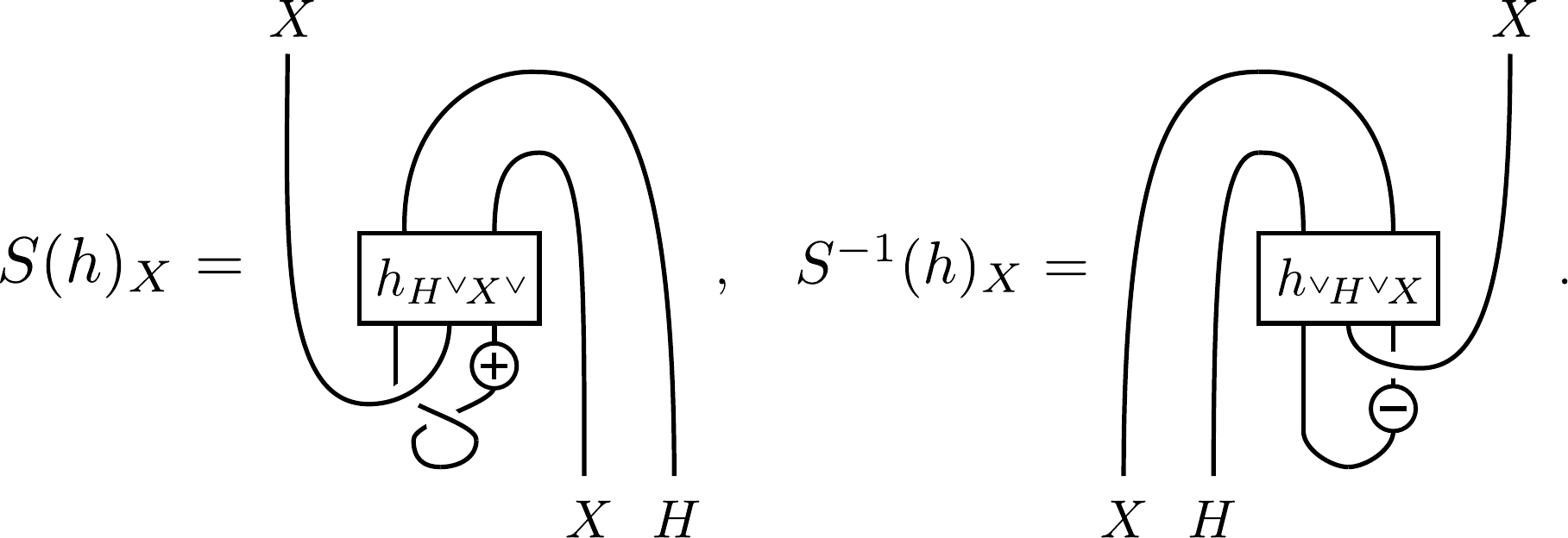
}
	\caption{Formulas for \( S^{\pm 1}(h) \) for \( h \in \Nat(1_\catV, 1_\catV \otimes H) \).}
	\label{Fig:mnd-antipode}
\end{figure}

\begin{remark}
As mentioned in the paragraph preceding \cite[Lemma 3.18]{BV-07:Hopf-mnd}, explicit formulas for \( S \) and \( S^{-1} \) are given by Figure~\ref{Fig:mnd-antipode}.
\end{remark}

\subsection{Monadic pivotal elements}
\label{Sec:mnd-piv}

In this section, we will also assume that \( \catV \) is a braided pivotal category with pivotal structure \( \psi \), or equivalently by Lemma~\ref{Res:Drinfeld-bij}, a balanced structure \( \theta = \bar{\nu} \psi \), and \( H \) is a Hopf algebra in \( \catV \). We present the theory of monadic pivotal elements as first introduced in \cite{BV-07:Hopf-mnd}, however, we will later present a slightly different view in Section~\ref{Sec:remove-pivotal}.

\subsubsection{The square of the antipode}

Recall that a pivotal element of a \( \dsk \)-Hopf algebra \( H \) is a grouplike element \( g \) such that \( S^2(h) = ghg^{-1} \) for all \( h \in H \). We can rewrite this as \( L_g(h) = R_g S^2(h) \), where \( L_g \) and \( R_g \) are left and right multiplication maps by \( g \),  respectively. The correct monadic generalization of the endomorphism \( S^2 \) is given by the following definition.

\begin{figure}[ht]
    \centering
    \scalebox{.24}{%
    \def\svgwidth{\columnwidth}
    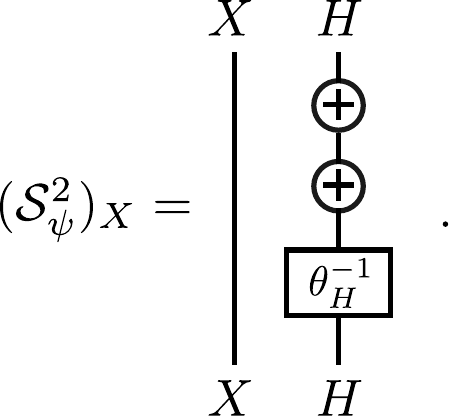
}
	\caption{The square of the antipode, \( \mathcal{S}^2_\psi \)}
	\label{Fig:mnd-S^2-psi}
\end{figure}

\begin{definition}
The \emph{square of the antipode} \( \mathcal{S}^2_\psi \in \End(1_\catV \otimes H) \) \emph{with respect to} \( \psi \) is defined as in Figure~\ref{Fig:mnd-S^2-psi}.
\end{definition}

\begin{lemma}[{\cite[Lemma~7.5]{BV-07:Hopf-mnd}}]
\label{Res:mnd-S^2}
For all \( g \in \Nat(1_\catV, 1_\catV \otimes H) \), the following are equivalent:
\begin{enumerate}[label=\upshape(\roman*), itemsep=0.5ex]
\item \( L_g = R_g \mathcal{S}^2_\psi \),
\item \( (\psi g^\sharp)_M: M \to M \dual \dual \) is \( H \)-linear for all \( M \in \catV_H \), i.e., \( \psi g^\sharp \) lifts to an element of \( \Nat(1_{\catV_H}, (-)^{\odual\odual}_{\catV_H}) \). \qed
\end{enumerate}
\end{lemma}

\subsubsection{Monadic pivotal elements}

\begin{definition}
A \emph{monadic \( \psi \)-pivotal element} of a braided Hopf algebra \( H \) is a monadic grouplike element \( p \in \mathsf{MGrp}(H) \) such that \( L_p = R_p \mathcal{S}^2_\psi \). The set of monadic \( \psi \)-pivotal elements of \( H \) is denoted by \( \mathsf{MPiv}_\psi(H) \). 
\end{definition}

As a consequence of Lemma~\ref{Res:mnd-grplike} and Lemma~\ref{Res:mnd-S^2}, we obtain:

\begin{theorem}[{\cite[Proposition 7.6]{BV-07:Hopf-mnd}}]
\label{Res:mnd-piv-1-1}
Let \( \catV \) be a braided pivotal category with pivotal structure \( \psi \), and let \( H \) be a Hopf algebra in \( \catV \). We have mutually inverse bijections 
\begin{equation*}
\psi(-)^\sharp: \mathsf{MPiv}_\psi(H) \leftrightarrows \mathsf{Piv}(\catV_H): (\psi^{-1}(-))^\flat \\
\end{equation*}
between the set \( \mathsf{MPiv}_\psi(H) \) of monadic \( \psi \)-pivotal elements and the set \( \mathsf{Piv}(\catV_H) \) of pivotal structures of \( \catV_H \).
\end{theorem}

\begin{proof}
A pivotal structure \( \phi \) on \( \catV_H \) is a natural collection of morphisms \( M \to M \dual \dual \) on the underlying object \( M \) of each \( H \)-module \( (M, r) \) that is (1) monoidal, and (2) \( H \)-linear. By composing with \( \psi^{-1} \), we see that \( \phi \) satisfies (1) if and only if \( (\psi^{-1} \phi)^\flat \) is monadic grouplike by Lemma~\ref{Res:mnd-grplike}. Similarly, \( \phi \) satisfies (2) if and only if \( L_p = R_p S^2_\psi \) by Lemma~\ref{Res:mnd-S^2}.
\end{proof}

\begin{remark}
If \( 1_\catV \otimes u_H \) is a monadic pivotal element of \( H \), then \( H \) is said to be \emph{involutory} \cite[Section~7.5]{BV-07:Hopf-mnd}. In this case, the corresponding pivotal structure on \( \catV_H \) is a lift of \( \psi \), i.e., the forgetful functor \( U_H \) is a pivotal functor.
\end{remark}

\subsection{Monadic R-matrices}

\begin{figure}[ht]
    \centering
    \scalebox{.7}{%
    \def\svgwidth{\columnwidth}
    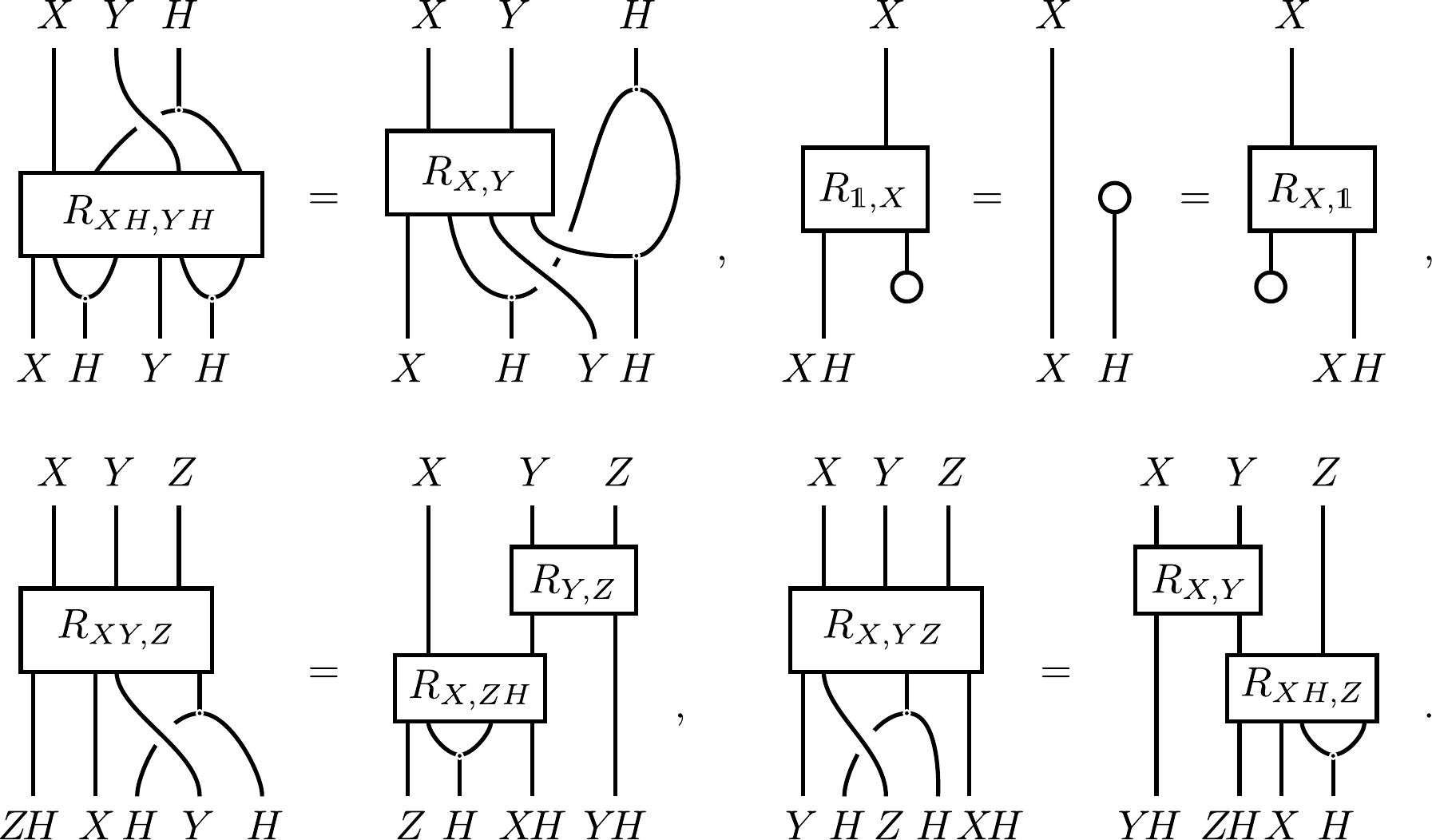
}
	\caption{The axioms of a monadic \( R \)-matrix \( R \).}
	\label{Fig:mnd-R-mat}
\end{figure}

\begin{definition}
	Let \( \catV \) be a braided monoidal category, and let \( H \) be a  bialgebra in \( \catV \). A \emph{monadic \( R \)-matrix} \( R \) for \( H \) is a natural transformation
\begin{equation*}
	\{ R_{X,Y}: X \otimes Y \to Y \otimes H \otimes X \otimes H \}_{X, Y \in \catV}
\end{equation*}
satisfying the axioms in Figure~\ref{Fig:mnd-R-mat}. A \emph{monadic quasitriangular bialgebra} (resp. \emph{Hopf algebra}) is a braided bialgebra (resp. Hopf algebra) \( H \) in \( \catV \) equipped with a monadic \( R \)-matrix.
\end{definition}

\begin{theorem}[{\cite[Theorem 8.5]{BV-07:Hopf-mnd}}]
\label{Res:mnd-R-1-1}
Let \( \catV \) be a braided monoidal category, and let \( H \) be a bialgebra in \( \catV \). Given a monadic \( R \)-matrix \( R \) for \( H \), we obtain a lax braiding \( \sigma^R \) on \( \catV_H \) by \( \sigma^R_{(M, r), (N, s)} = (s \otimes r) R_{M, N} \). Conversely, given a lax braiding \( \sigma \) on \( \catV_H \), we obtain a monadic \( R \)-matrix \( R^\sigma \) by \( R^\sigma_{X, Y} = \sigma_{X H, Y H} (\var{id} \otimes u_H \otimes \var{id} \otimes u_H) \). The correspondences \( R \mapsto \sigma^R \) and \( \sigma \mapsto R^\sigma \) are mutually inverse operations. \qed
\end{theorem}

\begin{remark}
\label{Rem:mnd-R-invertible}
When \( \catV \) is rigid and \( H \) is a Hopf algebra, any \( R \)-matrix \( R \) gives rise to an (invertible lax) braiding \( \sigma \). As a consequence, \( R \) is convolution invertible, and its convolution inverse \( \overline{R} \) corresponds to the reverse braiding \( \overline{\sigma} \). Using Figure~\ref{Fig:braiding}, we can express \( \overline{R} \) in terms of \( R \) as in Figure~\ref{Fig:mnd-R-rev} below, see also \cite[Corollary 8.7]{BV-07:Hopf-mnd}. Furthermore, in this case the axiom represented by the second diagram in Figure~\ref{Fig:mnd-R-mat}, which corresponds to the axiom \( \sigma_{X, \ds1} = \var{id}_\ds1 = \sigma_{\ds1, X} \) for all \( X \), is no longer needed, see Remark~\ref{Rem:braiding-rev}.
\end{remark}

\begin{figure}[ht]
    \centering
    \scalebox{.52}{%
    \def\svgwidth{\columnwidth}
    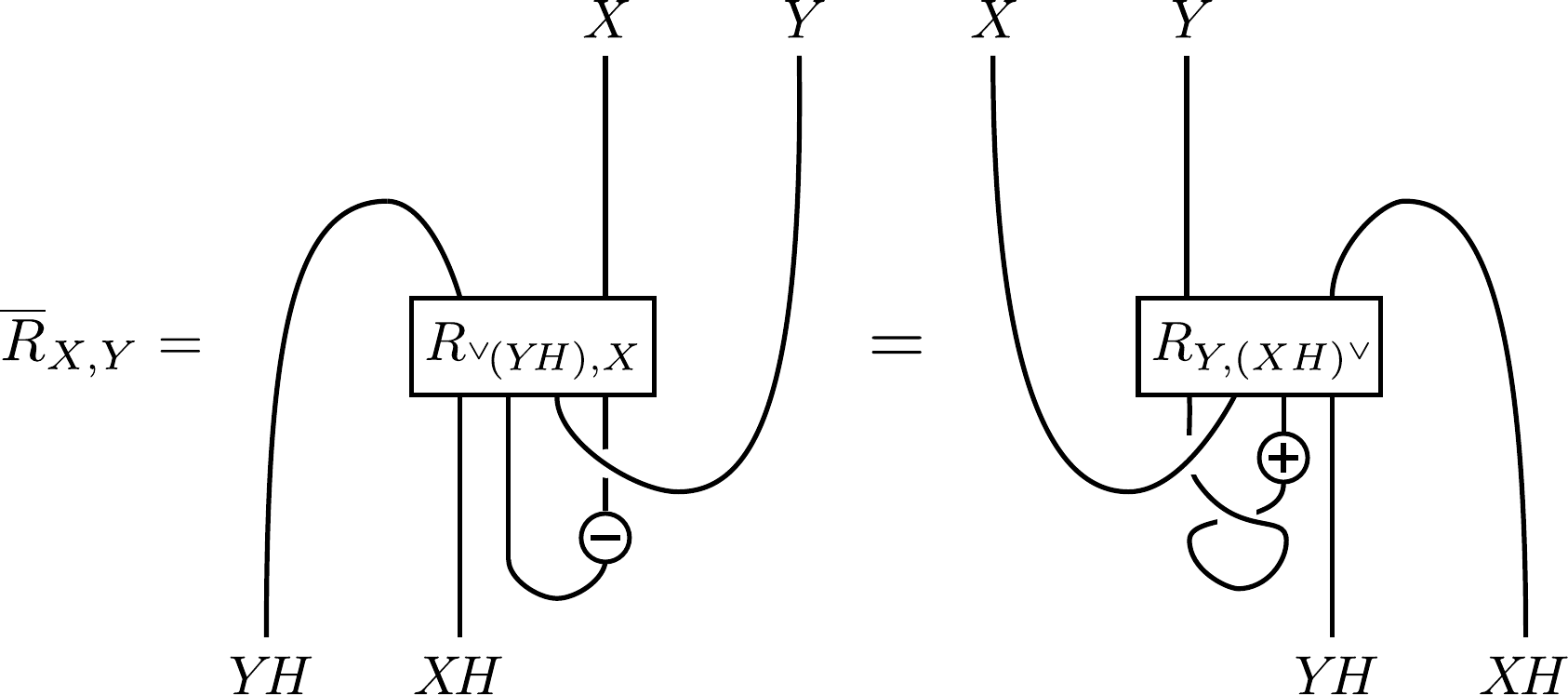
}
	\caption{The monadic reverse \( R \)-matrix \( \overline{R} \) corresponding to \( \overline{\sigma} \).}
	\label{Fig:mnd-R-rev}
\end{figure}

\subsection{Monadic ribbon elements}

\begin{figure}[ht]
    \centering
    \scalebox{.4}{%
    \def\svgwidth{\columnwidth}
    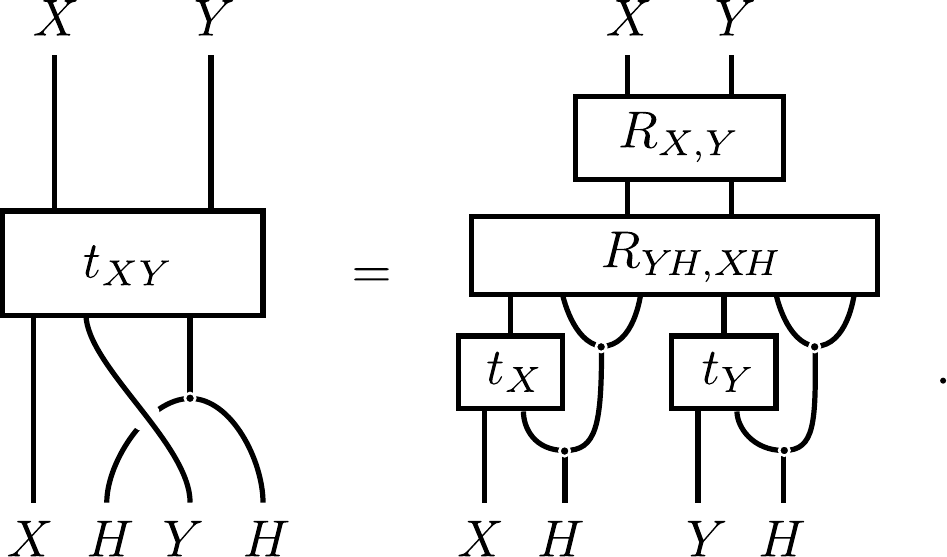
}
	\caption{The axiom (MT3) of a monadic balanced element.}
	\label{Fig:mnd-twist}
\end{figure}

\begin{definition}
Let \( (H, R) \) be a monadic quasitriangular bialgebra in a braided monoidal category \( \catV \). A \emph{monadic balanced element} (or \emph{monadic twist}) of \( H \) is a monadic element \( t \in \Nat(1_\catV, 1_\catV \otimes H) \) satisfying the following properties:
\begin{enumerate}
\item[(MT1)] \( t \) is central (in the sense of Definition~\ref{Def:mnd-central-elt}).
\item[(MT2)] \( \epsilon_H t_\ds1 = \var{id}_\ds1 \).
\item[(MT3)] \( t \) satisfies the diagram in Figure~\ref{Fig:mnd-twist}.
\end{enumerate}
Furthermore, when \( \catV \) is rigid and \( H \) is a Hopf algebra, we say that \( t \) is \emph{ribbon} if \( t \) further satisfies
\begin{enumerate}
\item[(MT4)] \( S(t) = t \).
\end{enumerate}
The set of monadic balanced (resp. ribbon) elements of \( H \) is denoted by \( \mathsf{MBal}(H) \) (resp. \( \mathsf{MRbn}(H) \)).
\end{definition}

\begin{theorem}[{\cite[Theorem 8.13]{BV-07:Hopf-mnd}}]
\label{Res:mnd-rbn-1-1}
Let \( \catV \) be a braided monoidal category, and let \( (H, R) \) be a quasitriangular bialgebra in \( \catV \). We have mutually inverse bijections
\begin{equation*}
\mathsf{MBal}(H) \cong \mathsf{Bal}(\catV_H), \quad t \mapsto t^\sharp, \quad \theta^\flat \mapsfrom \theta,
\end{equation*}
between the set \( \mathsf{MBal}(H) \) of monadic balanced elements of \( H \) and the set \( \mathsf{Bal}(\catV_H) \) of balanced structures of \( \catV_H \), which further restricts to a bijection \( \mathsf{MRbn}(H) \cong \mathsf{Rbn}(\catV_H) \) when \( \catV \) is rigid and \( H \) is a Hopf algebra in \( \catV \). \qed
\end{theorem}

\begin{remark}
\label{Rem:mnd-twist-invertible}
When \( \catV \) is rigid and \( H \) is a Hopf algebra in \( \catV \), any monadic balanced element is convolution invertible and the axiom (MT2) can be omitted, see Remark~\ref{Rem:twist-invertible}.
\end{remark}

\subsection{The monadic Drinfeld element}
\label{Sec:mnd-Drinfeld}

Let \( \catV \) be a braided pivotal category with pivotal structure \( \psi \), and let \( H \) be a monadic quasitriangular Hopf algebra in \( \catV \) with monadic \( R \)-matrix \( R \). We present the theory of the monadic Drinfeld element as first introduced in \cite{BV-07:Hopf-mnd}, however, we will later present some modification to the theory in Section~\ref{Sec:remove-pivotal}.

By Proposition~\ref{Res:mnd-R-1-1}, the category \( \catV_H \) is a braided rigid monoidal category with braiding \( \sigma = \sigma^R \). Hence we can consider the right Drinfeld morphism \( \nu \) in this category, see Section~\ref{Sec:Drinfeld-maps}. To utilize the bijection \( \Nat(1_\catV, 1_\catV \otimes H) \cong \End(U_H) \), we compose the Drinfeld morphism \( \nu \) with \( \psi^{-1} \) to obtain an endomorphism \( \psi^{-1} \nu \in \End(U_H) \).

\begin{figure}[ht]
    \centering
    \scalebox{.3}{%
    \def\svgwidth{\columnwidth}
    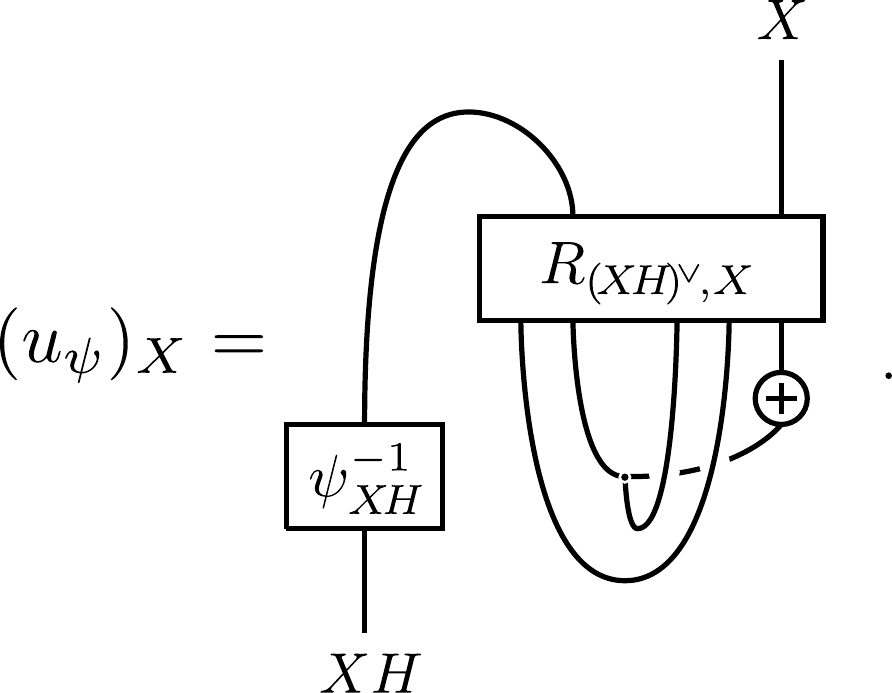
}
	\caption{Monadic \( \psi \)-Drinfeld element \( u_\psi \)}
	\label{Fig:mnd-Drinfeld}
\end{figure}

\begin{definition}
\label{Def:mnd-Drinfeld-elt}
Define the monadic \( \psi \)-Drinfeld element \( u_\psi \) as in Figure~\ref{Fig:mnd-Drinfeld}.
\end{definition}

\begin{lemma}[{\cite[Theorem 8.10]{BV-07:Hopf-mnd}}]
\label{Res:mnd-Drinfeld-1-1}
The monadic \( \psi \)-Drinfeld element \( u_\psi \) corresponds to \( \psi^{-1} \nu \) under the bijection \( \Nat(1_\catV, 1_\catV \otimes H) \cong \End(U_H) \) in \textup{(\ref{Eqn:mnd-bijection})}. \qed
\end{lemma}

\begin{theorem}[{\cite[Theorem 8.14]{BV-07:Hopf-mnd}}]
\label{Res:mnd-Drinfeld-bal-piv}
Let \( (H, R) \) be a monadic quasitriangular Hopf algebra on a braided pivotal category \( \catV \). The monadic \( \psi \)-Drinfeld element \( u_\psi \) induces a bijection
\begin{equation}
\label{Eqn:mnd-Drinfeld-bal-piv}
	u_\psi * (-): \mathsf{MBal}(H) \cong \mathsf{MPiv}_\psi(H)
\end{equation}
between the set of monadic balanced elements of \( H \) and the set of monadic \( \psi \)-pivotal elements of \( H \). 
\end{theorem}

\begin{proof}
By design, the diagram
\begin{equation*}
\begin{tikzcd}[column sep=large, row sep = large]
	\mathsf{Bal}(\catV_H) & \mathsf{Piv}(\catV_H) \\
	\mathsf{MBal}(H) & \mathsf{MPiv}_\psi(H)
	\arrow["\nu \circ (-)", "\cong" swap, from=1-1, to=1-2]
	\arrow["(-)^\sharp", "\cong" swap, from=2-1, to=1-1]
	\arrow["u_\psi \ast (-)" swap, from=2-1, to=2-2]
	\arrow["(\psi^{-1}(-))^\flat", "\cong" swap, from=1-2, to=2-2],
\end{tikzcd}
\end{equation*}
commutes, and three out of four arrows are bijections by Lemma \ref{Res:Drinfeld-bij}, Theorem \ref{Res:mnd-piv-1-1}, and Theorem  \ref{Res:mnd-rbn-1-1}.
\end{proof}

\begin{definition}
\label{Def:mnd-q-c} 
Define monadic elements \( q_\psi := u_\psi * S(\overline{u}_\psi) \) and \( c := u_\psi * S(u_\psi) \), where \( \overline{u}_\psi \) is the convolution inverse of \( u_\psi \) in \( \Nat(1_\catV, 1_\catV \otimes H) \).
\end{definition}

\begin{remark}
The element \( c \) does not depend on \( \psi \), in fact \( \psi \) no longer appears in \( c \) since it is cancelled by \( \psi^{-1} \). 
\end{remark}

\begin{theorem}[{\cite[Theorem 8.14, Corollary 8.15]{BV-07:Hopf-mnd}}]
\label{Res:mnd-Drinfeld-rbn-piv}
Under the bijection \textup{(\ref{Eqn:mnd-Drinfeld-bal-piv})}, a monadic balanced element \( t \) is ribbon if and only if one of the following equivalent conditions hold: 
\begin{enumerate}[label=\upshape(\roman*)]
\item The corresponding \( \psi \)-pivotal element \( p \) satisfies \( p^2 = q_\psi \).
\item \( t^{-2} = c \).
\end{enumerate}
\end{theorem}

\begin{proof}
Under the bijection \( \Nat(1_\catV, 1_\catV \otimes H) \cong \End(U_H) \), \( q_\psi \) and \( c \) correspond to  \( \psi^{-2} \kappa \) and \( \gamma \) respectively, see Definition~\ref{Def:Drinfeld-maps-2}. The result now follows from Lemma~\ref{Res:Drinfeld-rbn-piv}.
\end{proof}

\subsection{Removing the pivotal assumption}
\label{Sec:remove-pivotal}

Having described the theory of monadic elements as presented in \cite{BV-07:Hopf-mnd}, we introduce some small modifications to the theory. A general Hopf monad can be defined on any rigid monoidal category, not necessarily one equipped with a braiding. Therefore, the definitions of monadic pivotal elements and the monadic Drinfeld element require a pivotal structure \( \psi \) on \( \catV \), in order to utilize the bijection~(\ref{Eqn:mnd-bijection}). However, when the Hopf monad comes from a Hopf algebra, \( \catV \) is necessarily braided, and we have a canonical isomorphism \( X \to X \dual \dual \), namely the Drinfeld isomorphism. By utilizing the Drinfeld isomorphism, we can remove the pivotal assumption on \( \catV \). A small drawback is that the Drinfeld isomorphism is not monoidal, which requires us to slightly modify certain definitions and results.

\begin{remark}
We use \( \mu \) to denote the Drinfeld morphism in \( \catV \) and reserve \( \nu \) for the Drinfeld morphism in \( \catV_H \) coming from an \( R \)-matrix for \( H \). 
\end{remark}

\subsubsection{Monadic twisted grouplike elements}

\begin{figure}[ht]
    \centering
    \scalebox{.57}{%
    \def\svgwidth{\columnwidth}
    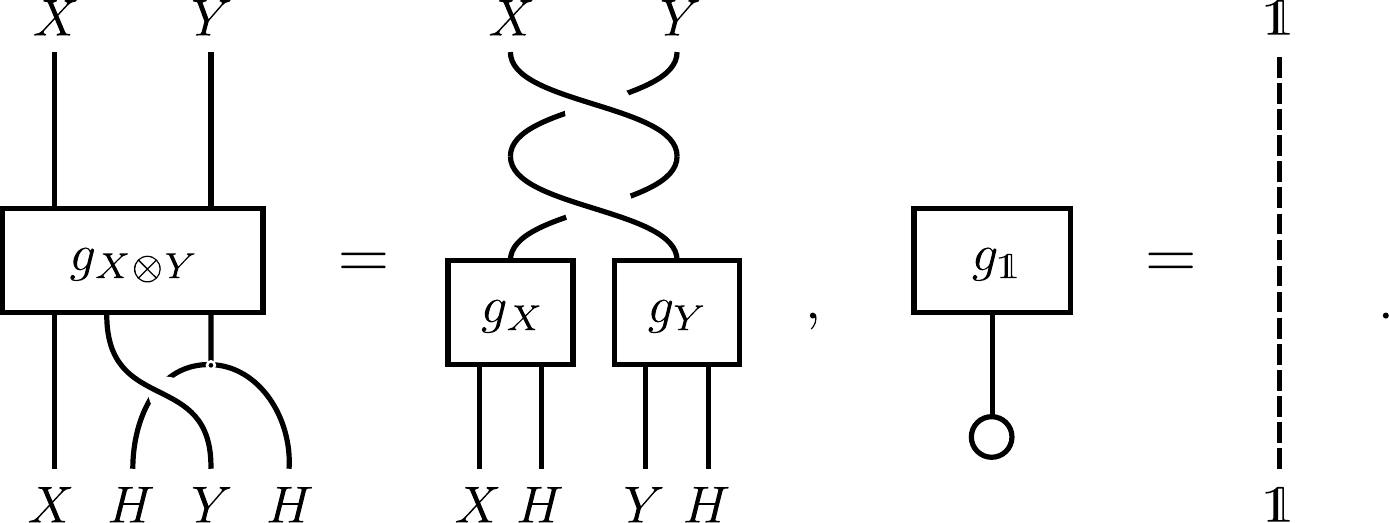
}
	\caption{Axioms for a monadic twisted grouplike element \( g \).}
	\label{Fig:mnd-twisted-grplike}
\end{figure}

\begin{definition}
\label{Def:mnd-twisted-grplike}
Let \( H \) be a bialgebra in a braided monoidal category \( (\catV, \sigma) \). A monadic element \( g \in \Nat(1_\catV, 1_\catV \otimes H) \) is \emph{twisted grouplike} if \( g \) satisfies the two diagrams in Figure~\ref{Fig:mnd-twisted-grplike}. 
\end{definition}

\begin{lemma}
\label{Res:mnd-twisted-grplike}
Under the isomorphism \( (-)^\sharp: \Nat(1_\catV, 1_\catV \otimes H) \cong \End(U_H) \), a monadic element \( g \) is twisted grouplike if and only if the corresponding natural endomorphism \( g^\sharp \) satisfies \( g^\sharp_{M \otimes N} = (g^\sharp_M \otimes g^\sharp_N) \, \sigma_{N, M} \sigma_{M, N} \) for all \( H \)-modules \( M \) and \( N \), and \( g^\sharp_\ds1 = \var{id}_\ds1 \).
\end{lemma}

\begin{proof}
Similar to the proof of Lemma~\ref{Res:mnd-grplike}, the only difference being the introduction of the double braiding.
\end{proof}

\subsubsection{The square of the antipode}

\begin{figure}[ht]
    \centering
    \scalebox{.17}{%
    \def\svgwidth{\columnwidth}
    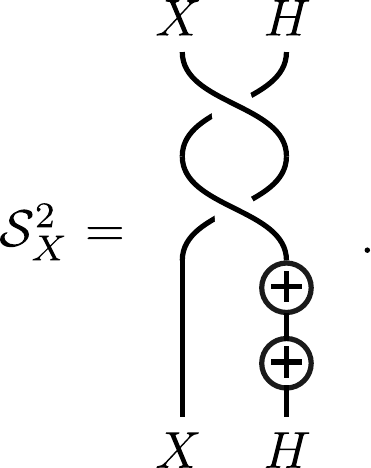
}
	\caption{The square of the antipode, \( \mathcal{S}^2 \)}
	\label{Fig:mnd-S^2}
\end{figure}

\begin{definition}
The \emph{square of the antipode} \( \mathcal{S}^2 \in \End(1_\catV \otimes H) \) is defined as in Figure~\ref{Fig:mnd-S^2}.
\end{definition}

\begin{lemma}
\label{Res:mnd-S^2-can}
For all \( g \in \Nat(1_\catV, 1_\catV \otimes H) \), the following are equivalent:
\begin{enumerate}[label=\upshape(\roman*), itemsep=0.5ex]
\item \( L_g = R_g \mathcal{S}^2 \),
\item \( (\mu g^\sharp)_M: M \to M \dual \dual \) is \( H \)-linear for all \( M \in \catV_H \), i.e., \( \mu g^\sharp \) lifts to an element of \( \Nat(1_{\catV_H}, (-)^{\odual\odual}_{\catV_H}) \).
\end{enumerate}
\end{lemma}

\begin{proof}
Similar to Lemma~\ref{Res:mnd-S^2}. Note that the double braiding comes from the fact that \( \mu_X \otimes \mu_H = \mu_{XH} \sigma_{H, X} \sigma_{X,H} \), see (\ref{Eqn:balancing-relation}), and the fact that the \( H \)-action on \( M \dual \dual \) can be described in terms of the double dual of the \( H \)-action on \( M \), the Drinfeld morphism, and the antipode of \( H \) squared. 
\end{proof}

\subsubsection{Monadic pivotal elements}

\begin{definition}
\label{Def:mnd-piv-elt-can}
A \emph{monadic pivotal element} of a braided Hopf algebra \( H \) is a monadic twisted grouplike element \( p \) such that \( L_p = R_p \mathcal{S}^2 \). The set of monadic pivotal elements of \( H \) is denoted by \( \mathsf{MPiv}(H) \). 
\end{definition}

As a consequence of Lemma~\ref{Res:mnd-twisted-grplike} and Lemma~\ref{Res:mnd-S^2-can}, we obtain:

\begin{theorem}
\label{Res:mnd-piv-1-1-can}
Let \( \catV \) be a braided rigid monoidal category, and let \( H \) be a Hopf algebra in \( \catV \). We have mutually inverse bijections 
\begin{equation*}
\mu(-)^\sharp: \mathsf{MPiv}(H) \leftrightarrows \mathsf{Piv}(\catV_H): (\bar{\mu}(-))^\flat \\
\end{equation*}
between the set \( \mathsf{MPiv}(H) \) of monadic pivotal elements and the set \( \mathsf{Piv}(\catV_H) \) of pivotal structures of \( \catV_H \). \qed
\end{theorem}

\begin{proof}
Same as the proof of Theorem~\ref{Res:mnd-piv-1-1}, however, since we compose a monoidal natural isomorphism with \( \bar{\mu} \), a double braiding is introduced, and therefore the corresponding monadic element is no longer grouplike but twisted grouplike.
\end{proof}

\subsubsection{The monadic Drinfeld element}

Let \( (H, R) \) be a monadic quasitriangular Hopf algebra in a braided rigid monoidal category \( \catV \), and let \( \mu \) and \( \nu \) denote the Drinfeld morphism in \( \catV \) and \( \catV_H \), respectively. We compose \( \nu \) with \( \bar{\mu} \) to obtain an endomorphism \( \bar{\mu} \nu \in \End(U_H) \), then use the isomorphism~(\ref{Eqn:mnd-bijection}) to obtain the following definition.

\begin{figure}[ht]
    \centering
    \scalebox{.3}{%
    \def\svgwidth{\columnwidth}
    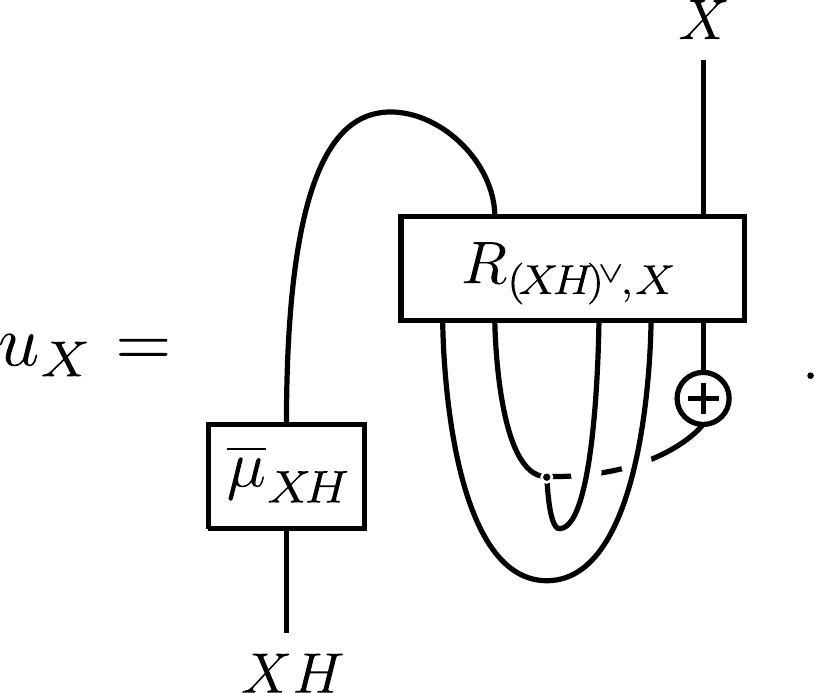
}
	\caption{The monadic Drinfeld element \( u \)}
	\label{Fig:mnd-Drinfeld-can}
\end{figure}

\begin{definition}
\label{Def:mnd-Drinfeld-elt-can}
Define the monadic Drinfeld element \( u \) as in Figure~\ref{Fig:mnd-Drinfeld-can}.
\end{definition}

Since the Drinfeld morphism is composed with \( \bar{\mu} \) instead of \( \psi^{-1} \), the formula for \( u \) is essentially the same as that of \( u_\psi \) but with \( \bar{\mu}_{XH} \) replacing \( \psi^{-1}_{XH} \), and therefore the proof of the following lemma is completely analogous to that of Lemma~\ref{Res:mnd-Drinfeld-1-1}.

\begin{lemma}
\label{Res:mnd-Drinfeld-1-1-can}
The monadic Drinfeld element \( u \) corresponds to \( \bar{\mu} \, \nu \) under the bijection \( \Nat(1_\catV, 1_\catV \otimes H) \cong \End(U_H) \) in \textup{(\ref{Eqn:mnd-bijection})}. \qed
\end{lemma}

\begin{theorem}
\label{Res:mnd-Drinfeld-bal-piv-can}
Let \( (H, R) \) be a monadic quasitriangular Hopf algebra on a braided rigid monoidal category \( \catV \). The monadic Drinfeld element \( u \) induces a bijection
\begin{equation}
\label{Eqn:mnd-Drinfeld-bal-piv-can}
	u * (-): \mathsf{MBal}(H) \cong \mathsf{MPiv}(H)
\end{equation}
between the set of monadic balanced elements of \( H \) and the set of monadic pivotal elements of \( H \). \qed
\end{theorem}

\begin{proof}
By design, the diagram
\begin{equation*}
\begin{tikzcd}[column sep=large, row sep = large]
	\mathsf{Bal}(\catV_H) & \mathsf{Piv}(\catV_H) \\
	\mathsf{MBal}(H) & \mathsf{MPiv}(H)
	\arrow["\nu \circ (-)", "\cong" swap, from=1-1, to=1-2]
	\arrow["(-)^\sharp", "\cong" swap, from=2-1, to=1-1]
	\arrow["u \ast (-)" swap, from=2-1, to=2-2]
	\arrow["(\bar{\mu}(-))^\flat", "\cong" swap, from=1-2, to=2-2],
\end{tikzcd}
\end{equation*}
commutes, and three out of four arrows are bijections by Lemma \ref{Res:Drinfeld-bij}, Theorem~\ref{Res:mnd-rbn-1-1}, and Theorem~\ref{Res:mnd-piv-1-1-can}. 
\end{proof}

\begin{definition}
\label{Def:mnd-q-c-can}
Let \( c_0 \) denote the automorphism \( \mu^! \mu  \), and let \( \bar{c}_0 \) denote its inverse, i.e., \( \bar{c}_0 = \bar{\mu}^!\bar{\mu} \). Also define two monadic elements \( q_\mu = (u * S(\bar{u})) \, \bar{c}_0 \) and \( c_\mu = (u * S(u)) \, c_0 \).
\end{definition}

\begin{theorem}
\label{Res:mnd-Drinfeld-rbn-piv-can}
Under the bijection \textup{(\ref{Eqn:mnd-Drinfeld-bal-piv-can})}, a monadic balanced element \( t \) is ribbon if and only if one of the following equivalent conditions hold: 
\begin{enumerate}[label=\upshape(\roman*)]
\item The corresponding pivotal element \( p \) satisfies \( p^2 = q_\mu \).
\item \( t^{-2} = c_\mu \). \qed
\end{enumerate}
\end{theorem}

\begin{proof}
Under the bijection (\ref{Eqn:mnd-Drinfeld-bal-piv-can}), \( p^2 \) corresponds to \( \bar{\mu}^2 \phi^2 \), while \( u * S(\bar{u}) \) corresponds to \( \bar{\mu} \mu^! \kappa \). When \( t \) is ribbon, we have \( \phi^2 = \kappa \) by Lemma~\ref{Res:Drinfeld-rbn-piv}, hence the difference between \( p^2 \) and \( u * S(\bar{u}) \) is given by \( \bar{c}_0 \). Similarly, the difference between \( t^{-2} \) and \( u * S(u) \) is given by \( c_0 \). 
\end{proof}

\section{Theory of coend elements}
\label{Sec:C-elts}

In this section, we fix a braided rigid monoidal category \( \catV \), which we further assume to admit a coend \( C \), see Section~\ref{Sec:coend}. We also fix a bialgebra \( H \) in \( \catV \). Our goal is to study certain (pivotal, ribbon, etc.) structures on the category \( \catV_H \) of right \( H \)-modules in terms of certain morphisms \( C \to H \). 

\subsection{\texorpdfstring{The monoid \( \var{Hom}(C, H) \) of coend elements of \( H \)}{The monoid Hom(C, H) of coend elements of H}}

\begin{figure}[ht]
    \centering
    \begin{subfigure}{0.45\textwidth}
        \centering
        \scalebox{0.43}{%
    \def\svgwidth{\columnwidth}
    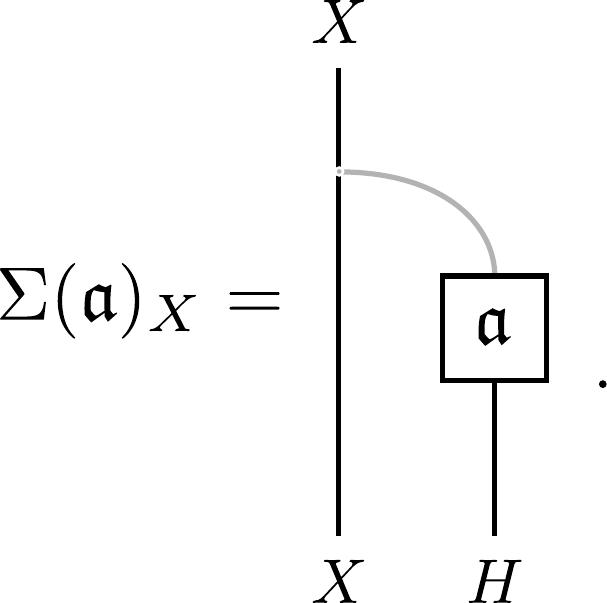
}
        \caption{}
    \end{subfigure}\hfill
    \begin{subfigure}{0.55\textwidth}
    	\centering
    	\scalebox{.75}{%
    \def\svgwidth{\columnwidth}
    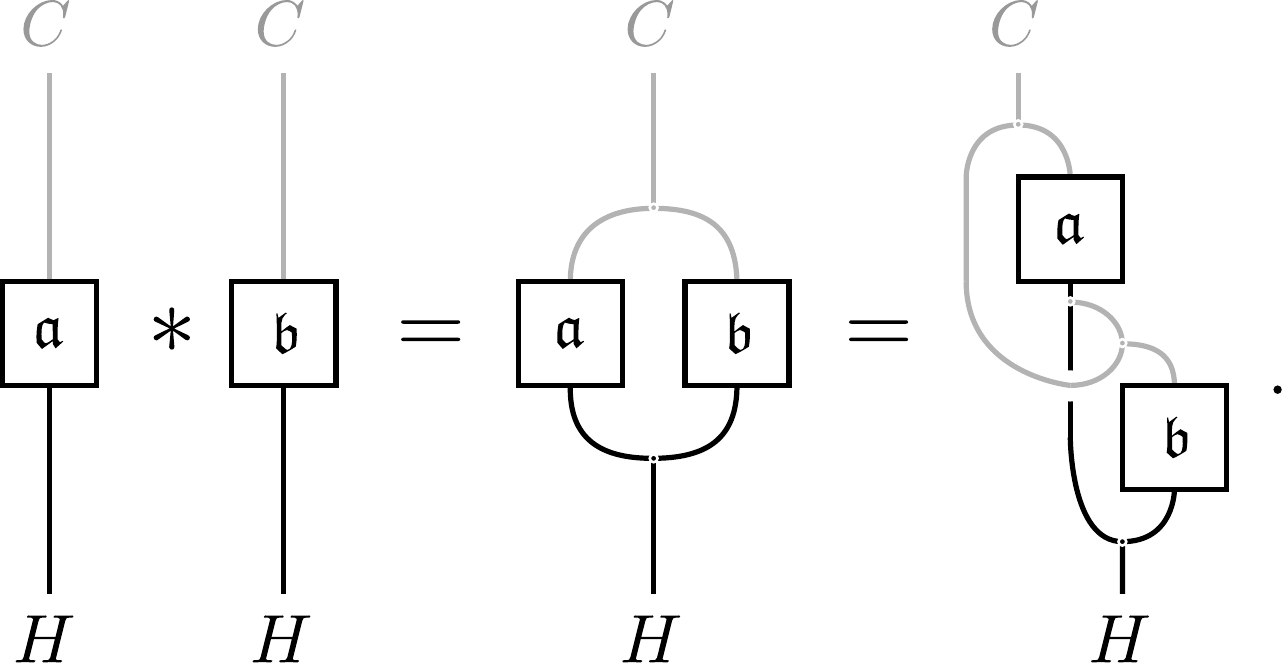
}
		\caption{}
    \end{subfigure}
    \caption{The map \( \Sigma \) and the convolution product of \( \var{Hom}(C,H) \).}
	\label{Fig:C-conv-prod}
\end{figure}

In this section, it suffices for \( H \) to be an algebra. Recall the factorization property of \( C \), in particular the bijection in~\eqref{Eqn:coend-fp},
\begin{align}
\label{Eqn:Sigma}
\begin{split}
	\Sigma = \Sigma^{(1)}_H: \var{Hom}(C, H) &\cong \Nat(1_\catV, 1_\catV \otimes H) \\
	\mathfrak{a} &\mapsto \Sigma(\mathfrak{a})_X := (\var{id}_X \otimes \mathfrak{a}) \, \delta_X
\end{split}
\end{align}
where \( \delta_X: X \to X \otimes C \) is the universal coaction. The map \( \Sigma \) is displayed in Fig~\ref{Fig:C-conv-prod}(a). 

\begin{definition}
\label{Def:C-elts}
The set \( \var{Hom}(C, H) \) is called the set of \emph{coend elements} (or simply \emph{\( C \)-elements}) of \( H \). 
\end{definition}

The bijection \( \Sigma \) is therefore a correspondence between the set of coend elements and the set of monadic elements described in Section 3.1. 

Since the set \( \Nat(1_\catV, 1_\catV \otimes H) \) of monadic elements is in fact a monoid (Lemma~\ref{Res:mnd-bij-End(U)}), we can use \( \Sigma \) to transport the opposite monoid structure of \( \Nat(1_\catV, 1_\catV \otimes H) \)  to obtain two formulas for the product \( * \) on \( \var{Hom}(C, H) \) as shown in Figure~\ref{Fig:C-conv-prod}(b), with unit \( \epsilon_C u_H \). In Hopf algebra literature, the first formula is the well-known convolution product of \( \var{Hom}(C, H) \). As a consequence of our setup, we obtain the following lemma.

\begin{lemma}
\label{Res:C-conv-prod}
The map \( \Sigma \) in \textup{(\ref{Eqn:Sigma})} is antimultiplicative and unital. Hence, \( \var{Hom}(C, H) \) with the product defined in \textup{Figure~\ref{Fig:C-conv-prod}(b)} with unit \( \epsilon_C u_H \) is a monoid, and \( \Sigma \) is an antiisomorphism of monoids. \qed
\end{lemma}

\begin{remark}
From now on, we will simply use concatenation for the convolution product on \( \var{Hom}(C, H) \). 
\end{remark}

\subsubsection{Central elements}

\begin{figure}[ht]
    \centering
    \scalebox{.3}{%
    \def\svgwidth{\columnwidth}
    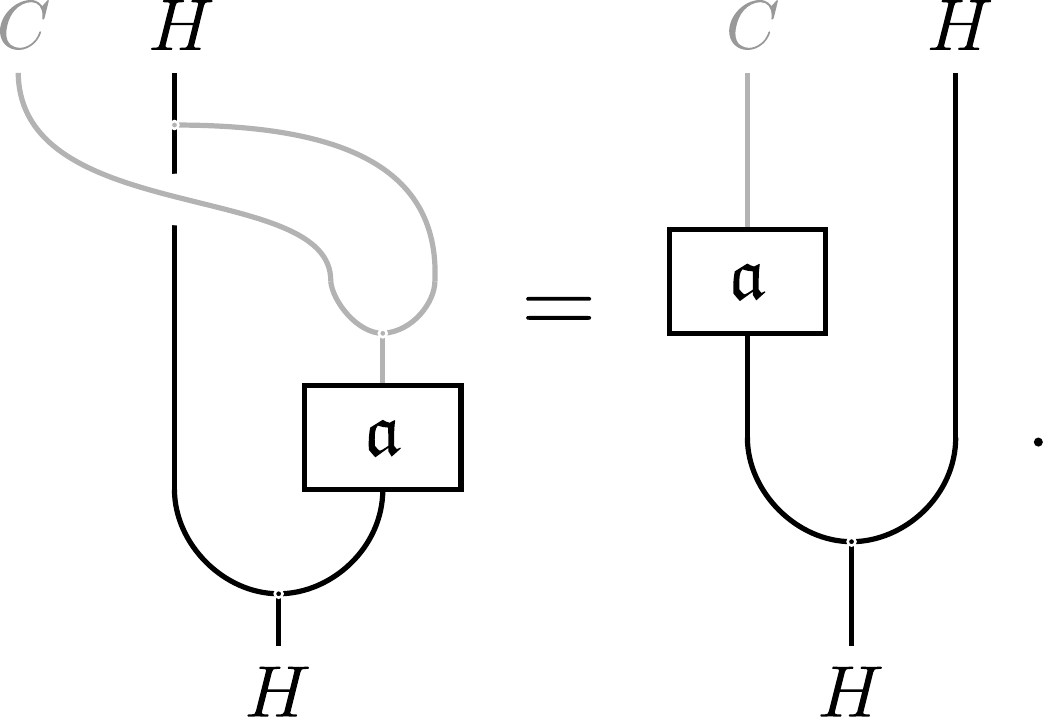
}
	\caption{Axiom of a central element in \( \var{Hom}(C, H) \).}
	\label{Fig:C-central-elt}
\end{figure}

\begin{definition}
\label{Def:C-central-elt}
A \( C \)-element \( \mathfrak{a}: C \to H \) is \emph{central} if \( \mathfrak{a} \) satisfies the diagram in Figure~\ref{Fig:C-central-elt}.
\end{definition}

\begin{lemma}
\label{Res:C-central-elt}
A \( C \)-element \( \mathfrak{a}: C \to H \) is central if and only if the corresponding monadic element \( a = \Sigma(\mathfrak{a}) \) is central, see \textup{Definition~\ref{Def:mnd-central-elt}}. 
\end{lemma}

\begin{proof}
Let \( \mathfrak{a}: C \to H \) be a \( C \)-element, and let \( a =\Sigma(\mathfrak{a}) \) be the corresponding monadic element. Recall that \( a  \) is central if and only if \( L_a = R_a \), where \( L_a \) and \( R_a \) are defined as in Figure~\ref{Fig:mnd-L-R-mult}. We can bend the incoming \( H \)-strand on both sides of the equation \( L_{\Sigma(\mathfrak{a})} = R_{\Sigma(\mathfrak{a})} \) downwards using the coevaluation morphism \( \var{coev}^l_H \), in order to apply \( \Sigma^{-1} \). Finally, by using \( \var{ev}^l_H \), we obtain Figure~\ref{Fig:C-central-elt}.
\end{proof}

\subsubsection{Coend grouplike elements}

Recall from Definition~\ref{Def:mnd-grplike} the notion of a monadic twisted grouplike element. The corresponding coend version of these elements is as follows.

\begin{figure}[ht]
    \centering
    \scalebox{.5}{%
    \def\svgwidth{\columnwidth}
    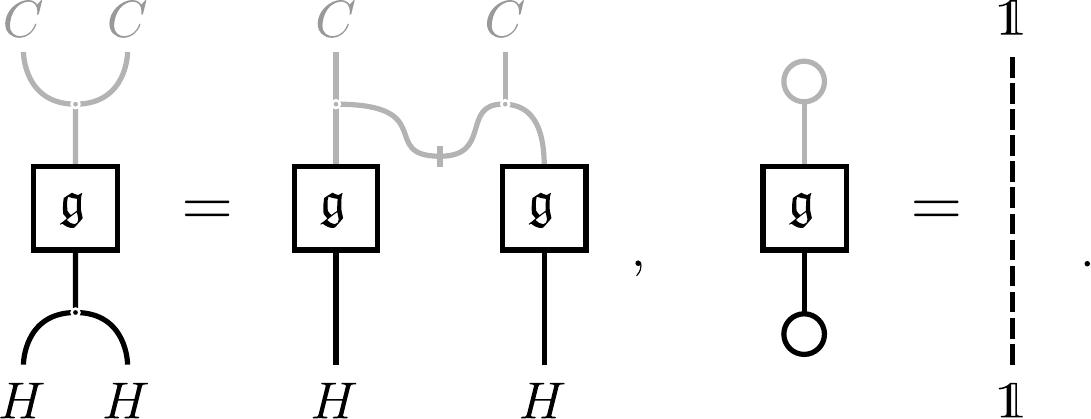
}
	\caption{Axioms of a grouplike element in \( \var{Hom}(C, H) \).}
	\label{Fig:C-grplike}
\end{figure}

\begin{definition}
A \( C \)-element \( \mathfrak{g}: C \to H \) is \emph{grouplike} if it satisfies the axioms presented in Figure~\ref{Fig:C-grplike}. Note that the first diagram uses the universal coaction of \( C \) on itself, for which there is a formula as shown in Remark~\ref{Rem:coend-self-coaction}.
\end{definition}

\begin{lemma}
\label{Res:C-grplike-M}
A \( C \)-element \( \mathfrak{g}: C \to H \) is grouplike if and only if the corresponding monadic element \( g = \Sigma(\mathfrak{g}) \) is monadic grouplike.
\end{lemma}

\begin{figure}[ht]
    \centering
    \scalebox{.85}{%
    \def\svgwidth{\columnwidth}
    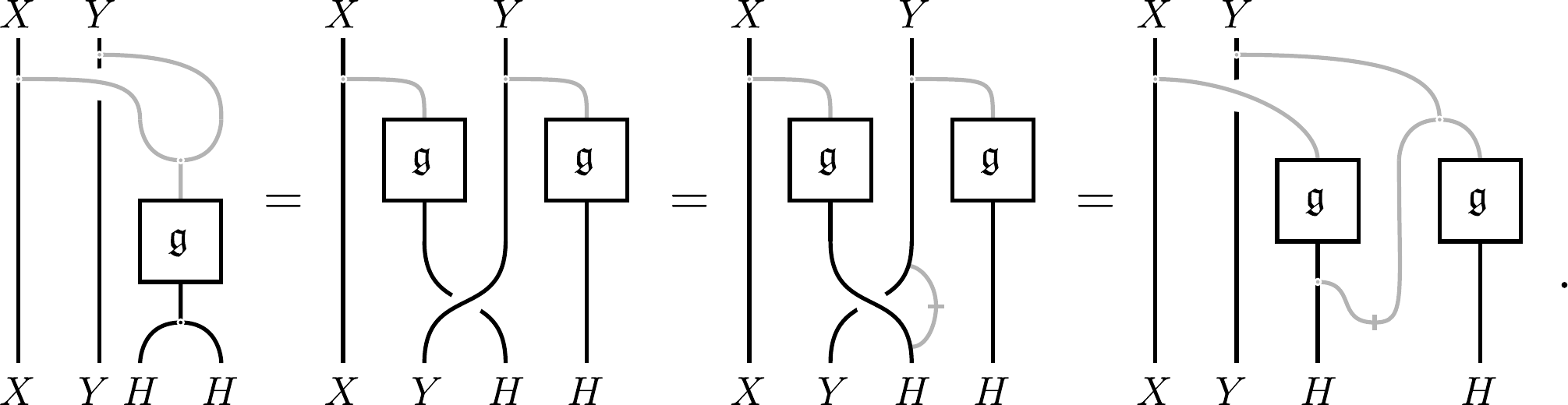
}
	\caption{Proof of the first defining property of coend grouplike elements.}
	\label{Fig:C-grplike-pf}
\end{figure}

\begin{proof}
Recall that the monadic element \( g = \Sigma(\mathfrak{g}) \) is grouplike if it satisfies the diagrams in Figure~\ref{Fig:mnd-grplike}. Clearly, the second diagram in Figure~\ref{Fig:C-grplike} corresponds to the second diagram in Figure~\ref{Fig:mnd-grplike}. It remains to show that the first diagrams in the two figures correspond to each other. The proof for this equivalence is given in Figure~\ref{Fig:C-grplike-pf}. Note that the second equality comes from Lemma~\ref{Res:braid-switch}.
\end{proof}

\subsubsection{Coend twisted grouplike elements}

\begin{figure}[ht]
    \centering
    \scalebox{.5}{%
    \def\svgwidth{\columnwidth}
    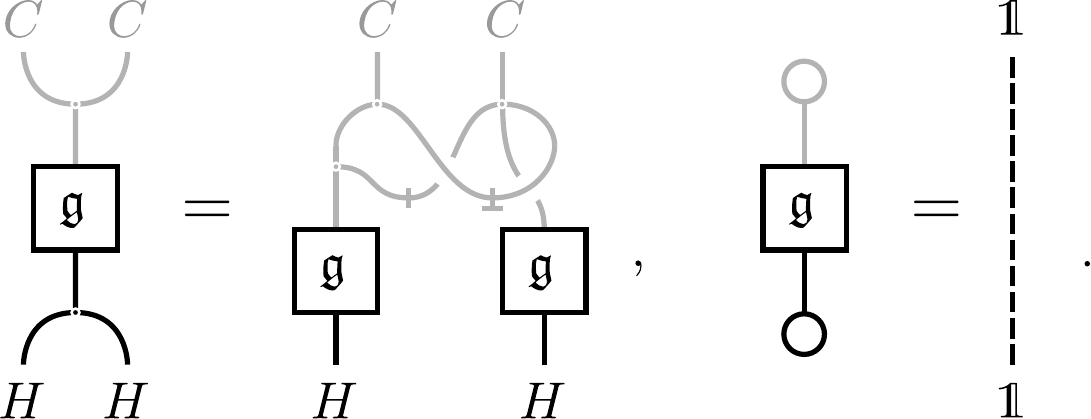
}
	\caption{Axioms of a twisted grouplike element in \( \var{Hom}(C, H) \).}
	\label{Fig:C-twisted-grplike}
\end{figure}

\begin{definition}
\label{Def:C-twisted-grplike}
A \( C \)-element \( \mathfrak{g}: C \to H \) is \emph{twisted grouplike} if it satisfies the axioms presented in Figure~\ref{Fig:C-twisted-grplike}.
\end{definition}

Since the definition of a monadic twisited grouplike element (Definition~\ref{Def:mnd-twisted-grplike}) only differs from that of a monadic grouplike element by a double braiding, the proof of the following result is a straightforward extension of the proof of Lemma~\ref{Res:C-grplike-M}.

\begin{lemma}
\label{Res:C-twisted-grplike-M}
A \( C \)-element \( \mathfrak{g}: C \to H \) is twisted grouplike if and only if the corresponding monadic element \( g = \Sigma(\mathfrak{g}) \) is twisted grouplike. \qed
\end{lemma}

\begin{remark}
If \( H \) is a Hopf algebra and not merely a bialgebra, then the second axiom \( \epsilon_H \mathfrak{g} u_C = \var{id}_\ds1 \) in Figure~\ref{Fig:C-grplike} and Figure~\ref{Fig:C-twisted-grplike} is not needed, see Remark~\ref{Rem:mnd-grplike-invertible}.
\end{remark}

\subsubsection{\texorpdfstring{The antipode \( S \) on \( \var{Hom}(C,H) \)}{The antiautomorphism S on Hom(C,H)}}

In this section \( H \) is always a Hopf algebra. Recall from Definition~\ref{Def:mnd-antipode} the antipode map \( S \) on \( \Nat(1_\catV, 1_\catV \otimes H) \). By using the bijection \( \Sigma \), we transport this map to \( \var{Hom}(C, H) \) with the following definition. 

\begin{figure}[ht]
    \centering
    \scalebox{.25}{%
    \def\svgwidth{\columnwidth}
    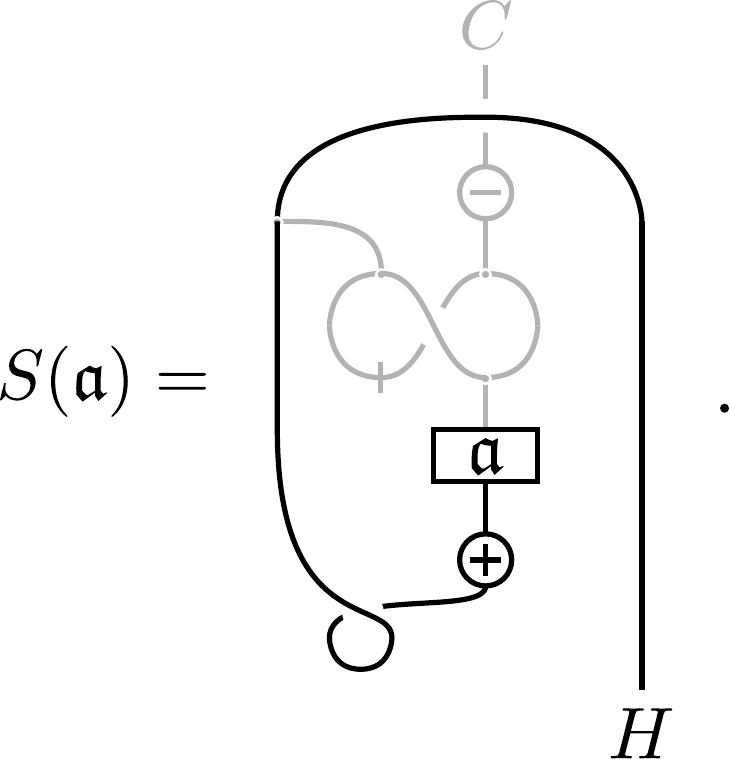
}
	\caption{The definition of the endomorphism \( S \) on \( \var{Hom}(C,H) \).}
	\label{Fig:C-antipode}
\end{figure}

\begin{definition}
Define an endomorphism \( S \) on \( \var{Hom}(C,H) \) as in Figure~\ref{Fig:C-antipode}.
\end{definition}

\begin{lemma}
For all \( \mathfrak{a} \in \var{Hom}(C, H) \), we have \( \Sigma( S(\mathfrak{a})) = S(\Sigma(\mathfrak{a})) \), where \( S(\mathfrak{a}) \) is defined as in \textup{Figure~\ref{Fig:C-antipode}}, and \( S(a) \) is defined as in \textup{Figure~\ref{Fig:mnd-antipode}} for \( a = \Sigma(\mathfrak{a}) \). As a consequence, \( S \) is an antiautomorphism on \( \var{Hom}(C,H) \).
\end{lemma}

\begin{figure}[ht]
    \centering
    \scalebox{1}{%
    \def\svgwidth{\columnwidth}
    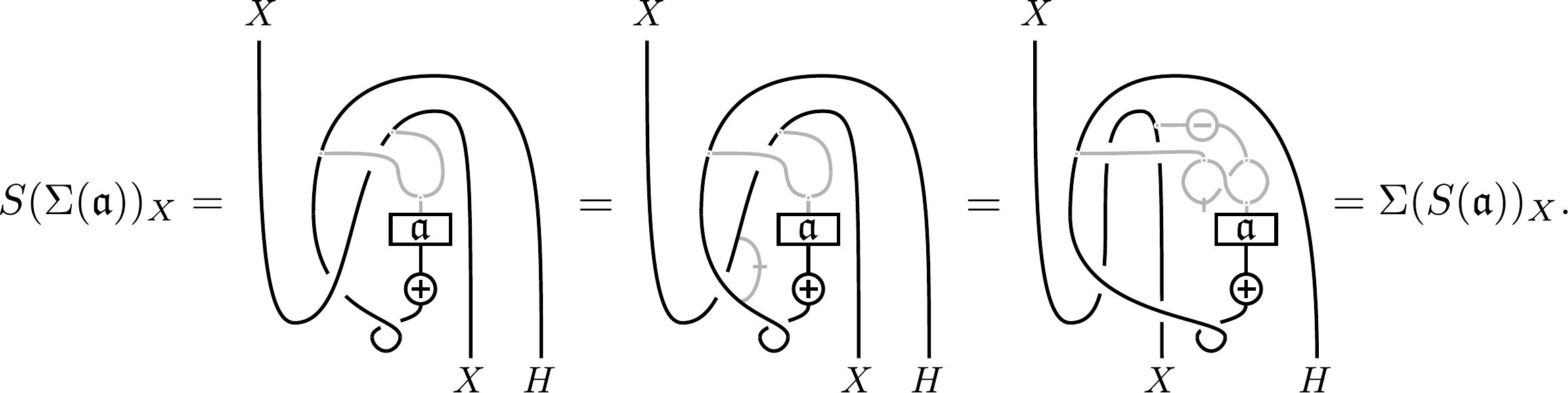
}
	\caption{Proof that \( \Sigma \circ S(\mathfrak{a}) = S \circ \Sigma(\mathfrak{a}) \).}
	\label{Fig:C-antipode-pf}
\end{figure}

\begin{proof}
See Figure~\ref{Fig:C-antipode-pf}. The first equality follows from the definition of \( S(\Sigma(\mathfrak{a})) \). The second equality follows from Lemma~\ref{Res:braid-switch}. Finally, we use the naturality of the braiding and the properties of \( C \) as in presented in Figure~\ref{Fig:coend-properties} of Corollary~\ref{Res:coend-EFP}. The second statement is a clear consequence of the first, since \( S \) is an antiautomorphism on \( \Nat(1_\catV, 1_\catV \otimes H) \) and \( \Sigma \) is antimultiplicative.
\end{proof}

\subsection{Coend pivotal elements}
\label{Sec:C-piv}

In this section, \( H \) is a Hopf algebra in a braided rigid monoidal category \( \catV \). Recall from Definition~\ref{Def:mnd-piv-elt-can} the notion of a monadic pivotal element.

\begin{figure}[ht]
    \centering
    \scalebox{.35}{%
    \def\svgwidth{\columnwidth}
    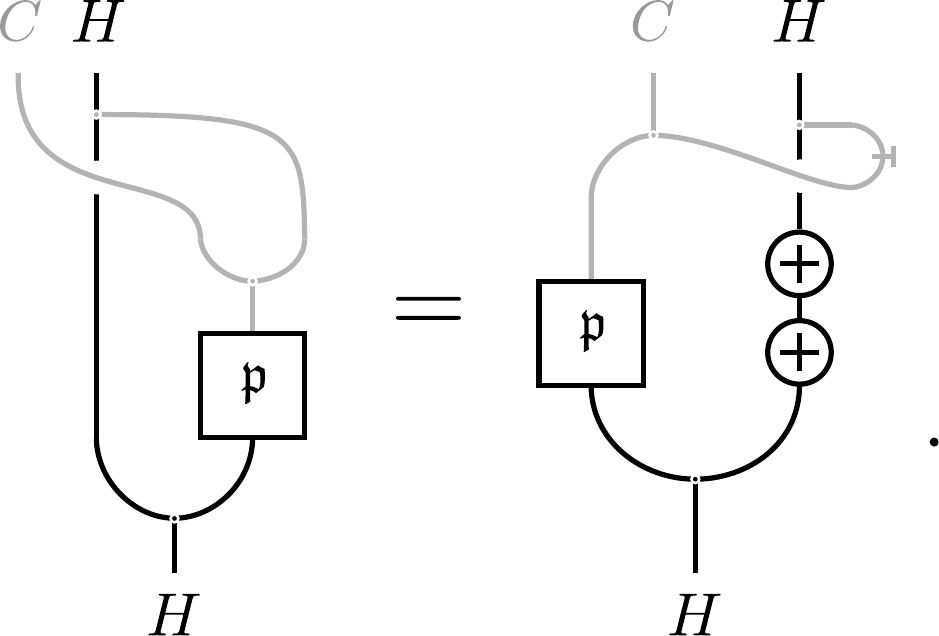
}
	\caption{The main axiom of a pivotal element in \( \var{Hom}(C, H) \).}
	\label{Fig:C-piv}
\end{figure}

\begin{definition}
\label{Def:C-piv}
A \emph{\( C \)-pivotal element} of \( H \) is a \( C \)-twisted grouplike element \( \mathfrak{p} \) (Definition~\ref{Def:C-twisted-grplike}) that satisfies the axiom in Figure~\ref{Fig:C-piv}. The set of \( C \)-pivotal elements of \( H \) is denoted by \( \mathsf{CPiv}(H) \). 
\end{definition}

\begin{lemma}
\label{Res:C-piv-M}
The bijection \( \Sigma \) in \textup{(\ref{Eqn:Sigma})} between coend and monadic elements of \( H \) restricts to a bijection \( \mathsf{MPiv}(H) \cong \mathsf{CPiv}(H) \) between the respective pivotal elements. 
\end{lemma}

\begin{proof}
We have already seen in Lemma~\ref{Res:C-twisted-grplike-M} that a \( C \)-twisted grouplike element \( \mathfrak{p} \) corresponds to a monadic twisted grouplike element \( p \). It remains to show that the condition \( L_p = R_p \mathcal{S}^2 \) translates to Figure~\ref{Fig:C-piv} for \( p = \Sigma(\mathfrak{p}) \), which is straightforward.
\end{proof}

\begin{theorem}
\label{Res:C-piv-1-1}
There is a bijection \( \mathsf{CPiv}(H) \cong \mathsf{Piv}(\catV_H) \) between \( C \)-pivotal elements of \( H \) and pivotal structures of \( \catV_H \).
\end{theorem}

\begin{proof}
Compose the two bijections \( \mathsf{CPiv}(H) \xrightarrow{\Sigma} \mathsf{MPiv}(H) \) in Lemma~\ref{Res:C-piv-M} and \( \mathsf{MPiv}(H) \xrightarrow{\mu (-)^{\sharp}} \mathsf{Piv}(\catV_H) \) in Proposition~\ref{Res:mnd-piv-1-1-can}.
\end{proof}

\subsection{Coend R-matrices}

We have seen that the results on coend \( R \)-matrices were already obtained by Bruguières and Virelizier in \cite{BV-12:double-Hopf-mnd}, see Section~\ref{Sec:C-qst-HA}. In particular, Theorem~\ref{Res:C-R-mat-1-1} is a consequence of the correspondence between monadic and coend \( R \)-matrices \cite[Section 8.6]{BV-12:double-Hopf-mnd}. 

\begin{remark}
\label{Rem:C-R-invertible}
By Remark~\ref{Rem:mnd-R-invertible}, when \( H \) is a Hopf algebra, any monadic \( R \)-matrix \( R \) for \( H \) is convolution invertible, and the axiom corresponding to the last diagram in Figure~\ref{Fig:C-R-mat} is automatically satisfied. Furthermore, the extended factorization property of \( C \) implies that there exists a coend \( R \)-matrix \( \overline{\mathfrak{R}} \) corresponding to the convolution inverse \( \overline{R} \) of \( R \). Using the first formula for \( \overline{R} \) in Figure~\ref{Fig:mnd-R-rev}, we can express \( \overline{\mathfrak{R}} \) in terms of \( \mathfrak{R} \) as in Figure~\ref{Fig:C-R-rev} below.
\end{remark}

\begin{figure}[ht]
    \centering
    \scalebox{.27}{%
    \def\svgwidth{\columnwidth}
    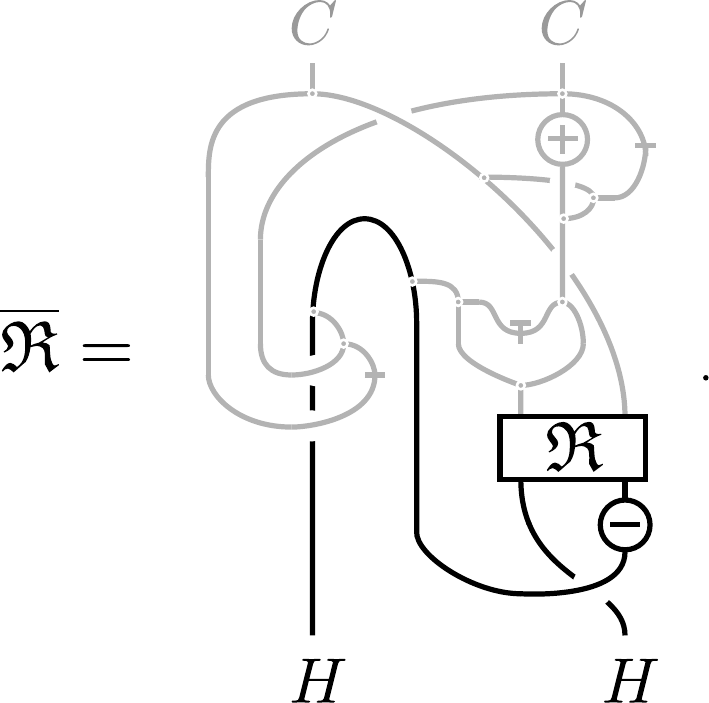
}
    \caption{The reverse coend \( R \)-matrix \( \overline{R} \) in terms of \( \mathfrak{R} \).}
    \label{Fig:C-R-rev}
\end{figure}

\subsection{Coend ribbon elements}
\label{Sec:C-rbn}

\begin{figure}[ht]
    \centering
    \scalebox{.4}{%
    \def\svgwidth{\columnwidth}
    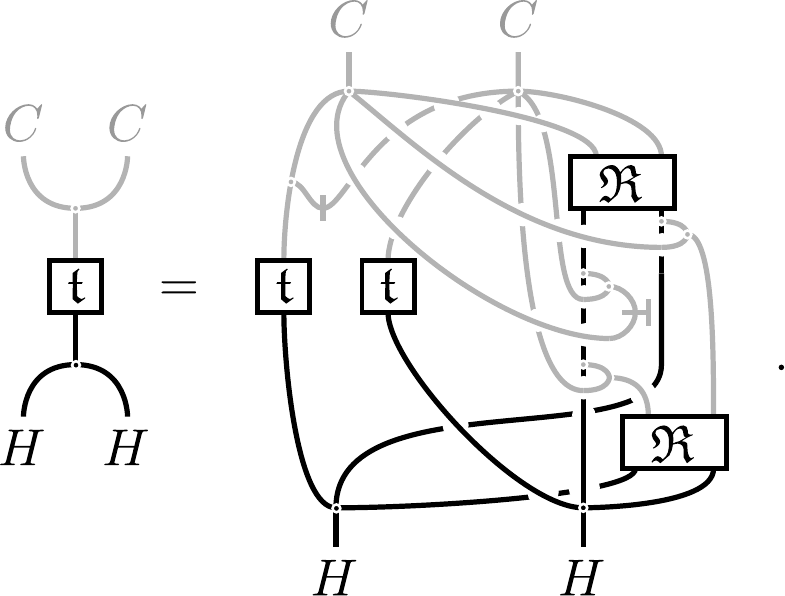
}
	\caption{The axiom (CT3) of a coend ribbon element \( \mathfrak{t}: C \to H \).}
	\label{Fig:C-twist}
\end{figure}

\begin{definition}
\label{Def:C-rbn}
	Let \( (H, \mathfrak{R}) \) be a quasitriangular bialgebra in a braided rigid monoidal category \( \catV \) admitting a coend \( C \). A \emph{\( C \)-balanced element} (or \emph{\( C \)-twist}) is a \( C \)-element \( \mathfrak{t}: C \to H \) satisfying the following properties:
\begin{enumerate}
\item[(CT1)] \( \mathfrak{t} \) is central (in the sense of Definition~\ref{Def:C-central-elt}).
\item[(CT2)] \( \epsilon_H \mathfrak{t} u_C = \var{id}_\ds1 \).
\item[(CT3)] \( \mathfrak{t} \) satisfies the axiom in Figure~\ref{Fig:C-twist}.
\end{enumerate}
Furthermore, when \( H \) is a Hopf algebra, \( \mathfrak{t} \) is \emph{ribbon} if it also satisfies
\begin{enumerate}
\item[(CT4)] \( S(\mathfrak{t}) = \mathfrak{t} \).
\end{enumerate}
The set of coend pivotal and ribbon elements are denoted by \( \mathsf{CBal}(H) \) and \( \mathsf{CRbn}(H) \), respectively.
\end{definition}

\begin{lemma}
\label{Res:C-rbn-M}
Under the bijection \( \Sigma \) in \textup{(\ref{Eqn:Sigma})}, the sets \( \mathsf{MBal}(H) \) and \( \mathsf{MRbn}(H) \) correspond to the sets \( \mathsf{CBal}(H) \) and \( \mathsf{CRbn}(H) \), respectively. 
\end{lemma}

\begin{figure}[ht]
    \centering
    \scalebox{1}{%
    \def\svgwidth{\columnwidth}
    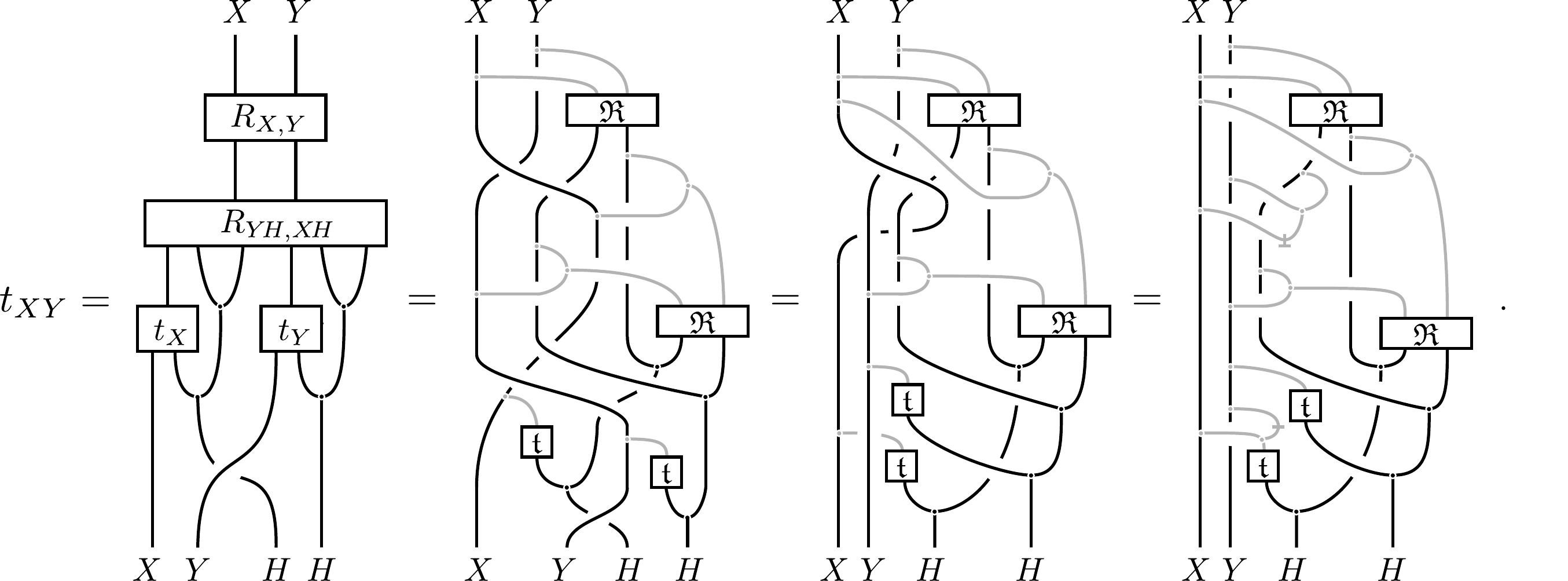
}
	\caption{Graphical proof of the axiom (CT3).}
	\label{Fig:C-twist-pf}
\end{figure}

\begin{proof}
Let \( t = \Sigma(\mathfrak{t}) \) be the corresponding monadic element. The equivalences of the pairs (CT1) and (MT1), (CT2) and (MT2), and (CT4) and (MT4) of axioms for \( \mathfrak{t} \)  and \( t \) are clear. It remains to show that (CT3) is equivalent to (MT3). The calculations are displayed in Figure~\ref{Fig:C-twist-pf}. In more detail, the first equality is obtained by applying the reverse braiding to the middle two strands in (MT3). For the second equality, the monadic elements \( R \) and \( t \) are expressed in terms of their coend counterparts \( \mathfrak{R} \) and \( \mathfrak{t} \). In the third equality, we use the naturality of the braiding to rearrange the morphisms. The last equality uses Lemma~\ref{Res:braid-switch}. 
\end{proof}

\begin{remark}
When \( H \) is a Hopf algebra, any coend balanced element \( \mathfrak{t} \) is convolution invertible, and the axiom (CT2) is automatically satisfied, see Remark~\ref{Rem:mnd-twist-invertible}.
\end{remark}

\begin{theorem}
\label{Res:C-rbn-1-1}
	Let \( \catV \) be a braided rigid monoidal category admitting a coend \( C \), and let \( (H, \mathfrak{R}) \) be a quasitriangular  bialgebra in \( \catV \).
\begin{enumerate}[label=\upshape(\alph*)]
\item Balanced structures \( \theta \) on \( \catV_H \) correspond bijectively to balanced \( C \)-elements \( \mathfrak{t}: C \to H \).
\item When \( H \) is a Hopf algebra, a balanced structure \( \theta \) is ribbon if and only if the corresponding \( C \)-balanced element \( \mathfrak{t} \) is ribbon.
\end{enumerate}	
\end{theorem}

\begin{proof}
The result follows from composing the bijections \( \Sigma \) and \( (-)^\sharp \) in Theorem~\ref{Res:mnd-rbn-1-1} and Lemma~\ref{Res:C-rbn-M}.
\end{proof}

\subsection{The coend Drinfeld element}

In this section, we assume that \( (H, \mathfrak{R}) \) is a quasitriangular Hopf algebra. Recall that in this case, we denote by \( \mu \) the Drinfeld morphism in \( \catV \), and the monadic Drinfeld element \( u \) corresponding to the Drinfeld morphism in \( \catV_H \) is defined as in Definition~\ref{Def:mnd-Drinfeld-elt-can}.

\begin{definition}
\label{Def:C-Drinfeld}
Let \( \mathfrak{u} \) be the coend element corresponding to the monadic Drinfeld element \( u \) under the bijection \( \Sigma \). We call \( \mathfrak{u} \) the \emph{coend Drinfeld element} of \( H \).
\end{definition}

\begin{figure}[ht]
    \centering
    \scalebox{.25}{%
    \def\svgwidth{\columnwidth}
    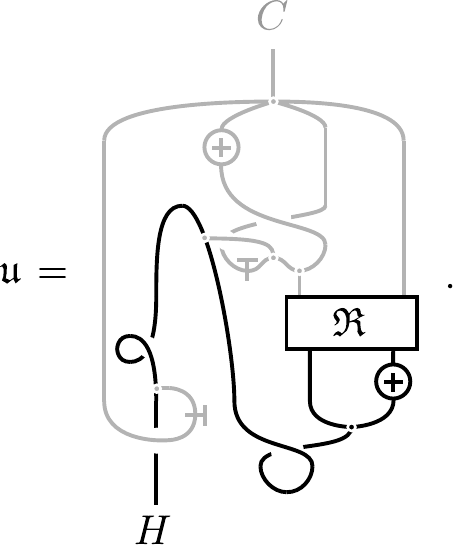
}
	\caption{The \( C \)-Drinfeld element \( \mathfrak{u} \).}
	\label{Fig:C-Drinfeld}
\end{figure}

\begin{remark}
The Drinfeld element \( \mathfrak{u} \) can be expressed using the \( R \)-matrix \( \mathfrak{R} \) as in Figure~\ref{Fig:C-Drinfeld}.
\end{remark}

\begin{theorem}
\label{Res:C-Drinfeld-bal-piv}
Let \( (H, \mathfrak{R}) \) be a quasitriangular Hopf algebra in a braided rigid monoidal category \( \catV \) admitting a coend \( C \). Via the coend Drinfeld element \( \mathfrak{u} \), there is a bijection \( \mathsf{CBal}(H) \cong \mathsf{CPiv}(H)  \), given by \( \mathfrak{t} \mapsto \mathfrak{t} \mathfrak{u} \). 
\end{theorem}

\begin{proof}
By design, the diagram
\begin{equation*}
\begin{tikzcd}[column sep=large]
	\mathsf{MBal}(H) & \mathsf{MPiv}(H) \\
	\mathsf{CBal}(H) & \mathsf{CPiv}(H)
	\arrow["u * (-)", "\cong" swap, from=1-1, to=1-2]
	\arrow["\Sigma", "\cong" swap, from=2-1, to=1-1]
	\arrow["\Sigma" swap, "\cong", from=2-2, to=1-2]
	\arrow["(-) \cdot \mathfrak{u}" swap, from=2-1, to=2-2]
\end{tikzcd}
\end{equation*}
commutes, and three out of four arrows are bijections by Proposition~\ref{Res:mnd-Drinfeld-bal-piv-can} and Lemmas~\ref{Res:C-piv-M}, \ref{Res:C-rbn-M}.
\end{proof}

Recall the definitions of \( c_0, q_\mu \) and \( c_\mu \) in Definition~\ref{Def:mnd-q-c-can}. We translate these monadic elements to their respective coend counterparts as follows.

\begin{definition}
\label{Def:C-q-c}
Define \( \mathfrak{q}_\mu = (\mathfrak{u} S(\bar{\mathfrak{u}}) \otimes \bar{\mathfrak{c}}_0) \Delta_C \) and \( \mathfrak{c}_\mu = (\mathfrak{u} S(\mathfrak{u})  \otimes \mathfrak{c}_0) \Delta_C \), where \( \bar{\mathfrak{u}} \) is the convolution inverse of \( \mathfrak{u} \), \( \mathfrak{c}_0: C \to \ds1 \) is the coend element of \( \ds1 \) corresponding to \( c_0 \) under the bijection \( \Sigma^{(1)}_\ds1 \) in (\ref{Eqn:coend-fp}), and \( \overline{\mathfrak{c}}_0 \) is the convolution inverse of \( \mathfrak{c}_0 \), see Definition~\ref{Def:mnd-q-c-can}.
\end{definition}

\begin{theorem}
\label{Res:C-Drinfeld-rbn-piv}
Under the correspondence \( \mathfrak{t} \mapsto \mathfrak{p} = \mathfrak{u} \mathfrak{t} \) between balanced and pivotal elements, \( \mathfrak{t} \) is ribbon if and only if \( \mathfrak{p}^2 = \mathfrak{q}_\mu \), or \( \mathfrak{t}^{-2} = \mathfrak{c}_\mu \). 
\end{theorem}

\begin{proof}
Clearly, \( \mathfrak{q}_\mu \) and \( \mathfrak{c}_\mu \) are defined to correspond to the monadic elements \( q_\mu \) and \( c_\mu \) in Definition~\ref{Def:mnd-q-c-can}. Consider the elements \( \mathfrak{t} \in \mathsf{CBal}(H) \), \( \mathfrak{p} = \mathfrak{u} \mathfrak{t} \), \( p = \Sigma(\mathfrak{p}) \), and \( t = \Sigma(\mathfrak{t}) \). Since the monadic elements \( p \) and \( t \) satisfy \( p^2 = q_\mu \) and \( t^{-2} = c_\mu \) by Theorem~\ref{Res:mnd-Drinfeld-rbn-piv-can}, and \( \Sigma \) is antimultiplicative, their coend counterparts \( \mathfrak{p} \) and \( \mathfrak{t} \) also satisfy the analogous properties.
\end{proof}

\begin{remark}
If \( \catV \) has a pivotal structure \( \psi \), then one can develop the theory of coend \( \psi \)-pivotal elements, the coend \( \psi \)-Drinfeld element \( \mathfrak{u}_\psi \), and the elements \( \mathfrak{q}_\psi = \mathfrak{u}_\psi S(\overline{\mathfrak{u}}_\psi) \) and \( \mathfrak{c}_\psi = \mathfrak{u}_\psi S(\mathfrak{u}_\psi) \) in a completely analogous way, from the theory of the corresponding monadic elements in Sections~\ref{Sec:mnd-piv} and~\ref{Sec:mnd-Drinfeld}. 
\end{remark}

\section{Concluding remarks}

\subsection{\texorpdfstring{Special cases of \( H \)}{Special cases of H}}
\label{Sec:H=1-or-C}

Two special cases of the Hopf algebra \( H \) are worth mentioning. When \( H = \ds1 \), the \( C \)-elements are morphisms \( C \to \ds1 \) and correspond to endomorphisms of the identity functor on \( \catV_H = \catV \). Note that the \( R \)-matrix corresponding to the default braiding \( \sigma \) on \( \catV \) is simply \( \epsilon_C \otimes \epsilon_C \otimes u_H \otimes u_H \). Our work therefore provides another description for the set of pivotal and ribbon structures on any braided rigid monoidal category admitting a coend, see Section~\ref{Sec:C-piv} and Section~\ref{Sec:C-rbn}. When \( H = C \), there is a canonical \( R \)-matrix \( C \otimes C \to C \otimes C \) which is given by \( u_C \otimes \epsilon_C \otimes \var{id}_C \), such that we have a braided isomorphism \( \catV_C \cong \catZ(\catV) \) \cite[Example~2.9]{BV-13:double-BHA}. Therefore, we can study structures on \( \catZ(\catV) \) in terms of endomorphisms of \( C \). 

\subsection{On factorizable Hopf algebras}
\label{Sec:factorizable-HA}

Let \( H \) be a Hopf algebra in a braided finite tensor category \( \catV \). By \cite{Shimizu-19:non-deg-BFTC}, the nondegeneracy of \( \catV_H \) is equivalent to the nondegeneracy of the canonical Hopf pairing of the coend of \( \catV_H \). We have seen that the coend of \( \catV_H \) is \( \dual H \otimes C \) from Example~\ref{Eg:coend-H-mod}. Its structure morphisms are displayed in \cite[Section~8.5]{BV-12:double-Hopf-mnd}.

\begin{definition}
Let \( (H, \mathfrak{R}) \) be a quasitriangular Hopf algebra in a braided finite tensor category \( \catV \). We say that \( (H, \mathfrak{R}) \) is \emph{factorizable} if the canonical Hopf pairing \( \omega \) for \( \dual H \otimes C \) is nondegenerate.
\end{definition}

\begin{remark}
We can now conclude that if \( H \) is a Hopf algebra in a braided finite tensor category \( \catV \), then \( \catV_H \) is a (non-semisimple) modular tensor category precisely when there is an \( R \)-matrix \( \mathfrak{R} \) and a ribbon element \( \mathfrak{t} \) for \( H \) such that \( (H, \mathfrak{R}) \) is factorizable.
\end{remark}

\subsection{On ribbon structures of the Drinfeld center}

A construction that frequently gives rise to a modular tensor category is the Drinfeld center \( \catZ(\catC) \) of a finite tensor category \( \catC \). By \cite[Theorem~1.1]{Shimizu-19:non-deg-BFTC} and \cite[Proposition~8.6.3]{EGNO-16:tensor-cat}, \( \catZ(\catC) \) is a nondegenerate braided finite tensor category, and hence \( \catZ(\catC) \) is modular if and only if \( \catZ(\catC) \) has a ribbon structure. As a result, the set of ribbon structures of \( \catZ(\catC) \) is a subject of great interest. For example, when \( \catC \) is semisimple and spherical, the pivotal structure on \( \catC \) gives rise to a ribbon structure on \( \catZ(\catC) \) \cite[Lemma~5.2]{TV-17:monoidal-cat-TFT}. This result extends to the nonsemisimple case in \cite[Theorem 5.11]{Shimizu-21:rbn-str-ZC} with the nonsemisimple notion of sphericality defined in \cite[Definition 3.5.2]{DSPS-20:dualizable-tensor-cat}.  

Note that when \( \catC = \catM_H \) for a finite dimensional Hopf algebra \( H \) over \( \dsk \), there is a braided isomorphism \( \catZ(\catC) \cong \catM_{D(H)} \), where \( D(H) \) is the Drinfeld double of \( H \). Thus, the modularity of \( \catZ(\catM_H) \) now translates to the set of ribbon elements of \( D(H) \). In \cite[Theorem 3]{KR-93:rbn-elt-DH}, Kauffman and Radford describes the ribbon elements of the Drinfeld double \( D(H) \) of a Hopf algebra \( H \) in terms of certain square roots of the distinguished grouplike elements of \( H \) and \( H^* \). The main result of \cite{Shimizu-21:rbn-str-ZC}, which describes the set of ribbon structures of \( \catZ(\catC) \) in terms of certain ``square roots'' of the distinguished invertible object of \( \catC \), can be seen as a tensor categorical generalization of this result.

In \cite{BV-12:double-Hopf-mnd}, Bruguières and Virelizier constructed the Drinfeld double \( D(H) \) of a Hopf algebra \( H \) in a braided rigid monoidal category \( \catV \), such that \( D(H) = H \otimes \dual H \otimes C \) is a quasitriangular Hopf algebra in \( \catV \) and there is an isomorphism \( \catZ(\catV_H) \cong \catV_{D(H)} \) of braided monoidal categories. For instance, when \( H = \ds1 \), \( D(H) \) is the quasitriangular Hopf algebra \( C \) described in Section~\ref{Sec:H=1-or-C} above. We now have two more descriptions for the set of ribbon structures of \( \catZ(\catC) \) when \( \catC = \catV_H \): as ribbon \( C \)-elements of \( D(H) \) (Theorem~\ref{Res:C-rbn-1-1}), and as pivotal \( C \)-elements of \( D(H) \) that are square roots of \( \mathfrak{q}_\mu \) (Theorem~\ref{Res:C-Drinfeld-rbn-piv}).

\section*{Acknowledgments}
The author would like to thank Chelsea Walton for valuable feedback and guidance throughout the project. The author also thanks Harshit Yadav for helpful correspondences and suggestions that improve several results in the paper.

\bibliography{references-tensor-categories}

\end{document}